\documentclass[leqno, 11pt]{article}

\usepackage{amsfonts}
\usepackage{amscd, amsxtra,float, tikz, subfigure, array, caption}
\usepackage{mathrsfs}
\usepackage{graphics}
\usepackage{subfigure}
\usepackage{amssymb}
\usepackage{amsmath, amsthm}
\usepackage{easybmat}
\usepackage{etex}
\usepackage{url}
\usepackage{texdraw}
\usepackage{epsfig}
\usepackage{tikz}
\usepackage{xcolor, soul}
\usepackage{mathrsfs}

\numberwithin{equation}{section}

\allowdisplaybreaks[4]

\begin{document}

\title{Spectral invariants of the Stokes problem}

\author{  Genqian Liu\thanks{E-mail: liugqz@bit.edu.cn} \\ 
\textit{\small School of Mathematics and Statistics, Beijing Institute of Technology, Beijing 100081, China }
}
\date{}
\maketitle

\begin{abstract}
For a bounded domain $\Omega\subset \mathbb{R}^n$ with smooth boundary, we  explicitly calculate the first two coefficients of the asymptotic expansion for the integral of the trace of the Stokes semigroup $e^{-t S}$ as $t\to 0^+$. These coefficients (i.e., spectral invariants) provide precise information for the volume  of the domain $\Omega$ and the surface area of the boundary $\partial \Omega$ by the spectrum of the  Stokes  problem. As an application, we show that an $n$-dimensional ball is uniquely determined by its Stokes spectrum among all Euclidean bounded domains with smooth boundary. 
\end{abstract}

\vspace{0.3cm}
{\bf Mathematics Subject Classification (2010):}~76D07; 35Q30; 35P20;  35S05   

\vspace{0.3cm} 
{\bf\em Keywords:}~ Stokes eigenvalues; Spectral invariants; Asymptotic expansion; Pseudodifferential operator

\section{ Introduction}

\vskip 0.45 true cm

Let $\Omega\subset \mathbb{R}^n$ ($n\ge 2$) be
 a bounded domain with smooth boundary $\partial \Omega$.
 We consider the following Stokes eigenvalue problem:
  \begin{eqnarray} \label{m1-1} \left\{ \begin{array}{ll}   -\mu \Delta {\mathbf{u}}+ \nabla p  = \lambda {\mathbf{u}}  \;\; &\mbox{in}\;\;  \Omega,\\
    \mbox{div}\; {\mathbf{u}}=0 \;\; & \mbox{in}\;\;  \Omega, \\
    {\mathbf{u}}=0 \;\;& \mbox{on}\;\;  \partial \Omega. \\
    \end{array}  \right.\end{eqnarray}
 Here $\mu$ is a positive constant (the kinematic coefficient of viscosity), the ``pressure'' term $p$ is not known a priori but is determined a posteriori from the solution itself.

It is well-known (see, p.$\,$457 of \cite{Pay}) that the problem (1.1) has nontrivial solutions ${\mathbf u}$ only for a discrete set of $\lambda =\lambda_k$, which are called
 Stokes eigenvalues. Let us enumerate the eigenvalues in
increasing order: $0< \lambda_1 \le \lambda_2 \le \cdots \le \lambda_k \le \cdots \to +\infty$,
where each eigenvalue is counted as many times as its multiplicity. The corresponding eigenvectors ${\mathbf{u}}_1, {\mathbf{u}}_2, \cdots, {\mathbf{u}}_k, \cdots$ form a complete orthonormal basis.

$\,$The eigenvalue problem (\ref{m1-1}) stems from $\,$the initial-boundary problem for the Stokes equations
 \begin{eqnarray} \label{m1-2}  \left\{ \begin{array}{ll}  \frac{\partial \mathbf{v}}{\partial t} -\mu \Delta \mathbf{v} +\nabla p =0 \;\; &\mbox{in}\;\;  (0,+\infty)\times \Omega,\\
    \mbox{div}\; \mathbf{v}=0 \;\; & \mbox{in}\;\; (0,+\infty) \times \Omega, \\
         \mathbf{v}=0 \;\;& \mbox{on}\;\; (0,+\infty)\times \partial \Omega, \\
     \mathbf{v}(0,x)= \mathbf{v}_0 \;\; & \mbox{on}\;\; \{0\}\times \Omega, \end{array}  \right.\end{eqnarray}
because  (\ref{m1-1}) can be immediately obtained by looking for the separated solutions in the Stokes equations (\ref{m1-2}). The Stokes equations (\ref{m1-2}) play an important role in fluid dynamics (see,  \cite{CM}, \cite{CF}, \cite{Lad}, \cite{SS} and \cite{Te}).
 More importantly,  the solutions of the Stokes equations provide a good approximation to the solutions of
  nonlinear Navier-Stokes equations:
   \begin{eqnarray} \label{m1-3}  \left\{ \begin{array}{ll}  \frac{\partial \mathbf{w}}{\partial t} -\mu \Delta \mathbf{w} + ({\mathbf{w}}\cdot \nabla){\mathbf{w}} +\nabla p =0 \;\; &\mbox{in}\;\;  (0,+\infty)\times \Omega,\\
    \mbox{div}\; \mathbf{w}=0 \;\; & \mbox{in}\;\; (0,+\infty) \times \Omega, \\
         \mathbf{w}=0 \;\;& \mbox{on}\;\; (0,+\infty)\times \partial \Omega, \\
     \mathbf{w}(0,x)= \mathbf{w}_0 \;\; & \mbox{on}\;\; \{0\}\times \Omega. \end{array}  \right.\end{eqnarray}
The Stokes eigenvalues are physical quantities. They just are the frequencies of the vibration of a Stokes (i.e., incompressible slow velocity, large viscosity, or small bodies) flow, which can be measured experimentally.
In theory of elasticity, the Stokes eigenvalue problem may be used to describe the vibration modes of an incompressible elastic body (see \cite{Pay}) but it also related with the classical buckling eigenvalue problem for a clamped plate. In 1986, Girault and Raviart \cite{GiR} (see also \cite{CL}) showed that all the Stokes eigenvalues coincide with all the buckling eigenvalues in the two-dimensional case ($\Lambda_k$ is said to be the $k$-th buckling eigenvalue for a clamped plate $\Omega\subset \mathbb{R}^n$  corresponding to the eigenfunction $\psi_k$ if and only if
\begin{eqnarray} \label{2020.12.8-1} \left\{ \begin{array} {ll} \mu\Delta^2 \psi_k +\Lambda_k \Delta \psi_k =0 \;\; &\mbox{in}\;\; \Omega,
\\  \psi_k =\frac{\partial \psi_k}{\partial {\boldsymbol{\nu}}} =0\;\; &\mbox{on}\;\; \partial \Omega,\end{array}\right.\end{eqnarray}
where ${\boldsymbol{\nu}}$ is the vector of the unit inner normal to $\partial \Omega$, see \cite{CH}, \cite{PS}, or (1.7) of \cite{Liu3}.

For the Stokes equations, one of the most important problems is to study the geometry of the flow region from the
physical quantities of the flow, because the geometric properties reflect the true behavior of the flow.
 An interesting question, which is similar to the famous Kac question for the Dirichlet-Laplacian
  (see \cite{Kac}, \cite{Lo} or \cite{We1} ), is: ``can one hear the shape of a Stokes flow  region by hearing the pitches (or frequencies) of its vibration?''
More precisely, for a Stokes flow region, one slaps the fluid of the flow region, and then listening to (i.e., measuring) the
vibration frequencies of the flow. Can one determine the shape (or the geometric quantities) of the flow region?

In order to explain the motivation and the importance of the above question, we briefly review the historical
background for the case of Dirichlet-Laplacian on domains. In 1910, H. A. Lorentz  conjectured that
for a two-dimensional domain $\Omega\subset \mathbb{R}^2$, the
 asymptotics of the counting function of the Dirichlet eigenvalues $\{\beta_k\}$ are given by:
  \begin{eqnarray} \label{116}  N_D(\tau)= \mbox{max}\, \{k\big|\beta_k \le \tau\}
  = \frac{|\Omega|}{2\pi}\tau +o(\tau)\quad \,\,
  \mbox{as}\;\; \tau\to \infty,\end{eqnarray}
  where $|\Omega|$ is the two-dimensional volume of $\Omega$ ($\beta_k$ is called the $k$-th eigenvalue of the Dirichlet-Laplacian if and only if
  $\Delta \phi_k+\beta_k\phi_k =0$ in $\Omega$ and $\phi_k=0$ on $\partial \Omega$).
This asymptotic in particular implies that $|\Omega|$ is a spectral invariant. Lorentz's conjecture was proved in 1913 by Hermann
Weyl  (see \cite{We1} and \cite{We2}).
With the Weyl formula as a starting point, Pleijel \cite{Ple1} in 1954 obtained more terms in the asymptotic expansion.
For a simply connected domain $\Omega$ in $\mathbb{R}^2$ with two-dimensional volume $|\Omega|$ and one-dimensional length $|\partial\Omega|$ of boundary $\partial\Omega$ he established the formula \begin{eqnarray} \label{0.2}  \sum_{k=1}^\infty e^{-\beta_k t} \sim \frac{|\Omega|}{2\pi t} -\frac{1}{4} \cdot
 \frac{|\partial \Omega|}{\sqrt{2\pi t}} +\frac{1}{6} \quad \; \mbox{as} \;\; t\to 0^+, \end{eqnarray}
 and in fact he showed that an additional term may be added to  the right side of (\ref{0.2}), one involving the curvature of the boundary of $\Omega$. Clearly, (\ref{0.2}) implies that one can get the area and the length by the spectrum of the Dirichlet problem in $\Omega$.
 In particular, by a Tauberian theorem the asymptotic formula for the first term on the right side of (\ref{0.2})
is equivalent to  Weyl's formula (\ref{116}).
 Kac \cite{Kac} used a combination of probability techniques and heat equation methods to establish the first two terms of (\ref{0.2}) for convex domains, and he obtained (\ref{0.2})
as a limiting case of convex polygonal domains. Kac also conjectured that for  multiply connected domains in $\mathbb{R}^2$  with $r$ holes,
the number $\frac{1}{6}$ in (\ref{0.2}) should be replaced by $\frac{1}{6}(1-r)$. McKean and Singer  in a celebrated paper \cite{MS} gave an
affirmative answer to the conjecture of Kac with respect to the third term for multiply connected domains in $n$-dimensional Riemannian manifold (with or without boundary). McKean and Singer \cite{MS} also obtained information about the curvature of the boundary of
 $\Omega$, which showed that the Euler characteristic
 $\chi(\Omega)$  is also a spectral invariant. Gilkey \cite{Gil} explicitly calculated the first four coefficients of the expansion of the trace of the heat kernel (see also \cite{vdB}). Furthermore, Branson and Gilkey in \cite{BGi} gave the first five coefficients of the asymptotic expansions
  for  the Dirichlet and Neumann boundary problems.

 Let us come back to the Stokes eigenvalue problem. We denote by $J$ and $J^1$ the closures in $[L^2(\Omega)]^n$  and the Sobolev space $[H^1(\Omega)]^n$ respectively of the set $\{ \mathbf{u}\in [C_0^\infty(\Omega)]^n\big| \mbox{div}\, \mathbf{u} =0\,\, \mbox{in}\,\, \Omega\}$
  of all smooth solenoidal vectors with compact supports in $\Omega$, where $[L^2(\Omega)]^n:=L^2(\Omega)\times  \cdots \times L^2(\Omega)$. Let $P_J$ be the orthogonal projection $[L^2(\Omega)]^n\to J$.
 We introduce the Stokes operator $S:= -\mu P_J \Delta$ with domain $\mathcal{D}(S)= J^1\cap [H^2(\Omega)]^n$, where
 $\Delta$  is the Laplace operator. It is easy to verify that
 the domain $\mathcal{D}(S)$ of the Stokes operator $S$ is dense in the Hilbert space $J$ with the inner product of $[L^2(\Omega)]^n$, and  the Stokes operator $S$ is an unbounded, self-adjoint, positive definite operator with respect to the $[L^2(\Omega)]^n$ inner product. Thus, the Stokes eigenvalue problem can be rewritten as
     \begin{eqnarray} \label{1-2} S{\mathbf{u}}_k =\lambda_k {\mathbf{u}}_k,\end{eqnarray} where ${\mathbf{u}}_k\in J^1\cap [H^2(\Omega)]^n$ are the orthogonal eigenvectors corresponding to the  Stokes eigenvalues $\lambda_k$.
 In 1986, Kozhevnikov \cite{Ko} gave an asymptotic formula  with sharp remainder estimate for the Stokes eigenvalues:
   \begin{eqnarray} \label{1-5} N_S(\tau)  =\frac{(n-1) \omega_n\,|\Omega|}{(2\pi)^{n}\mu^{n/2}} \tau^{n/2} +O(\tau^{(n-1)/2}) \quad \;\; \mbox{as}\;\; \tau\to +\infty,\end{eqnarray}
   where  $N_S(\tau)=\max\{k\big|\lambda_k\le \tau\}$ is the number of the Stokes eigenvalues less than or equal to $\tau$,
   $\,\omega_n$ denotes the volume of the unit ball in $\mathbb{R}^n$ and $|\Omega|$ denotes the volume of the domain $\Omega$.
  Weaker estimates of the remainder in the formula (\ref{1-5}) for $N_S(\tau)$ of the form $o(\tau^{n/2})$ and $O(\tau^{n/2}/\ln \tau)$ (the latter for $n=3$) were proved by Metivier \cite{Me} and Babenko \cite{Ba} respectively.
     The formula (\ref{1-5}) implies that one can ``hear'' (i.e., obtain) the volume of the domain $\Omega$ if one ``hears'' (i.e., knows) all the Stokes eigenvalues. For more  geometric quantities of the Stokes spectrum for $\Omega$, it has been a long-standing open problem.

In this paper, some surprising and interesting results are obtained $\,$by considering the Stokes operator semigroup $U(t) =e^{-t S}$ and by using some new methods of pseudodifferential  operators. The following theorem is the main result of this paper:

\vskip 0.25 true cm

  \noindent{\bf Theorem 1.1.} \ {\it Let  $\Omega\subset \mathbb{R}^n$ ($n\ge 2$) be
 a bounded domain with smooth boundary $\partial \Omega$, and let $0<  \lambda_1 \le \lambda_2 \le \cdots \le \lambda_k \le \cdots$ be the eigenvalues
 of the Stokes operator $S$.  Then
 \begin{eqnarray} \label{1-7} \quad \;  &&\sum_{k=1}^\infty e^{-\lambda_k t} = \int_\Omega \big(\mbox{Tr}\,(e^{-tS})\big) dx=\frac{(n-1)}{(4\pi \mu t)^{n/2}}|\Omega| - \frac{1}{4}\cdot\frac{(n-1)}{(4\pi \mu t)^{(n-1)/2}}|\partial \Omega|
 +O(t^{1-n/2})\;\; \mbox{as}\;\; t\to 0^+.\end{eqnarray}
 Here $|\Omega|$ denotes the $n$-dimensional volume of $\Omega$, and $|\partial \Omega|$ denotes the $(n-1)$-dimensional volume of  $\partial \Omega$.}

 \vskip 0.25 true cm

  Our result shows that not only the volume $|\Omega|$  but also the surface area $|\partial \Omega|$ can be known if  we know all Stokes eigenvalues. Roughly speaking,  one can ``hear'' the volume of the domain and the surface area of its boundary $\partial \Omega$ by ``hearing'' all the pitches of the vibration of a Stokes flow.

   The main ideas of this paper are as follows. It follows from Giga \cite{Gig},  Abe and Giga \cite{AG}, Solonnikov \cite{Sol} that the Stokes operator $S$ generates a $C_0$-analytic semigroup $e^{-t S}$ in the space $J$ or the space $C_{0,\sigma} (\Omega)$ (the $L^\infty(\Omega)$-closure of the space of all smooth solenoidal vector fields with compact supports in $\Omega$). Moreover, there exists an integral kernel (function matrix) $\mathbf{K}(t, x, y)$ such that
    \begin{eqnarray*}  e^{-t S}{\mathbf{f}}(x)= \int_{\Omega} {\mathbf{K}}(t,x, y){\mathbf{f}}(y) \, dy, \quad \forall \,\mathbf{f}\in J.\end{eqnarray*}
        If ${\mathbf{u}}_k$ is the normalized eigenvector of Stokes eigenvalue problem with eigenvalue $\lambda_k$, then the Stokes integral kernel ${\mathbf{K}}(t, x, y)$ is  given by \begin{eqnarray} \label{1-0a-1} \mathbf{K}(t,x,y) =\sum_{k=1}^\infty e^{-t \lambda_k} {\mathbf{u}}_k(x)\otimes {\mathbf{u}}_k(y).\end{eqnarray}
Thus the integral of the trace of ${\mathbf{K}}(t,x,y)$ is actually a spectral invariant: by  (\ref{1-0a-1}), we can get
\begin{eqnarray} \label{1-0a-2}\int_\Omega  \Big(\mbox{Tr}\big(\mathbf{K}(t,x,x)\big)\Big) dx=\sum_{k=1}^\infty e^{-t \lambda_k}.\end{eqnarray}
To further analyze the geometric contents of the spectrum, we calculate the same
integral of trace of the semigroup by another completely different way: we construct the Stokes semigroup $e^{-tS}$ by the Cauchy integral formula:
 \begin{eqnarray*} e^{-tS} =\frac{i}{2\pi} \int_{\mathcal{C}} e^{-t\lambda} (S-\lambda\, \mathbf{I})^{-1} d\lambda,\end{eqnarray*}
where $\mathcal{C}$ is a suitable curve in the complex plane in the positive direction around the spectrum of $S$. This leads us to discuss the resolvent operator $(S-\lambda\, \mathbf{I})^{-1}$ of $S$. The main difficulty (for proving theorem 1.1) is that unlike the Laplace operator, the Stokes operator is not a differential operator and it has not an explicit expression. Even more difficultly, neither McKean-Singer's classical method (see \cite{MS}) nor Gilkey-Seeley's calculus method of pseudodifferential operators (see \cite{See}, \cite{Gre} or \cite{BG}) can be applied and this may be a reason why the more (spectral) geometric quantities had been left open for a long time. Fortunately, we have a representation information for the inverse $S^{-1}$ of the Stokes operator $S$ by pseudodifferential (Green) operators that was given by Kozhevnikov \cite{Ko}.
Note that the  integral kernel of the Stokes semigroup is just the action of this semigroup on the Dirac Delta function $\delta(x-y)$. A key point in this paper is to find out the principal symbol and then estimate the lower order symbols of pseudodifferential (Green) operator. In \cite{Monv}  Boutet de Monvel proved that \begin{eqnarray} \label{2019.12.30-1} \begin{matrix} \left(\begin{matrix} M & K\\
 T & Q \end{matrix} \right)  :\;\; \, \begin{matrix} C^\infty(\bar \Omega, E) & {} &  C^\infty (\bar \Omega, E')\\
\oplus  & {\longrightarrow}  & \oplus \\
C^\infty (\partial \Omega, F) & {} & C^\infty (\partial \Omega, F)\end{matrix} \end{matrix}\end{eqnarray}
form an ``algebra''---i.e., the sum and the composition of two matrices such as (\ref{2019.12.30-1}) is another one if it is defined,  where $M$ is a Green operator (i.e., $M$ is the sum of a pseudodifferential operator $P$ (satisfying the transmission property) and a singular Green operator $G$); $K$ is a Poisson operator; $T$ is a trace operator; and $Q$ is a pseudodifferential operator on the boundary $\partial \Omega$. Such operators were posed in classical boundary problem. Applying this theory, we will prove that $\mbox{div}\, G_1 \mbox{grad}-I$ (which is closely related to the $S^{-1}$) is the sum of a pseudodifferential operator and a singular Green operator, each of which has order $-1$.
 Therefore, we extend $S^{-1}$  to all of $[L^2(\Omega)]^n$ that is denoted by $\mathbf{A}$ and given by the same expression as $S^{-1}$. Note that $[L^2(\Omega)]^n$ can be decomposed into the direct sum of three spaces $J$, $F$, $E$, and that $\mathbf{A}E=0$.
Thus, on $J\oplus F$ we will show $\mathbf{A}=\left(\begin{matrix} S^{-1} & 0\\
0 & \mathbf{A}_{FF}\end{matrix} \right)$. But $S^{-1}$ has non-negative eigenvalues, while $\mathbf{A}_{FF}$ has non-positive eigenvalues ($\mbox{dim}\, (\mbox{ker}\, \mathbf{A})= \mbox{dim}\, (\mbox{ker}\, \mathbf{A}_{FF})$ is a finite number).
In order to apply the asymptotic expansions of the heat kernels,  we will consider the Green operators
$({\mathbf{A}}^2+\mathbf{R})^{-1/2}$ and $(-{\tilde{A}}_{FF}-\tilde{R}_F)^{-1}$ as well as the strongly continuous semigroups
$e^{-t({\mathbf{A}}^2+\mathbf{R})^{-1/2}}$ and $e^{-t(-\tilde{A}_{FF}- \tilde{R}_F)^{-1}}$ because
 $(-\mathbf{A}_{FF}-\mathbf{{R}}_F)^{-1}$ can be reduced into an operator $(-{\tilde{A}}_{FF}-\tilde{R}_{F})^{-1}$ defined on $H^1_0(\Omega)$ and
  $(\mathbf{A}^2+ \mathbf{R})^{-1/2}=\left(\begin{matrix} S & 0\\
0 & (-\mathbf{A}_{FF}-\mathbf{{R}}_{F})^{-1}\end{matrix} \right)$ .
 Obviously, $\sum_{k=1}^\infty e^{-t\lambda_k}= \int_{\Omega}\! \big[\mbox{Tr} (e^{-t({\mathbf{A}}^2+\mathbf{R})^{-1/2}})\big]dx  - \int_{\Omega}\!\big[\mbox{Tr} (e^{-t(-\tilde{A}_{FF}- \tilde{R}_F)^{-1}})\big] dx.$
By calculating the  symbols of these two pseudodifferential (Green) operator semigroups, and then applying a technique of ``method of images'' which originates from McKean-Singer (see \cite{MS}), we  show (see Theorem 1.1) that the integral of the trace of the Stokes semigroup $e^{-tS}$ has an asymptotic expansion
\begin{eqnarray} \int_{\Omega}\Big( \mbox{Tr}\big({\mathbf{K}}(t,x,x)\big) \Big) dx \sim a_0t^{-n/2}+ a_1 t^{-(n-1)/2} +\cdots\quad \quad \mbox{as}\;\; t\to 0^+,\end{eqnarray}
   where $a_0= \frac{(n-1)|\Omega|}{(4\pi \mu )^{n/2}}$, $ \;a_1=-\frac{1}{4}\cdot\frac{(n-1)|\partial \Omega|}{(4\pi \mu )^{(n-1)/2}}$. More precisely,
 we first convert $e^{-t(\mathbf{A}^2 +\mathbf{R})^{-\frac{1}{2}}}$ and $e^{t(\tilde{A}_{FF}+ {\tilde{R}_F})}$ into
  $e^{-t(\mu^2 \Delta^2 \mathbf{I}-\mathbf{B}_1 -\mathbf{C}_1)^{1/2}}$ and $e^{-t(-\mu \Delta -B_2 -C_2)}$ by considering their symbols, respectively.
To obtain the coefficients $a_0$ and $a_1$, we will approximate the integral kernel near the boundary locally by the ``method of images.''  Let   $\mathcal{M}=\Omega \cup (\partial \Omega) \cup \Omega^*$ be the closed double of $\Omega$ and let us choose the geodesic normal coordinates in a collar neighborhood of the boundary $\partial \Omega$. Locally, the $\partial \Omega$ looks like the hyperplane $\{x=(x_1, \cdots, x_n)\big| x_n = 0\}$ in $\mathbb{R}^n$; letting $x\to x^*$ be
the reflection $(x_1, \cdots, x_{n-1}, x_n) \to (x_1, \cdots, x_{n-1},-x_n)$, the
${\mathbf{K}}_1^-(t, x, y) = {\mathbf{K}}_1(t, x, y) - {\mathbf{K}}_1(t, x, y^*)$ and ${{K}}_2^-(t, x, y) = {{K}}_2(t, x, y) - {{K}}_2(t, x, y^*)$ vanish on $\{x\in \mathbb{R}^n\big|x_n=0\}$, where ${\mathbf{K}}_1(t,x,y)$ and $K_2(t,x,y)$ are the integral kernels of $e^{-t(\mu^2 \Delta^2 \, \mathbf{I} -\mathbf{B}_1 -\mathbf{C}_1)^{1/2}}$
   and $e^{-t(-\mu \Delta  -B_2 -C_2)}$ defined on $(0,\infty)\times \mathcal{M}\times \mathcal{M}$, respectively.
    Next, for example, we consider $K(t,x,y)= e^{-t(-\mu \Delta -B_2 -C_2)} \delta(x-y)$, which is equal to $\frac{i}{2\pi}
  \int_{\mathcal{C}} (-\mu\Delta -\lambda)^{-1} \sum_{k=0}^\infty \big[(B_2+C_2) (-\mu\Delta-\lambda)^{-1} \big]^k \delta(x-y)$.
 By calculating the symbols of pseudodifferential operators and singular Green operators associated with $(-\mu\Delta -\lambda)^{-1} \big[(B_2+C_2) (-\mu\Delta-\lambda)^{-1} \big]^k$ according to the method provided by Boutet de Monvel in \cite{Monv}, we can show that for any $\Omega'\subset \Omega$,
    \begin{eqnarray*} \int_{\Omega'} \bigg\{ \frac{i}{2\pi} \int_{\mathcal{C}} e^{-t\lambda} (-\mu\Delta -\lambda)^{-1} \delta(x-x) \,d\lambda\bigg\} dx  =a_0 t^{-\frac{n}{2}} +O(t^{1-\frac{n}{2}})\quad \mbox{as}\, \; t\to 0^+;\\
    \int_{\Omega'} \bigg\{ \frac{i}{2\pi} \int_{\mathcal{C}}  e^{-t\lambda} (-\mu\Delta -\lambda)^{-1} \delta(x-\overset{\ast}{x}) \,d\lambda\bigg\} dx  =a_1 t^{-\frac{n-1}{2}} +O(t^{1-\frac{n}{2}})\quad \mbox{as}\, \; t\to 0^+\end{eqnarray*}
    and  \begin{eqnarray*} \int_{\Omega'} \bigg\{ \frac{i}{2\pi} \int_{\mathcal{C}} e^{-t\lambda} (-\mu\Delta -\lambda)^{-1} \sum_{k=1}^\infty \big[ (B_2+C_2) (-\mu \Delta -\lambda)^{-1} \big]^{k} \delta(x-x) \,d\lambda\bigg\} dx  = O(t^{1-\frac{n}{2}})\quad \mbox{as}\, \; t\to 0^+;\\
   \int_{\Omega'} \bigg\{ \frac{i}{2\pi} \int_{\mathcal{C}} e^{-t\lambda} (-\mu\Delta -\lambda)^{-1} \sum_{k=1}^\infty \big[ (B_2+C_2) (-\mu \Delta -\lambda)^{-1} \big]^{k} \delta(x-\overset{\ast}{x}) \,d\lambda\bigg\} dx  = O(t^{1-\frac{n}{2}})\quad \mbox{as}\, \; t\to 0^+.    \end{eqnarray*}
 A quite important fact is that
\begin{eqnarray*} && \frac{i}{2\pi} \int_{\mathcal{C}}e^{-t\lambda} \bigg[ \frac{1}{(2\pi)^n} \int_{\mathbb{R}^n} e^{i(x-y)\cdot \xi }\, l(x, \xi,\lambda) d\xi \\
&& \;\, \;\,\;\;\, + \frac{1}{(2\pi)^{n+1}} \int_{\mathbb{R}^{n-1}} d\xi' \int^+ d\xi_n  \int^+  r (x',\xi',\xi_n, \eta_n, \lambda)d\eta_n
\bigg] d\lambda\end{eqnarray*} just is the elementary solution of the Green operator semigroup $e^{-t(L+R)}$, where $(l (x, \xi,\lambda); r (x', \xi',\xi_n,\eta_n, \lambda))$ is  the  symbol of the resolvent (Green) operator $(L+P-\lambda)^{-1}$.  Thirdly, we apply the Gilkey-Seeley's calculus in the interior of $\Omega$ (see \cite{See}) to the principal symbols of the two kinds of pseudodifferential (Green) operator semigroups (i.e., $e^{-t(\mu^2 \Delta^2 \, \mathbf{I} -\mathbf{B}_1 -\mathbf{C}_1)^{1/2}}$ and $e^{-t(-\mu \Delta -B_2 -C_2)}$)
to get the  coefficient $a_0$, in which the corresponding terms for the lower symbols are put to $a_2$. By further estimating the integral of the traces of these two integral kernels and again applying the calculus method of Green operators near the boundary $\partial \Omega$, we finally obtain the coefficient $a_1$.

\vskip 0.10 true cm

\vskip 0.15 true cm

As an application of theorem 1.1, we can prove the following spectral rigidity result:
\vskip 0.25 true cm

  \noindent{\bf Corollary 1.2.} \ {\it Let $\Omega \subset \mathbb{R}^n$ be a bounded domain with smooth boundary $\partial \Omega$.
   Suppose that its Stokes spectrum is equal to that of $B_r$, a ball of radius $r$. Then $\Omega=B_r$. }

\vskip 0.23 true cm

  Corollary 1.2 also shows that a ball is uniquely determined by its Stokes spectrum among all Euclidean bounded domains with smooth boundary.

\vskip 1.49 true cm

\section{Pseudodifferential representation of the inverse $S^{-1}$}

\vskip 0.45 true cm

Denote $D^\alpha =D_1^{\alpha_1}\cdots D_n^{\alpha_n}$, $D_j=\frac{1}{i} \frac{\partial}{\partial x_j}$.
A pseudodifferential operator is an extension of the concept of differential operator. If $W$ is an open subset  of $\mathbb{R}^n$, we denote  by $S^m_{1,0}=S^m_{1,0} (W,
\mathbb{R}^n)$ the set of all $p\in C^\infty (W, \mathbb{R}^n)$ such
that for every compact set $O\subset W$ we have
 \begin{eqnarray} \label {-2.3} |D^\beta_x D^\gamma_\xi p(x,\xi)|\le C_{O,
 \gamma,\beta}(1+|\xi|)^{m-|\gamma|}, \quad \; x\in O,\,\, \xi\in \mathbb{R}^n\end{eqnarray}
 for all ${\gamma}, {\beta}\in {\mathbb{N}}^n$, where ${\mathbb{N}}^n$ is the set of ${\gamma}=(\gamma_1, \cdots, \gamma_n)$ with $\gamma_k=\mbox{integer}\ge 0$, $|{\gamma}|=\gamma_1+\cdots +\gamma_n$ and $|\xi|= \big( \sum_{j=1}^n \xi_j^2\big)^{1/2}$.
 The  elements  of $S^m_{1,0}$  are  called  symbols (or full symbols) of order $m$.
This class is defined by Kohn and Nirenberg in \cite{KN}. It is clear that $S^m_{1,0}$ is a
Fr\'{e}chet space with semi-norms given by the smallest constants
which can be used in (\ref{-2.3}) (i.e.,
\begin{eqnarray*} \|p\|_{O, \gamma,\beta}=
 \,\sup_{x\in O}\big|\left(D_x^\beta D_\xi^\gamma
 p(x, \xi)\right)(1+|\xi|)^{|\gamma|-m}\big|).\end{eqnarray*}
 For $p(x, \xi)\in S^m_{1,0}$, a pseudodifferential operator
in an open set $W\subset \mathbb{R}^n$ is defined by:
\begin{eqnarray} \label{-2.1}  p(x,D) f(x) = \frac{1}{(2\pi)^{n}} \int_{\mathbb{R}^n} p(x,\xi)
 e^{i  x\cdot\xi} \hat{f} (\xi)d\xi,\end{eqnarray}
and is denoted by $p(x,D)\in OPS^m_{1,0}$ (cf. \cite{Ho4}, \cite{Ta2}, \cite{Gr}).
Here $f\in C_0^\infty (W)$
and $\hat{f}(\xi) =\int_{\mathbb{R}^n}  f(x)e^{-ix\cdot\xi}\,d\xi$ is the Fourier transform of $f$ on $\mathbb{R}^n$.
   If there are smooth $p_{m-j} (x, \xi)$, homogeneous in $\xi$ of degree $m-j$ for $|\xi|\ge 1$, that is,
  $p_{m-j} (x, r\xi) =r^{m-j} p_{m-j} (x, \xi)$ for $r, |\xi|\ge 1$, and if
 \begin{eqnarray} \label{-2.2} p(x, \xi) \sim \sum_{j\ge 0} p_{m-j} (x, \xi)\end{eqnarray}
 in the sense that
 \begin{eqnarray}  p(x, \xi)-\sum_{j=0}^k p_{m-j} (x, \xi) \in S^{m-k-1}_{1, 0}, \end{eqnarray}
 for all $k$, then we say $p(x,\xi) \in S_{cl}^m$, or just $p(x, \xi)\in S^m$. We call
 $p_m(x, \xi)$ the principal symbol of $p(x, D)$, sometimes is denoted by $\sigma(p(x,D))$. We also write $OPS^{-\infty}= \bigcap_{m} OPS^m_{1,0}$.

\vskip 0.10 true cm

 It is well-known (see p.$\,$3 of  \cite{Ta2} or Theorem 18.1.6 of \cite{Ho3}) that if $p(x, \xi)\in S^m_{1,0}$, then $p(x, D):\mathscr{S}(\mathbb{R}^n)\to \mathscr{S} (\mathbb{R}^n)$,
 where $\mathscr{S}(\mathbb{R}^n)= \{f\in C^\infty (\mathbb{R}^n)\big| x^\beta D^\alpha f\in L^\infty (\mathbb{R}^n)\,\,\mbox{for all}\,\, \alpha,\beta\ge 0\}$
 is the space of all rapidly decreasing functions. Furthermore, according to Proposition 5.5 on p.$\,$20 of \cite{Ta2} we see that if $p(x, \xi)\in S^m_{1,0}$,  then \begin{eqnarray} \label{2020.11.29-1}  p(x, D): H^s(\Omega) \to H^{s-m} (\Omega)\end{eqnarray}  for any $s\in \mathbb{R}^1$, where $\Omega\subset \mathbb{R}^n$ is a bounded domain with smooth boundary.
  Let $p_j(x, D)\in OPS^{m_j}_{1,0}$, $j=1,2$, be the pseudodifferential operators of order $m_j$ with symbols $p_j(x, \xi)$. Then $ p_1(x, D) \,p_2(x,D)= q(x, D)$ is a pseudodifferential operator of order $m_1+m_2$ with symbol (see (3.23) of Chapter 7 in \cite{Ta2})
\begin{eqnarray}\label{17-9-12.1}  q(x,\xi)= \sum_{\alpha\ge 0} \frac{i^{|\alpha|}}{\alpha !} D^\alpha_\xi p_1(x,\xi) D^\alpha_x p_2 (x,\xi).\end{eqnarray}

   An operator $P$ with symbol $p(x,D)$ is said to be an elliptic pseudodifferential operator of order $m$ if
  for every compact $O\subset \Omega$ there exists a  positive constant $c=c(O)$ such that \begin{eqnarray*} |p(x, \xi)|\ge c|\xi|^m, \,x\in O,\, |\xi|\ge 1.\end{eqnarray*}
    If $q(x, D)\in OPS^{-m}_{1,0}$ is a pseudodifferential operator of order $-m$ such that
\begin{eqnarray*} q(x,D)p(x,D)=I\;\; \mbox{mod}\;\; OPS^{-\infty},\\
 p(x,D)q(x,D)=I\;\; \mbox{mod}\;\; OPS^{-\infty}, \end{eqnarray*}
 then we say that $q(x,D)$ is a (two-sided) parametrix for $p(x,D)$.
 Furthermore, if $P$ is a non-negative elliptic pseudodifferential operator of order $m$, then the spectrum of $P$ lies in a right half-plane and has a finite lower bound $\tau(P) =
\inf\{\mbox{Re}\, \lambda\big| \lambda\in \sigma(P)\}$. We can modify the principal symbol $p_m(x, \xi)$ for small $\xi$ such that
$p_m(x,\xi)$  has a positive lower bound throughout and lies in $\{\lambda=re^{i\theta} \big| r>0, |\theta|\le \theta_0\}$, where $\theta_0\in (0,\frac{\pi}{2})$.
According to \cite{Gr}, the resolvent $(P-\lambda)^{-1}$ exists and is holomorphic in $\lambda$ on a neighborhood of a set
\begin{eqnarray*}  W_{r_0,\epsilon} =\{\lambda\in {\mathbb{C}} \big| |\lambda|\ge r_0, \,\mbox{arg}\, \lambda \in [\theta_0+\epsilon, 2\pi -\theta_0-\epsilon], \,\mbox{Re}\, \lambda \le \tau (P)-\epsilon\}\end{eqnarray*}
(with $\epsilon>0$). There exists a parametrix $Q'_\lambda$ on a neighborhood of a possibly larger set
(with $\delta>0,\epsilon>0$)
\begin{eqnarray*} V_{\delta,\epsilon} =\{ \lambda \in {\mathbb{C}} \big| |\lambda| \ge \delta \;\;\mbox{or arg}\, \lambda \in [\theta_0+\epsilon, 2\pi-\theta_0-\epsilon]\}\end{eqnarray*}
such that this parametrix coincides with $(P-\lambda)^{-1}$ on the intersection. Its symbol $q(x,\xi, \lambda)$
in local coordinates is holomorphic in $\lambda$ there and has the form (cf. Section 3.3 of \cite{Gr})
\begin{eqnarray} \label {-2.6} q(x, \xi,\lambda) \sim \sum_{l\ge 0} q_{-m-l} (x, \xi, \lambda),\end{eqnarray}
where \begin{eqnarray} \label{-2.7} & q_{-m}= (p_m(x, \xi) -\lambda)^{-1}, \quad \;  q_{-m-1} = \vartheta_{1,1}(x, \xi) q^2_{-m},\\
& \, \cdots, \,  q_{-m-l} = \sum_{k=1}^{2l} \vartheta_{l,k} (x, \xi) q^{k+1}_{-m}, \cdots \nonumber\end{eqnarray}
with symbols $\vartheta_{l,k}$ independent of $\lambda$  and homogeneous of degree $mk-l$ in $\xi$ for $|\xi|\ge1$.
 The semigroup $e^{-tP}$ can be defined from $P$ by the Cauchy integral formula (see p.$\,$4 of \cite{GG}):
 \begin{eqnarray*} e^{-tP} =\frac{i}{2\pi} \int_{\mathcal{C}} e^{-t\lambda} (P-\lambda)^{-1} d\lambda,\end{eqnarray*}
where $\mathcal{C}$ is a suitable curve in the complex plane in the positive direction around the spectrum of $P$.
From (\ref{-2.6}) and the above formula, we  get the symbol $\frac{i}{2\pi} \int_{\mathcal{C}} e^{-t\lambda} \big[
\sum_{l\ge 0} q_{-m-l} (x, \xi, \lambda)\big]\, d\lambda$ of the semigroup $e^{-tP}$, and furthermore we can  obtain the trace of $e^{-tP}$.
\vskip 0.12 true cm

 The following decomposition lemma is well-known:

 \vskip 0.2 true cm

 \noindent{\bf Lemma 2.1} (see p.$\,$37 of \cite{CM}). \ {\it Any vector field ${\mathbf{u}}$ on $\Omega$ can  be  uniquely decomposed in the form}:
 \begin{eqnarray} \label{m2-1} {\mathbf{u}} ={\mathbf{w}}+\mbox{grad}\, p,\end{eqnarray}
 {\it where ${\mathbf{w}}$ satisfies}
 \begin{eqnarray} \label{m2-2} \left\{ \begin{array}{ll} \mbox{div}\, {\mathbf{w}}=0 \;\; \mbox{in}\;\; \Omega,\\
 {\mathbf{w}}\cdot {\boldsymbol{\nu}}=0 \;\; \mbox{on}\;\; \partial \Omega,\end{array}\right. \end{eqnarray}
 {\it and ${\boldsymbol{\nu}}$ is  the unit inward normal to $\partial \Omega$.}

 \vskip 0.2 true cm

  Recall that $J$  is the closure in $[L^2(\Omega)]^n$  of the set of all smooth solenoidal vectors with compact supports in $\Omega$.  As is known (see \cite{Te}),  the space $J$ of the vector-valued functions $\mathbf{u}=(u_1, \cdots, u_n)$ can be rewritten as
  \begin{eqnarray} \label{2-1} J=\{\mathbf{u}\in [L^2(\Omega)]^n\big| \mbox{div}\, \mathbf{u}=0,\, \gamma_{{}_{\boldsymbol{\nu}}} \mathbf{u}=(\mathbf{u}\cdot  {\boldsymbol{\nu}})\big|_{\partial \Omega}=0\}. \end{eqnarray}
   It follows from  \cite{Te} that  the operator $\gamma_{\boldsymbol{\nu}} \!\mathbf{u} \equiv (\mathbf{u}\,\cdot \,{\boldsymbol{\nu}}) \big|_{\partial \Omega}$ continuously maps the Hilbert space
 $\{ \mathbf{u}\in [L^2(\Omega)]^n\big| \mbox{div}\, \mathbf{u} \in L^2(\Omega)\}$  with inner product $\langle\langle \mathbf{u}, \mathbf{v}\rangle\rangle \equiv
 \langle\mathbf{u},\mathbf{v}\rangle+ \langle\mbox{div}\, \mathbf{u}, \mbox{div}\, \mathbf{v}\rangle$ into the space $H^{-1/2}(\partial \Omega) $. Here
 $\langle\cdot, \cdot\rangle$ is the inner product in $[L^2(\Omega)]^n$ (or $[L^2(\Omega)]$). Note also that (see Chapter 1, \cite{Te})
 $$J^1=\big\{ \mathbf{u} \in [H_0^1(\Omega)]^n \big| \mbox{div}\, \mathbf{u} =0\big\}.$$
  It follows from Lemma 2.1 that  $[L^2(\Omega)]^n$ is the orthogonal sum of $J$ and the space $\{\mathbf{u}\in [L^2(\Omega)]^n\big| \mathbf{u}=\mbox{grad}\, p,\, p\in H^1(\Omega)\}$.

 We also introduce the following spaces
 \begin{eqnarray} \label{2-2}  && F:=  \{ \mathbf{u}\in [L^2(\Omega)]^n \big|  \mathbf{u}=\mbox{grad}\, p, \  p\in H^1_0(\Omega)\}, \\
\label{2-3} && E := \{ \mathbf{u}\in [L^2(\Omega)]^n \big|  \mathbf{u}=\mbox{grad}\, p, \  p\in H^1(\Omega), \, \Delta p=0\}.\end{eqnarray}
 Then, the following orthogonal Weyl-Sobolev decomposition holds:

  \vskip 0.25 true cm

  \noindent{\bf Lemma 2.2} (see Chapter 1 of \cite{Te}). \ {\it The space $[L^2(\Omega)]^n$ can be orthogonally decomposed into the  sum of
  $J$, $F$ and $E$:  \begin{eqnarray} \label {2-4} [L^2(\Omega)]^n = J \oplus F\oplus E,\end{eqnarray}
i.e., any vector-valued function $\mathbf{f}\in [L^2(\Omega)]^n$ here admits a unique orthogonal decomposition
\begin{eqnarray} \label{2-5} \mathbf{f}= {\mathbf{f}}_J +{\mathbf{f}}_F + {\mathbf{f}}_E \quad \mbox{with}\;\; {\mathbf{f}}_J\in J,\, \; {\mathbf{f}}_F \in F,\; \, {\mathbf{f}}_E\in E.\end{eqnarray}  }

 \vskip 0.20 true cm

Denote by  $\gamma_k \mathbf{u}$ the boundary value on $\partial \Omega$
 of the derivative \begin{eqnarray}\label{a-2.1}  D^k_{\boldsymbol{\nu}} := (\frac{1}{i})^k \frac{\partial^k}{\partial x_{\boldsymbol{\nu}}^k}\end{eqnarray} of $\mathbf{u}$ in the direction of the inner normal ${\boldsymbol{\nu}}$ to the boundary $\partial \Omega$.
  We denote by $G_1$ and $G_2$ the operators solving the Dirichlet problems for the Poisson
and Laplace equations (see \cite{LM}):
 \begin{eqnarray} \label{2-6} & G_1: f\to v, \; \mbox{where}\;\; \Delta v= f \;\; \mbox{in}\;\; \Omega, \;\; v= 0\;\;\mbox{on}\;\; \partial \Omega, \quad \;\, G_1: L^2(\Omega)\to H^2(\Omega),\\
 \label{2-7} & \left. \begin{array}{ll}  G_2: g \to w, \;\;\mbox{where}\;\; \Delta w= 0\;\; \mbox{in}\;\; \Omega, \;\; w=g \;\;\mbox{on}\;\; \partial \Omega, \\
  G_2: H^s (\partial \Omega) \to H^{s+1/2}(\Omega) \quad \;\; (s\ge 1).\end{array}\right.\end{eqnarray}
It follows from Chapter 1 of \cite{Te} that the projection $P_F$ onto the subspace $F$ of (\ref{2-2})
has the form:
\begin{eqnarray} \label{2-8} P_F= \mbox{grad}\, G_1 \, \mbox{div}.\end{eqnarray}

  The following result due to Kozhevnikov (see, \cite{Ko}), shows that the inverse $S^{-1}$ of the Stokes operator has an explicit expression in terms of operators $G_1$ and $G_2$ solving the classical Dirichlet problems for the Poisson and Laplace equations:

   \vskip 0.2 true cm

\noindent  {\bf Lemma 2.3.} \   {\it The Stokes operator $S$ is continuously invertible in the space $J$, and there exists a pseudodifferential operator $K_{-1}\in OPS^{-1}_{1,0}$ defined on $\partial \Omega$  such that}
\begin{eqnarray} \label{2-12} S^{-1} \mathbf{f} = -(1/\mu)G_1 [\mathbf{I}-2\, \mbox{grad}\,(I+G_2 K_{-1}\gamma_0)\mbox{div}\, G_1]\mathbf{f},\quad \; \, \forall \,\mathbf{f}\in J, \end{eqnarray}
{\it where the operators $G_1, G_2$ and $\gamma_0$ are defined in (\ref{2-6}), (\ref{2-7}) and (\ref{a-2.1}),
 and $\mathbf{I}$ is an identity matrix.}

\vskip 0.2 true cm

Let us point out that the above lemma was obtained by eliminating the pressure $p$ from the Stokes equations $-\mu \Delta \mathbf{u} +\mbox{grad}\, p =\mathbf{f}$ in terms of the vector-valued function ${\mathbf{f}}$ and by calculating the principal symbol of $S^{-1}$ (see \S4 of \cite{Ko}).

\vskip 1.19 true cm

\section{Transmission property, Poisson operators, trace operators and  Green operators}

\vskip 0.45 true cm

In order to study the boundary value problems, one needs to separately consider variable $\xi_n$ in the definition of a symbol. We first introduce some complex valued functions on the real line.

 Let $H$ be the  space of all complex valued functions $f(s)$ on the real line,
which are $C^\infty$ and have a regular pole at infinity, i.e.
$f$ is $C^\infty$ and has an asymptotic expansion
\begin{eqnarray} \label{2019.12.30-3}  f(s)\sim \sum\limits_{k\ge -m} a_k s^{-k} \quad \, (s\to \infty)\end{eqnarray}
and this expansion still holds after any number of differentiations (see \cite{Monv}).
 Let $H^+$ be the subspace consisting of those functions  $f\in H$ which can be extended
analytically in the lower complex half plane $\mbox{Im}\, s \le 0$, and vanish at infinity (for such functions, the asymptotic expansion (\ref{2019.12.30-3}) holds when $s\to \infty$, $\mbox{Im}\, s \le 0$, and $m=-1$, see \cite{Monv}).
 Let $H^-$ be the supplementary of $H^+$  in $H$ consisting of those functions which can
be extended analytically in the upper half plane $\mbox{Im}\, s\ge 0$: for such functions, the asymptotic
expansion (\ref{2019.12.30-3}) holds when $s\to \infty$, $\mbox{Im}\, s\ge 0$.
 We  will denote by $h^+$ (respectively, $h^-$) the projection on $H^+$ parallel to $H^-$ (respectively, $h^-=1-h^+$).
Thus if $f$ is analytic on the real line, meromorphic at infinity, we have (see \cite{Monv})
\begin{eqnarray}\label{2020.11.12-1} \begin{aligned} &h^+ f(s) =-\frac{1}{2\pi i} \int_{\gamma} \frac{f(\tau)} {\tau-s} d\tau \quad \, (\mbox{if}\;\, \mbox{Im}\, s\le 0)\\
&h^- f(s) =\frac{1}{2\pi i} \int_{\gamma} \frac{f(\tau)}{\tau-s} d\tau \quad \, (\mbox{if}\;\, \mbox{Im}\, s>0),\end{aligned}\end{eqnarray}
where $\gamma$ is a large circle in the upper half plane $\mbox{Im}\, \tau>0$, oriented in the usual way ($\gamma$  is
required to contain $s$ in its interior for the second formula).
If $f$  vanishes at infinity, we also have (see \cite{Monv})
\begin{eqnarray}\label{2020.11.12-2} \begin{aligned} &h^+ f(s) =\lim_{\epsilon\to 0}  -\frac{1}{2\pi i} \int_{-\infty}^{+\infty}  \frac{f(\tau)} {\tau-s+i\epsilon}  d\tau \\
&h^- f(s) =\lim_{\epsilon\to 0}  \frac{1}{2\pi i} \int_{-\infty}^{+\infty}  \frac{f(\tau)}{\tau-s-i\epsilon}  d\tau\end{aligned}\end{eqnarray}
If $f$ is analytic on the real line, and meromorphic at infinity, we set
\begin{eqnarray} \label{17-11-11-0} \int^{+} f= \int^{+} f(s) ds = \int_{\gamma} f(\tau)d\tau. \end{eqnarray}
This operator extends continuously to $H$.

\vskip 0.25 true cm

\noindent{\bf Remark 3.1}. \   {\it It is well-known (see p.$\,$15 of \cite{Monv} or p.$\,$154-155 of \cite{Ah}) that $\int^{+}\!f= \int_{-\infty}^{+\infty} f(s)ds$ is just the ordinary
integral if $f$ is integrable (i.e., vanishes to the second order at infinity). Also, $\int_{\gamma}\! e^{i\alpha \tau} f(\tau) d\tau = \int_{-\infty}^{+\infty}  e^{i\alpha s} f(s) ds$ is the ordinary integral if $f$
vanish to the first order at infinity}.

\vskip 0.2 true cm

Let $f\in H$. Then $f$ has a unique expansion (\cite{Monv})
\begin{eqnarray} \label{2019.12.30-5} f(s)= \sum\limits_{l=1}^m \alpha_l s^l + \sum\limits_{k=-\infty}^{+\infty} a_k \left( \frac{1-is}{1+is}\right)^k,\end{eqnarray}
where the coefficients $a_k$ form a rapidly decreasing sequence.  If $f$ vanishes at infinity,
 then $f$ belongs to
$H^+$ (respectively, $H^-$) if and only if $\alpha_l=0$ and $a_k= 0$  when $k < 0$ (respectively, $k\ge 0$).
In (\ref{2019.12.30-5}), one can replace $\big(\frac{1-is}{1+is}\big)^k$ (respectively,
$\frac{(1-is)^k}{(1+is)^{k+1}}$) by  $\big(\frac{\rho-is}{\rho+is}\big)^k$ (respectively,
$\frac{(\rho-is)^k}{(\rho+is)^{k+1}}$), where $\rho$  is any positive number. Note that $\frac{(\rho-is)^k}{(\rho+is)^{k+1}}\;\, (k\ge 0)$ is the Fourier transform of the product of a Laguerre polynomial by an exponential:
\begin{eqnarray*} \phi(x) = \left\{\!\! \begin{array}{ll} \big(\rho-\frac{d}{dx}\big)^k \big(\frac{x^k}{k!}e^{-\rho x}\big)\quad &\mbox{if}\;\, x>0\\
0 \quad &\mbox{if}\;\, x<0.\end{array}\right.\end{eqnarray*}

\vskip 0.29 true cm

   \noindent  {\bf Lemma 3.2} (see (1.9) of \cite{Monv}). \   {\it $H^+$ is the space of Fourier transforms of functions $\psi(x)$  which vanish for $x <0$ and
are $C^\infty({\bar {\mathbb{R}}}_+)$, rapidly decreasing at infinity for $x > 0$ (i.e. every derivative tends to zero at infinity, faster than any power of $x$, and has a limit when $x\to +0$)}.

\vskip 0.22 true cm

 Let $p(x, D)$ be a pseudodifferential operator defined in a neighborhood of the closed half space $\mathbb{\bar R}^n_+$. $p(x, D)$ is said to has the {\it transmission property with respect to the boundary} ${\mathbb{R}}^{n-1}$ (see p.$\,$21 of \cite{Monv}) if every derivative of its symbol admits the following  series expansion  along
the boundary:
\begin{eqnarray*} (\partial/\partial x)^\alpha p(x,\xi) = \sum\limits_{s=0}^d \alpha_s (x, \xi') \xi_n^s + \sum\limits_{m=0}^\infty
a_m(x, \xi') \frac{(\langle \xi'\rangle -i \xi_n)^m}{(\langle \xi'\rangle +i \xi_n )^{m+1}}  \quad \mbox{if} \;\, x_n=0,\end{eqnarray*}
where $\alpha_s \in S^{d-s}_{1,0}$, $a_m$ is a rapidly decreasing sequence in $S^{d+1}_{1,0}$ and $\langle \xi'\rangle =(1+|\xi'|^2)^{\frac{1}{2}}$.

\vskip 0.2 true cm

 Let $k(x', \xi)$ be a $C^\infty$ function on $\mathbb{R}^{n-1}\times \mathbb{R}^n$, admitting a series expansion (see p.$\,$26 of \cite{Monv}):
\begin{eqnarray} \label{17-11-11-2} k(x', \xi)=\sum_{m=0}^\infty a_m (x', \xi') \frac{(\langle \xi'\rangle -i\xi_n)^m}{ (\langle \xi'\rangle +i \xi_n)^{m+1}},\end{eqnarray}
where $a_m(x',\xi')$ is a rapidly decreasing sequence in $S_{1,0}^d$.
 The {\it Poisson operator} $K$ of order $d$ with symbol $k(x', \xi)$ is the operator $K: C_0^\infty (\mathbb{R}^{n-1})\to C^\infty ( \mathbb{\bar R}^n_+)$ defined by:
  \begin{eqnarray} \label{17-11-11-3} (Kf)(x) = \frac{1}{(2\pi)^{n}} \int_{-\infty}^{+\infty} d\xi_n \int_{\mathbb{R}^{n-1}} e^{i x\cdot \xi} k(x', \xi) \hat{f} (\xi') d\xi'.\end{eqnarray}

 \vskip 0.2 true cm

 \noindent  {\bf Example 3.3.} {\it A typical example of the Poisson operator is the operator that solves the Dirichlet problem
in the special half space $(\mathbb{R}_+^n, g)$ (see, p.$\,$2504 of \cite{Liu1}):
\begin{eqnarray} \label{17-11-12-1} \left\{ \begin{array}{ll} \sum_{j,k=1}^{n} g^{jk} \frac{\partial^2 u}{\partial x_j \partial x_k} =0,\quad &\mbox{in}\;\;
\mathbb{R}^n_+, \\
u(x',0)=f(x'), \quad &\mbox{on}\,\, \partial \mathbb{R}_+^n,\end{array}\right. \end{eqnarray}
where $(g^{jk})$ is the inverse of $(g_{jk})$, and
\begin{gather*} \label{3---1}  \left(g_{jk}\right)=\begin{pmatrix} g_{11} & \cdots  & g_{1,n-1} &0\\
           \vdots & \ddots & \vdots & \vdots \\
      g_{n-1,1}  & \cdots & g_{n-1,n-1}  & 0\\
      0 & \cdots   & 0 &  1 \end{pmatrix},\end{gather*}
    is a positive definite, real symmetric
  $n\times n$  constant matrix.
      If $f\in C_0^\infty (\mathbb{R}^{n-1})$, the unique bounded solution of (\ref{17-11-12-1}) is:
\begin{eqnarray} \label{17-11-12-3} \quad\quad\;  u(x',x_n)\! \!&\!\!=\!\!&\!\!\frac{\Gamma (n/2)}{\pi^{n/2}} \int_{\mathbb{R}^{n-1}} \frac{x_n}{\big(x_n^2+ \sum_{j,k=1}^{n-1} g_{jk} (x_j-y_j)(x_k-y_k)\big)^{n/2}}{f} (y')dy'\\
\! \!&\!\!=\!\!&\!\!\frac{1}{(2\pi)^{n-1}} \int_{\mathbb{R}^{n-1}} e^{ix'\cdot \xi'} e^{-x_n(\sum_{j=1}^{n-1} g^{jk} \xi_j\xi_k)^{\frac{1}{2}}} \hat{f} (\xi') d\xi'\nonumber\\
\! \!&\!\!=\!\!&\!\!\frac{1}{(2\pi)^{n-1}} \int_{\mathbb{R}^{n-1}} e^{ix'\cdot \xi'}\bigg(\frac{1}{2\pi i}\int^+ \frac{
e^{ix_n\cdot \xi_n}}{\;-i(\sum_{j=1}^{n-1} g^{jk} \xi_j\xi_k)^{\frac{1}{2}}+\xi_n}d\xi_n\bigg) \hat{f}(\xi') d\xi'\nonumber\\
\! \!&\!\!=\!\!&\!\!\frac{1}{(2\pi)^n} \int_{\mathbb{R}^{n}} e^{ix\cdot \xi} \bigg( \frac{
1}{(\sum_{j=1}^{n-1} g^{jk} \xi_j\xi_k)^{\frac{1}{2}}+i\xi_n}\bigg) \hat{f}(\xi') d\xi\nonumber.\end{eqnarray}}

\vskip 0.16 true cm

Let $t(x', \xi)$ be a $C^\infty$ function on $\mathbb{R}^{n-1} \times \mathbb{R}^n$ admitting the following series expansion (see p.$\,$29 of \cite{Monv}):
\begin{eqnarray} \label{17-11-11-4} t(x', \xi)= \sum_{l=0}^{r-1} \alpha_l (x', \xi') \xi^l_n + \sum_{m=0}^\infty a_m (x', \xi')
\frac{(\langle \xi'\rangle +i\xi_n)^m}{ (\langle \xi'\rangle -i \xi_n)^{m+1}}
,\end{eqnarray}
where $\alpha_l(x', \xi')$ belongs to $S^{d-l}_{1,0}$, and the $a_m(x',\xi')$ form a rapidly decreasing sequence in $S^{d+1}_{1,0}$.
The {\it trace operator} $T$
 of order $d$ with symbol $t(x',\xi)$ is the continuous operator: $C_0^\infty (\mathbb{\bar R}^n_+)\to C^\infty (\mathbb{R}^{n-1})$ defined by
   \begin{eqnarray} \label{17-11-11-5} (Tf)(x')=\frac{1}{(2\pi)^n}\int_{\mathbb{R}^{n-1}}  e^{ix'\cdot \xi'}\,d\xi' \int^{+} \,
   t(x', \xi) \hat{f}(\xi)\, d\xi_n  ,\end{eqnarray}
    where $\hat{f}$  is the Fourier transform of the extension of $f$ by $0$ for $x_n<0$. We will say that $T$ is of class $r$ if $r$ is the integer limiting the first sum in (\ref{17-11-11-4}).

\vskip 0.26 true cm

 Let $g(x',\xi', \xi_n, \eta_n)$ be a $C^\infty$ function on $\mathbb{R}^{n-1}\times \mathbb{R}^{n-1}\times \mathbb{R}\times \mathbb{R}$
admitting a series expansion:
\begin{eqnarray}  \label{2019.12.30-8}\;\;\quad \quad \; g(x', \xi', \xi_n, \eta_n) \!= \!\sum\limits_{s=0}^{r-1} \! \kappa_s (x', \xi', \xi_n) \eta_n^s \!+\!\!\sum_{
l\ge 0, m\ge 0}\!\! a_{lm} (x', \xi')\frac{ (\langle \xi'\rangle- i\xi_n)^l}{ ( \langle \xi'\rangle \!+\!i\xi_n)^{l+1}} \,\frac{(\langle \xi'\rangle +i \eta_n)^m}{( \langle \xi' \rangle \!-\!i\eta_n)^{m+1}},\end{eqnarray}
where $\kappa_s\in {\mathcal{K}}^{d-s}$ is the symbol of a Poisson operator of order $d-s$, and $a_{lm}$  a rapidly decreasing double sequence in $S_{1,0}^{d+1}$. The singular Green operator $G$ of order $d$ and class $r$ defined
by the symbol $g$ is the operator $G: C_0^\infty ({\mathbb{\bar R}}^n_{+}) \to  C^\infty ({ \mathbb{\bar R}}^n_{+})$ defined by (see p.$\,$30 of \cite{Monv}):
\begin{eqnarray} \label{2020.9.18-1} Gf(x) =\frac{1}{(2\pi)^{n+1} } \int_{\mathbb{R}^{n-1}} \! d\xi' \int^{+} e^{ix\cdot \xi} d\xi_n \int^{+} g(x', \xi', \xi_n, \eta_n) \hat{f} (\xi', \eta_n) d\eta_n \end{eqnarray}
(here $\hat{f}$ is the Fourier transform of the extension of $f$ by $0$ for $x_n<0$).

\vskip 0.15 true cm

Let us remark that since  \begin{eqnarray}\label{2020.10.27-1}  L(x,y) =\frac{1}{(2\pi)^{n+1} } \int_{\mathbb{R}^{n-1}}\! d\xi' \int^{+} e^{i(x-y)\cdot \xi} d\xi_n \int^{+} g(x', \xi', \xi_n, \eta_n) d\eta_n\end{eqnarray}  is the Schwartz kernel of the singular Green operator $G$ (cf. \cite{Gr} or (18.1.7) of \cite{Ho3}), we see that $Gf$ can also be written as \begin{eqnarray*}Gf(x) =\int L(x,y)f(y)dy. \end{eqnarray*}
Here  the Schwartz kernel is generated by the expression  \begin{eqnarray} \label{2020.9.28-2} G\delta(x-y)  =\frac{1}{(2\pi)^{n+1} } \int_{\mathbb{R}^{n-1}} d\xi' \int^{+} e^{i(x-y)\cdot \xi} d\xi_n \int^{+} g(x', \xi', \xi_n, \eta_n) \hat{\delta} (\xi', \eta_n) d\eta_n \end{eqnarray}
and by applying the fact the Fourier transform of the Dirac delta function $\delta$ is $1$. In particular, (\ref{2020.9.18-1}) still holds when
$f\in C_0^\infty (\mathbb{R}^n)$ or $f\in \mathscr{S} ( \mathbb{\bar R}^n)$ or $f\in C^1( \mathbb{R}^n) \cap \mathscr{S} (\mathbb{\bar R}^n_+) \cap \mathscr{S} (\mathbb{\bar R}^n_-)$.

\vskip 0.14 true cm

An operator $M$ is said to be a Green operator if it is the  sum of a pseudodifferential operator $P$ and a singular Green operator $G$. We denote by $\big(\phi(x,\xi); \psi(x', \xi',\xi_n, \eta_n)
\big)$ the symbol of this Green operator $M$, where $\phi(x, \xi)\in S^m_{1,0}$ and $\psi(x', \xi',\xi_n, \eta_n)$ is as in (\ref{2019.12.30-8}).

\vskip 0.14 true cm

By means of local coordinates and a partition of unity all operators above can immediately be defined in the same way in a bounded domain $\Omega$ with boundary $\partial\Omega$. As pointed out in \cite{Monv}, a singular Green operator $G$ of order $d$ and class $r$ is continuous $C_0^\infty
(\bar\Omega) \to C^\infty (\bar \Omega)$ and extends continuously $H_{{comp}}^s (\bar\Omega)\to H_{{loc}}^{s-d} (\bar\Omega)$ if $s>r+\frac{1}{2}$.

 \vskip 0.28 true cm

\noindent  {\bf Lemma 3.4.} \   {\it  Let $\Omega$ be a $C^\infty$-manifold with boundary $\partial \Omega$. For each $j\in \{1,2\}$, let $P_j$ (respectively, $G_j$, $K_j$, $T_j$, $Q_j$) be a pseudodifferential operator on $\Omega$ satisfying the transmission property (respectively, singular Green operator,  Poisson operator, trace operator, and the pseudodifferential operator on the boundary). Then the following conclusions hold}:

(i) \ \  $P_jP_l$,  ($1\le j,l\le 2$), {\it are all pseudodifferential operators on $\Omega$};

(ii) \ \  $P_j G_l$, $G_jP_l$, $G_jG_l$, $K_jT_l$, ($1\le j,l\le 2$), {\it are all Green operators on $\Omega$};

(iii) \ \  $Q_j Q_l$ {\it and} $ T_j K_l$,  ($1\le j,l\le 2$), {\it are all paseudodifferential operators on $\partial \Omega$.}

\vskip 0.18 true cm

  \noindent  {\it Proof.} \ It is well-known (see, for example,  Proposition 3.3 on p.$\,$13 of \cite{Ta2}) that the composition of two pseudodifferential operators is still a pseudodifferential operator, which implies (i). Now, the other conclusions directly follow from (i) and Boutet de Monvel's algebra construct theory (see \cite{Monv}), which states that the matrix operators \begin{eqnarray} \label{2020.9.16-1} W= \begin{matrix} \left(\begin{matrix} P_\Omega + G  & K\\
 T & Q \end{matrix} \right)\end{matrix} \end{eqnarray} form an ``algebra'' by the sum and the composition of two such kinds of matrices,
 where $P_\Omega$ is a pseudodifferential operator satisfying the transmission property, $K$ a Poisson operator, $T$ a trace operator, and $Q$ a pseudodifferential operator on the boundary. \qed

\vskip 0.29 true cm

 As pointed out in \cite{Monv},  two kinds of symbols are given: the first kind is the interior symbol in $W$.  This is just the
 symbol $p(x, \xi)$ of the pseudodifferential operator $P_\Omega$ in $W$,
another is the boundary symbol of $W$. This is a Wiener-Hopf operator (see p.$\, 33$ of \cite{Monv}), whose matrix is
\begin{eqnarray*} \begin{pmatrix} p(\xi_n) +g(\xi_n, \eta_n) & k(\xi_n) \\
 t(\xi_n) & q\end{pmatrix}, \end{eqnarray*}
here we have written $p(\xi_n)$ instead of $p (x',\xi', \xi_n)$ etc. ...

\vskip 0.19 true cm

The boundary and interior symbols of an operator $W$ in (\ref{2020.9.16-1}) are related: the
coefficient $p(\xi_n)$ in the matrix $\begin{pmatrix} p +g & k \\
 t & q\end{pmatrix}$  is the restriction to the boundary of the
interior symbol in $W$. The following composition formulas of boundary symbols hold (see p.$\,$17 of \cite{Monv}):
  \begin{eqnarray} \label{2020.11.12-4}&& p \circ g= h^+_{\xi_n}[ p(\xi_n)\cdot g(\xi_n, \eta)],\\
   \label{2020.11.12-5} && g \circ p= h^{-}_{\eta_n} [ g(\xi_n, \eta_n)\cdot p(\eta_n) ],\\
  \label{2020.11.12-6} && g \circ g= \frac{1}{2\pi} \int^+ g_1 (\xi_n, s) g_2 (s, \eta_n) ds,\end{eqnarray}
where $h^+_{\xi_n}$ (respectively, $h^-_{\eta_n}$) is the projection on $H^+_{\xi_n}$ parallel to $H^{-}_{\xi_n}$ (respectively, $h^{-}_{\eta_n}= 1-h^{+}_{\eta_n}$), and $H^{+}_{\xi_n}$ is the space $H^+$ of complex-valued functions with variable $\xi_n$.
In particular, it follows from \S 4 of \cite{Monv} that
\begin{eqnarray} \label{2020.10.10-1,} \left. \begin{array} {ll} &\sigma_{\Omega} (P+G) = \sigma_{\Omega} (P)+ \sigma_{\Omega} (G), \quad
\sigma_{\Omega} (P G) = \sigma_{\Omega} (P) \cdot \sigma_{\Omega} (G),\\
& \sigma_{\partial \Omega} (P+G) = \sigma_{\partial \Omega} (P)+ \sigma_{\partial \Omega} (G), \quad
\sigma_{\partial \Omega} (P G) = \sigma_{\partial \Omega} (P) \cdot \sigma_{\partial \Omega} (G),\end{array} \right.\end{eqnarray}
where $\sigma (P) $, $\sigma(G)$ and $\sigma(P+G)$ respectively are the principal interior symbols of the pseudodifferential operator $P$, the singular Green operator $G$ and the Green operator $P+G$; and $\sigma_{\partial \Omega} (P)$, $\sigma_{\partial \Omega}(G)$ and $\sigma_{\partial \Omega} (P+G)$ respectively  are principal boundary symbol of $P$, $G$ and $P+G$.

Another equivalent definition for the singular Green operator (respectively, Poisson operator, trace operator) and the symbol composition rules (i.e., rules of calculus) we refer the reader to \cite{Gr} or \cite{Gr5}.

\vskip 0.20 true cm

 To obtain the spectral asymptotic result, we need to consider the asymptotic expansion of the Stokes semigroup in the interior and near boundary for a bounded domain $\Omega\subset {\mathbb{R}}^n$ with smooth boundary.
   Suppose that locally the boundary $\partial \Omega$ of the domain $\Omega$ is given by $C^\infty$-smooth functions $z_j=z_j(y_1,\cdots, y_{n-1}),\, j=1,\cdots, n$, of the parameters $y_1, \cdots, y_{n-1}$, chosen so that $y_j=\mbox{const}$ is a line of curvature. In vector notation, this can be written as $\mathbf{z}= \mathbf{z}({\mathbf{y}}')$, where ${\mathbf{y}}'=(y_1,\cdots, y_{n-1})$.
Then the first and second quadratic forms on $\partial \Omega$  take the form
 \begin{eqnarray} I=\sum_{j=1}^{n-1} E_j(y')(dy')^2, \quad II=\sum_{j=1}^{n-1} L_j(y') (dy')^2.\end{eqnarray}
  It is well-known that
\begin{eqnarray} \label{17.9.14.1} \qquad \frac{\partial \mathbf{z}}{\partial y_j}\cdot \frac{\partial \mathbf{z}}{\partial y_k}= \sum_{l=1}^n \frac{\partial z_l}{\partial y_j}\,\frac{\partial z_l}{\partial y_k}=E_j\delta_{jk}, \quad \, k_j=\frac{L_j({\mathbf{y}}')}{E_j({\mathbf{y}}')},\end{eqnarray}
where $\delta_{jk}$ is the Kronecker symbol, and $k_j$ are the principal normal curvatures of the surface $\partial \Omega$.
 In a neighborhood of $\partial \Omega$ we introduce coordinates  $y_1,\cdots, y_n$, where $y_n$ is the distance from the point $\mathbf{x}=(x_1,\cdots, x_n)$ to $\partial \Omega$.
 Then \begin{eqnarray*} \mathbf{x}=\mathbf{z}({\mathbf{y}}') +y_n {\boldsymbol{\nu}}({\mathbf{y}}'), \quad \frac{\partial {\boldsymbol{\nu}}}{\partial y_j}=-k_j \,\frac{\partial \mathbf{z}}{\partial y_j} \quad (j=1,\cdots, n-1).\end{eqnarray*}
 By this way we obtain a Riemannian  metric $g=(g_{jk}(x))$  in this neighborhood  of $\partial \Omega$ which has the form
 (see p.$\,$1101 of \cite{LU} or p.$\,$532 of \cite{Ta2})
\begin{eqnarray} \label{18/a-1} \; \quad\;\big[g_{jk} (x',x_{n}) \big]_{n\times n} = \begin{bmatrix}
 g_{11} (x',x_{n}) & g_{12} (x',x_{n}) & \cdots & g_{1,n-1} (x',x_{n}) & 0\\
 \cdots\cdots& \cdots\cdots & \cdots &\cdots\cdots  & \cdots\\
 g_{n-1,1} (x',x_{n})  & g_{n-1,2} (x',x_{n}) & \cdots & g_{n-1,n-1} (x',x_{n}) & 0\\
 0& 0& 0& 0&1 \end{bmatrix}.  \end{eqnarray}
 In particular,  according to the choice of boundary coordinate frame we have that for each $x_0\in \partial \Omega$, $$ \frac{1}{2}\,\frac{\partial g_{jk}}{\partial x_n} (x_0) =\kappa_k\delta_{jk}  \;\;  \mbox{for all} \;\; 1\le j,k \le n-1.$$

    Let $(\Omega, g)$ be a Riemannian  manifold. It is well-known (see \cite{Ta2}) that in the local coordinate system the
     usual divergence operator, gradient operator and Laplace-Beltrami operator,
    respectively, take the forms
\begin{eqnarray*} \label{9-2.1} \mbox{div} \,{X} :=\sum\limits_{j=1}^n \frac{1}{\sqrt{|g|}}\, \frac{\partial(\sqrt{|g|} \,{{X}}^j)}{\partial x_j}\quad \, \mbox{if}\;\; {{X}}=(X^1, \cdots, X^n)\in T\Omega,\end{eqnarray*}
\begin{eqnarray*} \label{9-2.2} \mbox{grad}\, v = \bigg(\sum_{k=1}^n g^{1k} \frac{\partial v}{\partial x_k}, \cdots, \sum_{k=1}^n g^{nk} \frac{\partial v}{\partial x_k}\bigg)\quad \, \mbox{if}\;\; v\in C^\infty(\Omega), \end{eqnarray*}
and  \begin{eqnarray*} \label{18/12/22-2} \Delta_g:= \mbox{div}\; \mbox{grad} = \frac{1}{\sqrt{|g|}} \sum_{j,k=1}^n \frac{\partial}{\partial x_j} \bigg(\sqrt{g}\,g^{jk} \frac{\partial}{\partial x_k}\bigg). \end{eqnarray*}

  \vskip 0.20 true cm

For the operators $\mbox{div}\,G_1\,\mbox{grad}$, $G_1\,\mbox{grad}\, \mbox{div}$ and $\mbox{grad}\, \mbox{div}\, G_1$ defined on Riemannian manifold $(\Omega, g)$, we will show that the operators $\mbox{div}\, G_1\mbox{grad}-I$ and $G_1\mbox{grad}\, \mbox{div}-\mbox{grad}\, \mbox{div}\, G_1$ are Green operators of order $-1$ respectively defined on $H_0^1(\Omega)$ and $[H_0^1(\Omega)]^n$.

  \vskip 0.28 true cm

 \noindent  {\bf Lemma 3.5.} \   {\it (i) \  The operator} $\mbox{div}\, G_1\,\mbox{grad}-${\it $I=P+Q$ realizes a continuous mapping $H^1_0(\Omega)\to H^2(\Omega)$, where $P\in OPS^{-1}_{1,0}$ and $Q$ is a singular Green operator of order $-1$. In particular, in the representation of the principal symbol of the Green operator
 $P+Q$, its each term must contain a factor of the form $\frac{\partial g^{lm}(x)}{\partial x_r}$, where $(g^{lm}(x))$ is the inverse of the Reimannian metric  $(g_{lm}(x))$ on $\Omega$;

 (ii) $\,$  The operator} $G_1\,\mbox{grad}\, \mbox{div}-\mbox{grad}\, \mbox{div}\,${\it $G_1=R+U$  realizes a continuous mapping $[H_0^1(\Omega)]^n\to [H^2(\Omega)]^n$,  where $R\in OPS^{-1}_{1,0}$ and $U$ is a singular Green operator of order $-1$. In the representation of the principal symbol of the Green operator $R+U$, its each term must contain a factor of the form $\frac{\partial g^{lm}(x)}{\partial x_r}$.}

  \vskip 0.27 true cm

  \noindent  {\it Proof.} \  Let $\Phi$ be a parametrix of the Laplace-Beltrami operator $\Delta_g$ on Riemannian manifold $(\Omega,g)$.
Obviously, $\Phi$ is a pseudodifferential operator, its symbol has the form:
\begin{eqnarray*}&& \frac{-1}{ \sum_{j,k=1}^n g^{jk} (x)\xi_j\xi_k} + i\bigg[  \frac{-1}{(\sum_{j,k=1}^n g^{jk}(x) \xi_j\xi_k)^2} \sum_{j,k=1}^n \frac{1}{\sqrt{|g(x)|}}  \frac{\partial (\sqrt{|g(x)|} g^{jk}(x) )}{\partial x_j}\xi_k \\
&& \quad \;\;\quad \;\; +\frac{1}{(\sum_{j,k=1}^n g^{jk}(x) \xi_j\xi_k)^3} \sum_{j,k,l=1}^n 2 g^{jl} (x)   \frac{\partial g^{jk(x)}}{\partial x_l} \xi_j^2 \xi_k\bigg]+\cdots, \end{eqnarray*}
where the cdots denote the symbol which belongs to $S^{-4}_{1,0}$.
From symbol composition formula (\ref{17-9-12.1}), we get that the $s$-th component  of the full symbol of the pseudodifferential operator $\Phi\, \mbox{grad}$ is
\begin{eqnarray} \label{2020.10.11-5} \frac{ -i\sum_{m=1}^n g^{sm}\xi_m }{\sum_{j,k=1}^n g^{jk} \xi_j\xi_k} + \Theta_s (x, \xi)+\cdots, \end{eqnarray}
 where
\begin{eqnarray} \label{2020.10.11-1}   && \Theta_s (x, \xi)=
 \frac{1}{(\sum_{j,k=1}^n g^{jk}(x) \xi_j\xi_k)^2}\!\!
\sum_{k,m=1}^n\!\! \Big(\sum_{j=1}^n\!\frac{1}{\sqrt{|g(x)|}} \frac{\partial (\sqrt{|g(x)|} g^{jk}(x))}{\partial x_j} g^{sm}(x)\\
&&\quad + 2 g^{kr}(x) \sum_{r=1}^n \!\frac{\partial g^{sm}(x)}{\partial x_r} \! \Big) \xi_k \xi_m   - \frac{1}{ (\sum_{j,k=1}^n \!\!g^{jk} (x)\xi_j\xi_k )^3}\!\! \sum_{j,k,l,m=1}^n \!2 g^{jl} (x)\frac{\partial g^{jk}(x)}{\partial x_l} g^{sm}(x)  \xi_j^2\xi_k\xi_m,\nonumber \end{eqnarray}
and the cdots in (\ref{2020.10.11-5}) denote the pseudodifferential symbol of order $-3$.
Applying symbol composition formula (\ref{17-9-12.1}) again, we get the full symbol  of the operator $\mbox{div}\, \Phi\, \mbox{grad}$ is $I+p(x,\xi)$, where
\begin{eqnarray}\label{2020.11.22-1}  p(x, \xi)= \sum_{s=1}^n\big( i \xi_s + \sum_{h=1}^n \Gamma_{hs}^h(x) \big)  \Theta_s (x, \xi) +\sum_{h,m,s=1}^n \Gamma_{hs}^h \Big(\! -\frac{i g^{sm}\xi_m}{\sum_{j,k=1}^n g^{jk} \xi_j\xi_k}\Big)+q(x,\xi),\end{eqnarray}
 $\Gamma_{lk}^j= \frac{1}{2} \sum_{m=1}^n g^{jm} \big( \frac{\partial g_{km}}{\partial x_l} +\frac{\partial g_{lm}}{\partial x_k} -\frac{\partial g_{lk}}{\partial x_m}\big)$ are the Christoffel symbols associated with the metric $g$ (see, for example, \cite{Ta2}), and the $q(x,\xi)\in S^{-2}_{1,0}$.
Thus  \begin{eqnarray*} \mbox{div}\, \Phi \, \mbox{grad} \, f= f+Pf, \quad  \mbox{for all}\,\,
 f\in C_0^\infty(\Omega),\end{eqnarray*}
where the full symbol of the pseudodifferential operator $P$ is $p(x, \xi)$.
In view of $G_1= \Phi -G_2 \gamma_0 \Phi$ (here $G_1, G_2$ and $\gamma_0$ are defined in (\ref{2-6}), (\ref{2-7}) and (\ref{a-2.1})),
we get \begin{eqnarray*} \mbox{div}\, G_1\, \mbox{grad}= \mbox{div}\, \Phi\, \mbox{grad}-\mbox{div}\, G_2\gamma_0\Phi \, \mbox{grad}= I- \mbox{div}\, G_2 \gamma_0 \Phi \,\mbox{grad}+P\end{eqnarray*}
    so that
\begin{eqnarray} \label{17-11-16-10}(\mbox{div}\, G_1\, \mbox{grad}-I)f=(-\mbox{div}\, G_2\gamma_0 \Phi\, \mbox{grad}+P)f.\end{eqnarray}
It remains to show that $\mbox{div}\, G_2\gamma_0\Phi\,\mbox{grad}$ is a Green operator on $\Omega$ (i.e., the sum of a pseudodifferential operator and a singular Green operator) of order $-1$. It follows from \cite{Monv} that (\ref{2020.9.16-1}) forms an ``algebra''---i.e., the sum and composition of two such operators is another one. Therefore, by (ii) of Lemma 3.4 we see that $G_2\gamma_0$ is a Green operator on $\Omega$. In view of $\mbox{div}$, $\Phi$ and $\mbox{grad}$ are all pseudodifferential operators having transmission property (see, (2.4)---(2.5) of \cite{Monv}), it follows from  (ii) of Lemma 3.4 that $-\mbox{div}\, G_2\gamma_0 \Phi\, \mbox{grad}$ is a Green operator  on $\Omega$. From the right-hand side of (\ref{17-11-16-10}) we immediately see that $\,\mbox{div}\, G_1\, \mbox{grad}-I\,$ is still a Green operator on $\Omega$. Moreover, Kozhevnikov proved (see p.$\,$12 of \cite{Ko}) that  $\mbox{div}\, G_2\gamma_0 \Phi\, \mbox{grad}$ is an operator of order $-1$.
Applying (\ref{2020.10.10-1,}) and (\ref{2020.10.11-1}), we get that  the principal symbol of $\mbox{div}\, G_2\gamma_0 \Phi\, \mbox{grad}$ is a sum of finitely many terms, and its each term certainly contains either the factor $\Theta_s (x, \xi)$ or the factor $\Gamma_{sh}^s$. Consequently, $(\mbox{div}\, G_1 \,\mbox{grad}\, - I) f = Jf$  for all $f\in C_0^\infty (\Omega)$, where $J:=P-\mbox{div}\, G_2 \gamma_0 \Phi \, \mbox{grad}\,$ is a Green operator of  order $-1$. In particular, by (\ref{2020.10.11-1})--(\ref{17-11-16-10}) we see that the principal symbol of $J$ is a sum of finitely many terms, each term of which also contains a factor of the form $\frac{\partial g^{jk}(x)}{\partial x_m}$. The first assertion i) of Lemma 3.5 is proved, since $C_0^\infty(\Omega)$ is dense in $H_0^1(\Omega)$.

The proof of the second assertion ii) is analogous to that of the assertion i). $\quad \;\;\square$

\vskip 0.22 true cm

If a Green operator is composed of some pseudodifferential operators and singular Green operators, then we may calculate the symbol of this Green operator by applying symbol formulas (\ref{17-9-12.1}) (for pseudodifferential operators) and (\ref{2020.11.12-4})---(\ref{2020.11.12-6}) (for singular Green operators). We now calculate the symbol of Green operator $(-\mu \Delta- \lambda)^{-1} \big[(B_2+C_2)(-\mu \Delta- \lambda)^{-1}\big]^{k}$ for any $k\ge 1$.
 Let $p(x, \xi, \lambda)$, $b(x, \xi)$ and $g(x', \xi',\xi_n, \eta_n)$ be the symbols of pseudodifferential operator $(-\mu\Delta- \lambda)^{-1}$, pseudodifferential operator $B_2$ and singular Green operator $C_2$, respectively. For brevity, the symbols $p(x, \xi,\lambda)$, $b(x, \xi)$ and $g(x',\xi',\xi_n,\eta_n)$ respectively are written as $p(\xi_n)$, $b(\xi_n)$ and $g(\xi_n,\eta_n)$.

For $k=1$, it is clear by (\ref{17-9-12.1}) that the pseudodifferential operator $(-\mu \Delta- \lambda)^{-1} B_2(-\mu \Delta- \lambda)^{-1}$ has symbol
\begin{eqnarray} \label{2020.11.24-1} \phi_1(x, \xi, \lambda):=\sum_{\alpha\ge 0} \frac{(-i)^{|\alpha|}}{\alpha!} \partial^{\alpha}_\xi p(x, \xi,\lambda) \cdot \partial_x^\alpha \big(\sum_{\beta\ge 0} \frac{ (-i)^{|\beta|} } {\beta!} \partial_\xi^\beta b(x, \xi) \cdot \partial_x^\beta p(x, \xi, \lambda)  \big).\end{eqnarray}
By (\ref{2020.11.12-4}) and (\ref{2020.11.12-1}), we get that the symbol of the singular Green operator $(-\mu\Delta- \lambda)^{-1}C$ is
 \begin{eqnarray*} h_{\xi_n}^{+} [p(\xi_n)\cdot g(\xi_n, \eta_n)] = -\frac{1}{2\pi i} \int^+ \frac{p(\tau_1) g(\tau_1,\eta_n)}{\tau_1- \xi_n} d\tau_1:= c_1(\xi_n, \eta_n),\end{eqnarray*}
where $\int^+$ is defined by (\ref{17-11-11-0}).
 Furthermore, from (\ref{2020.11.12-5}) and (\ref{2020.11.12-1}) we know the symbol of the singular Green operator $(-\mu\Delta -\lambda)^{-1} C(-\mu\Delta- \lambda)^{-1}$ is
 \begin{eqnarray} \label{2020.11.17-4}  h^-_{\eta_n} [c_1 (\xi_n, \eta_n) \cdot p(\eta_n)]\, = \!\!\!\!&\!\!&\!\!\! \frac{1}{2\pi i} \int^\oplus
 \frac{c_1(\xi_n, \tau_2) p(\tau_2)}{ \tau_2-\eta_n} d\tau_2 \\
= \!\!\!\!&\!\!&\!\!\!\!- \frac{1}{(2\pi i)^2}  \int^+ d\tau_1 \int^\oplus \frac{ p(\tau_1) g(\tau_1, \tau_2) p(\tau_2) }{ (\tau_1- \xi_n) (\tau_2- \eta_n)} \,d\tau_2:= c_2 (\xi_n, \eta_n), \nonumber \end{eqnarray}
where  $\int^\oplus\!\frac{f(z)}{z-s}\,dz$ denotes the contour integral $\int_\gamma \!\frac{f(z)}{z-s}\,dz$ which is required to contain $s$ in the interior of $\gamma$, and $\gamma$ is a large circle in the upper plane $\mbox{Im}\, z>0$.

 For $k=2$, by (\ref{17-9-12.1}) we see that the symbol of the pseudodifferential operator $(-\mu \Delta -\lambda)^{-1} \big[B(-\mu\Delta -\lambda)^{-1}\big]^2$ is
   \begin{eqnarray} \label{2020.11.24-2} \phi_2(x, \xi, \lambda):=\sum_{\gamma \ge 0} \frac{(-i)^{|\gamma|}}{\gamma!} \partial^\gamma_\xi \phi_1(x, \xi, \lambda) \cdot \partial^\gamma_x
  \Big( \sum_{\iota\ge 0} \frac{(-i)^{|\iota|}}{\iota!} \partial^\iota_\xi b(x, \xi) \cdot \partial_x^\iota p(x, \xi,\lambda)\Big), \end{eqnarray}
  where $\phi_1(x, \xi, \lambda)$ is given by (\ref{2020.11.24-1}).  From (\ref{2020.11.12-5}) and the right-hand side of (\ref{2020.11.17-4}) we get that the symbol of the singular Green operator $(-\mu \Delta -\lambda)^{-1} C(-\mu\Delta -\lambda)^{-1} B$ is
   \begin{eqnarray*} \label{2020.11.15-1}  && h^{-}_{\eta_n} [ c_2(\xi_n, \eta_n) b(\eta_n)] =\frac{1}{2\pi i} \int^\oplus \frac{c_2(\xi_n, \tau_3) b(\tau_3)}{ \tau_3 -\eta_n }\, d\tau_3 \\
&& \qquad =  -\frac{1}{(2\pi i)^3}\int^+  d\tau_1 \int^\oplus d\tau_2 \int^\oplus \frac{p(\tau_1) g(\tau_1, \tau_2) p(\tau_2) b(\tau_3)}{ (\tau_1-\xi_n ) (\tau_2-\tau_3) (\tau_3- \eta_n)} d\tau_3  := c_3 (\xi_n, \eta_n), \nonumber\end{eqnarray*}
so that the symbol of the singular Green operator $(-\mu \Delta -\lambda)^{-1} C(-\mu\Delta -\lambda)^{-1} B(-\mu\Delta -\lambda)^{-1}$ is
\begin{eqnarray} \label{2020.11.15-2} &&  h_{\eta_n}^{-} [c_3 (\xi_n, \eta_n) p(\eta_n)] = \frac{1}{2\pi i} \int^\oplus \frac{ c_3 (\xi_n, \tau_4) p(\tau_4)}{\tau_4 -\eta_n} \, d\tau_4 \\
&&\qquad  \quad =-\frac{1}{(2\pi i)^4}\!  \int^+\!d\tau_1\! \int^\oplus \!d\tau_2 \!\int^\oplus \!d\tau_3 \!\int^\oplus\!    \frac{  p(\tau_1) g (\tau_1, \tau_2) p(\tau_2) b(\tau_3) p(\tau_4)}{ (\tau_1\!-\! \xi_n) (\tau_2\!-\! \tau_3) (\tau_3\!-\! \tau_4) (\tau_4\!- \!\eta_n)}  d\tau_4.\nonumber\end{eqnarray}
Replacing $\tau_1$ and $\tau_2$ by $\tau_3$ and $\tau_4$ in (\ref{2020.11.17-4}) respectively, we can rewrite $c_2(\xi_n, \eta_n)$ as
\begin{eqnarray} c_2(\xi_n,\eta_n)=   -\frac{1}{(2\pi i)^2}  \int^+ d\tau_3 \int^\oplus \frac{ p(\tau_3) g(\tau_3, \tau_4) p(\tau_4) }{ (\tau_3- \xi_n) (\tau_4- \eta_n)} \,d\tau_4.\nonumber \end{eqnarray}
It follows from this and (\ref{2020.11.12-4}) that the symbol of the singular Green operator $B(-\mu \Delta -\lambda)^{-1} C(-\mu\Delta -\lambda)^{-1}$ is
\begin{eqnarray*} \label{2020.11.17-1}  && h^{+}_{\xi_n} [b(\xi_n)  c_2(\xi_n, \eta_n)] =-\frac{1}{2\pi i} \int^+ \frac{b(\tau_2)c_2(\tau_2, \eta_n)}{ \tau_2 -\xi_n }\, d\tau_2 \\
&& \qquad =  \frac{1}{(2\pi i)^3}\int^+ d\tau_2  \int^+ d\tau_3 \int^\oplus   \frac{b(\tau_2) p(\tau_3)  g(\tau_3, \tau_4) p(\tau_4)}{ (\tau_2-\xi_n ) (\tau_3-\tau_2) (\tau_4- \eta_n)} d\tau_4  := c_3 (\xi_n, \eta_n), \nonumber\end{eqnarray*}
and hence the symbol of the singular Green operator $(-\mu \Delta -\lambda)^{-1} B(-\mu\Delta -\lambda)^{-1} C(-\mu\Delta -\lambda)^{-1}$ is
\begin{eqnarray} \label{2020.11.1-2} &&  h_{\xi_n}^{+} [p(\xi_n) c_3 (\xi_n, \eta_n)] = -\frac{1}{2\pi i} \int^+ \frac{p(\tau_1)  c_3 (\tau_1, \eta_n) }{\tau_1-\xi_n} \, d\tau_1 \\
&&\qquad  \quad =-\frac{1}{(2\pi i)^4} \int^+ d\tau_1 \int^+ d\tau_2 \int^+ d\tau_3 \int^\oplus    \frac{  p(\tau_1)b(\tau_2) p(\tau_3)  g (\tau_3, \tau_4) p(\tau_4)}{ (\tau_1- \xi_n) (\tau_1- \tau_2) (\tau_2- \tau_3) (\tau_4- \eta_n)}  d\tau_4.\nonumber\end{eqnarray}
Finally, according to (\ref{2020.11.17-4}), we see that the symbol of the singular Green operator $(-\mu \Delta -\lambda)^{-1} C(-\mu\Delta -\lambda)^{-1}C$ is
\begin{eqnarray*} \label{2020.11.15-1}  && \frac{1}{2\pi} \!\int^+\!\!\! \!c_2 (\xi_n, \tau_3) g(\tau_3, \eta_n) d\tau_3 \!=\!-\frac{1}{2\pi (2\pi i)^2} \!\int^+\! \!d\tau_1 \int^\oplus\! \!d\tau_2  \!\int^+  \frac{p(\tau_1) g(\tau_1, \tau_2) p(\tau_2) g(\tau_3, \eta_n)}{ (\tau_1-\xi_n)(\tau_2- \tau_3)} d\tau_3\! := c_4 (\xi_n, \eta_n),\end{eqnarray*}
so the symbol of the  singular Green operator $(-\mu \Delta -\lambda)^{-1} C(-\mu\Delta -\lambda)^{-1}C(-\mu\Delta -\lambda)^{-1}$ is
\begin{eqnarray} \label{2020.11.1-6,}
&& h^{-}_{\eta_n} [  c_4(\xi_n, \eta_n)p(\eta_n)] =\frac{1}{2\pi i} \int^\oplus \frac{c_4(\xi_n, \tau_4)p(\tau_4)}{ \tau_4 -\eta_n }\, d\tau_4 \\
&& \qquad =  -\frac{1}{(2\pi)(2\pi i)^3}\int^+ d\tau_1\int^\oplus d\tau_2  \int^+ d\tau_3 \int^\oplus   \frac{p(\tau_1) g(\tau_1, \tau_2) p(\tau_2)  g(\tau_3, \tau_4) p(\tau_4)}{ (\tau_1-\xi_n ) (\tau_2-\tau_3) (\tau_4- \eta_n)} d\tau_4. \nonumber\end{eqnarray}
By summing the results in (\ref{2020.11.15-2})---(\ref{2020.11.1-6,}), we can get the symbol of the singular Green operator $(-\mu \Delta -\lambda)^{-1} [(B+C) (-\mu\Delta -\lambda)^{-1}]^2- (-\mu\Delta -\lambda)^{-1} [B(-\mu\Delta -\lambda)^{-1}]^2$.

 For general integer $k$, using  similar way as above, we can get
that the symbol of the pseudodifferential operator $(-\mu \Delta -\lambda)^{-1}\big[B(-\mu\Delta -\lambda)^{-1}\big]^{k}$ is  \begin{align} \label{2020.11.24-2,0} &\phi_{k} (x,\xi,\lambda) =\sum_{\gamma \ge 0} \frac{(-i)^{|\gamma|}}{\gamma!} \partial^\gamma_\xi \phi_{k-1}(x, \xi, \lambda) \cdot \partial^\gamma_x
  \Big( \sum_{\iota\ge 0} \frac{(-i)^{|\iota|}}{\iota!} \partial^\iota_\xi b(x, \xi) \cdot \partial_x^\iota p(x, \xi,\lambda)\Big), \end{align}
where $\phi_{k-1}(x, \xi, \lambda)$ is the full symbol of the pseudodifferential operator $(-\mu \Delta -\lambda)^{-1}\big[B(-\mu\Delta -\lambda)^{-1}\big]^{k-1}$.
 Let $\psi_{k-1} (x', \xi', \xi_n, \eta_n,\lambda)$  be the  symbol of the singular Green operator $(-\mu \Delta -\lambda)^{-1} [(B+C)(-\mu\Delta -\lambda)^{-1}]^{k-1} -(-\mu \Delta -\lambda)^{-1} [B(-\mu\Delta -\lambda)^{-1}]^{k-1}$.
 We will simply write $\psi_{k-1}(x', \xi',\xi_n,\eta_n, \lambda)$ as $\psi_{k-1} (\xi_n, \eta_n)$.
 Then the  symbol of the singular Green operator $(-\mu \Delta -\lambda)^{-1} [(B+C)(-\mu\Delta -\lambda)^{-1}]^{k-1}B$ is
 \begin{eqnarray} \label{2020.12.4-6}  h^-_{\eta_n} [\psi_{k-1}(\xi_n, \eta_n) b(\eta_n)] =\frac{1}{2\pi i} \int^\oplus \frac{\psi_{k-1}(\xi_n, \tau_{2k-1}) b(\tau_{2k-1}) }{ \tau_{2k-1} -\eta_n} d\tau_{2k-1} := \varpi (\xi_n, \eta_n),\end{eqnarray}
 so that the symbol of the singular Green operator  $(-\mu \Delta -\lambda)^{-1} [(B+C)(-\mu\Delta -\lambda)^{-1}]^{k-1}B(-\mu\Delta -\lambda)^{-1}$ is
 \begin{eqnarray} \label{2020.12.4-6}  h^-_{\eta_n} [\varpi (\xi_n, \eta_n) p(\eta_n)]\!\!\!\! &&\!\!\!\!=\frac{1}{2\pi i} \int^\oplus \frac{\varpi(\xi_n, \tau_{2k}) p(\tau_{2k})}{ \tau_{2k} -\eta_n} d\tau_{2k}\\
 \!\!\!\! &&\!\!\!\!=
 \frac{1}{(2\pi i)^2} \int^\oplus d\tau_{2k-1} \int^\oplus \frac{\psi_{k-1} (\xi_n, \tau_{2k-1}) b(\tau_{2k-1}) p(\tau_{2k})}{ (\tau_{2k}-\eta_n)(\tau_{2k-1}-\tau_{2k})} d\tau_{2k}.\nonumber\end{eqnarray}
  The  symbol of the singular Green operator $(-\mu \Delta -\lambda)^{-1} [(B+C)(-\mu\Delta -\lambda)^{-1}]^{k-1}C$ is
 \begin{eqnarray} \label{2020.12.4-6} \frac{1}{2\pi} \int^+
   \psi_{k-1}(\xi_n, \tau_{2k-1})  g(\tau_{2k-1},\eta_n) d\tau_{2k-1} =\varsigma (\xi_n, \eta_n),\end{eqnarray}
    and hence the  symbol of the singular Green operator  $(-\mu \Delta -\lambda)^{-1} [(B+C)(-\mu\Delta -\lambda)^{-1}]^{k-1}C(-\mu\Delta -\lambda)^{-1}$ is
 \begin{eqnarray} \label{2020.12.4-6}  h^-_{\eta_n} [\varsigma (\xi_n, \eta_n) p(\eta_n)]\!\!\!\! &&\!\!\!\!=\frac{1}{2\pi i} \int^\oplus \frac{\varsigma(\xi_n, \tau_{2k}) p(\tau_{2k})}{ \tau_{2k} -\eta_n} d\tau_{2k}\\
 \!\!\!\! &&\!\!\!\!=
 \frac{1}{2\pi (2\pi i)}\int^+d\tau_{2k-1} \int^\oplus  \frac{\psi_{k-1} (\xi_n, \tau_{2k-1}) g(\tau_{2k-1},\tau_{2k}) p(\tau_{2k})}{ (\tau_{2k}-\eta_n)} d\tau_{2k}.\nonumber\end{eqnarray}
  Consequently, the symbol of the singular Green operator $(-\mu \Delta -\lambda)^{-1} [(B+C)(-\mu\Delta -\lambda)^{-1}]^{k}-(-\mu \Delta -\lambda)^{-1}\big[B(-\mu\Delta -\lambda)^{-1}\big]^{k} $ is
  \begin{eqnarray}\label{2020.12.4-4,} && \frac{1}{(2\pi i)^2} \int^\oplus d\tau_{2k-1} \int^\oplus \frac{\psi_{k-1} (\xi_n, \tau_{2k-1}) b(\tau_{2k-1}) p(\tau_{2k})}{ (\tau_{2k}-\eta_n)(\tau_{2k-1}-\tau_{2k})} d\tau_{2k}\\
 && \qquad \quad  + \frac{1}{2\pi (2\pi i)}\int^+d\tau_{2k-1} \int^\oplus  \frac{\psi_{k-1} (\xi_n, \tau_{2k-1}) g(\tau_{2k-1},\tau_{2k}) p(\tau_{2k})}{ (\tau_{2k}-\eta_n)} d\tau_{2k}.\nonumber\end{eqnarray}

\vskip 1.39 true cm

\section{Decomposition of operators}

\vskip 0.45 true cm

 Now, we denote by $\mathbf{A}$ the extension of $S^{-1}$  to all  of $[L^2(\Omega)]^n$ that is given by the same formula as $S^{-1}$:
  \begin{eqnarray} \label{3-1} \mathbf{A} \mathbf{f} = -(1/\mu)G_1 [\mathbf{I}-2 \,\mbox{grad}({I}+G_2 K_{-1}\gamma_0)\mbox{div}\, G_1]\mathbf{f},\quad \; \, \forall\, \mathbf{f}\in [L^2(\Omega)]^n. \end{eqnarray}

  Since $[L^2(\Omega)]^n= J\oplus F\oplus E$, any element $\mathbf{f}\in [L^2(\Omega)]^n$ is thus uniquely determined by the column of its ``coordinates'' $({\mathbf{f}}_J, {\mathbf{f}}_F, {\mathbf{f}}_E)^T$, where $T$ denotes the transposition and ${\mathbf{f}}_J\in J$, ${\mathbf{f}}_F\in F$, ${\mathbf{f}}_E\in E$. The operator $\mathbf{A}$ of (\ref{3-1}) can be written correspondingly in the matrix form
    \begin{gather} \label{3-2}
   \mathbf{A} \mathbf{f} =\begin{pmatrix} {\mathbf{A}}_{{}_{JJ}}  & {\mathbf{A}}_{{}_{JF}} & {\mathbf{A}}_{{}_{JE}}\\
    {\mathbf{A}}_{{}_{FJ}} &  {\mathbf{A}}_{{}_{FF}} & {\mathbf{A}}_{{}_{FE}}\\
      {\mathbf{A}}_{{}_{EJ}} & {\mathbf{A}}_{{}_{EF}} & {\mathbf{A}}_{{}_{EE}} \end{pmatrix}\begin{pmatrix}{\mathbf{f}}_{{}_J}\\
     {\mathbf{f}}_{{}_F}\\
  {\mathbf{f}}_{{}_E} \end{pmatrix}, \end{gather}
   where, for example, ${\mathbf{A}}_{{}_{JF}}$ maps $F$ into $J$. We have the following:

\vskip 0.24 true cm

   \noindent  {\bf Lemma 4.1.} \    {\it The matrix (\ref{3-2}) representing the operator $\mathbf{A}$ has the form
     \begin{gather} \label{3-3-00}
   \mathbf{A} =\begin{pmatrix} S^{-1}  & 0 & 0\\
    0 &  {\mathbf{A}}_{FF} & 0\\
     0 & {\mathbf{A}}_{EF} & 0 \end{pmatrix}. \end{gather} }

  \vskip 0.2 true cm

  \noindent  {\it Proof.} \  i) \ Since $\mathcal{D} (S)$ is a subset of $J$, the operator $S^{-1}$ maps $J$ into itself. Thus, $J$ is an invariant subspace of ${\mathbf{A}}$ which is an extension of $S^{-1}$. This implies ${\mathbf{A}}_{FJ}=0$, ${\mathbf{A}}_{EJ}=0$ and ${\mathbf{A}}_{JJ}=S^{-1}$.

ii) \  It follows from the proof of Lemma of \S6 in \cite{Ko} that $\mathbf{A}\mathbf{f}=0$ for any $\mathbf{f}\in E$, which shows that the last column of the matrix $\mathbf{A}$ consists of zeros.

     iii) \  We will prove that ${\mathbf{A}}_{JF}{\mathbf{f}}=0$ for any ${\mathbf{f}}\in F$. Put ${\mathbf{u}}:= {\mathbf{A}}{\mathbf{f}}$ with $\mathbf{f}\in F$. Obviously, $\mathbf{f}\in \big[L^2 (\Omega)\big]^n$.
  Since $G_1= \Phi - G_2 \gamma_0 \Phi$, we see by (\ref{3-1}) that for any $\mathbf{f}\in F$,  \begin{eqnarray*}  \mathbf{A} \mathbf{f} \!\!\!\!\!&\!&\!\!\!\!
  =-(1/\mu)G_1 [\mathbf{I}-2 \,\mbox{grad}({I}+G_2 K_{-1}\gamma_0)\mbox{div}\, G_1]\mathbf{f}\\
 \!\!\!\!\!&\!&\!\!\!\!  = - (1/\mu) G_1 [\mathbf{I}-2(\mbox{grad}  + \mbox{grad}\, G_2 \, K_{-1} \gamma_0 ) (\mbox{div}\, \Phi -\mbox{div}\, G_2 \gamma_0 \Phi ) \big]\mathbf{f}\\
 \! \! \!\!\!&\!&\!\!\!\!= -(1/\mu) G_1 \big[\mathbf{I}\!-\!2\, \mbox{grad}\, \mbox{div}\, \Phi \!+\! 2\,\mbox{grad}\, \mbox{div}\, G_2 \, \gamma_0 \Phi \!- \!2\, \mbox{grad}\, G_2 \, K_{-1} \gamma_0 \,\mbox{div}\, \Phi \!+\! 2\,\mbox{grad}\, G_2 \, K_{-1}  \gamma_0 \, \mbox{div}\, G_2 \gamma_0  \Phi\big] \mathbf{f},\end{eqnarray*}
where $\Phi$ is the parametrix of the Laplace-Beltrami operator $\Delta_g$ on Riemannian manifold $(\Omega,g)$ (see, the proof of Lemma 3.5).
   Recall that $\mbox{div}$, $\mbox{grad}$ and $\Phi$ are all pseudodifferential operators having transmission property (see (2.4)---(2.5) of \cite{Monv}). It follows that $\mbox{div}\, \Phi$ and $\,-2\,\mbox{grad}\, \mbox{div}\, \Phi$ are  pseudodifferential operators  with orders $-1$ and $0$, respectively, each of which has transmission property.
      Furthermore, since $\gamma_0$ is the trace operator of order $0$ and $G_2$ is the Poisson operator of order $-1$, it follows from (ii) of Lemma 3.4 that $G_2 \gamma_0$ is a Green operator (i.e., the sum of a pseudodifferential operator and a singular Green operator) of order $-1$, so is $2\,\mbox{grad}\, \mbox{div}\, G_2  \gamma_0 \Phi$.
   Applying Proposition (3.11) on p.$\,$30 in \cite{Monv}, we obtain that $K_{-1} \gamma_0 (\mbox{div}\, \Phi)$ is a trace operator of order $-2$. It follows from (ii) of Lemma 3.4 (or the ``algebra'' property of
(\ref{2019.12.30-1}))  that $G_2 K_{-1} \gamma_0 \,\mbox{div}\, \Phi$ is a Green operator (on $\Omega$) of order $-3$. Using (ii) of Lemma 3.4 again, we get
 that   $\mbox{grad}\,G_2K_{-1} \gamma_0\,\mbox{div}\,\Phi$ is a Green operator (on $\Omega$) of order $-2$.
 By a similar method, we can also get that $2\,\mbox{grad}\, G_2 \, K_{-1} \, \gamma_0 \, \mbox{div}\, G_2 \gamma_0  \Phi$ is a Green operator of order $-3$.
 Consequently, $[\mathbf{I}-2 \,\mbox{grad}({I}+G_2 K_{-1}\gamma_0)\mbox{div}\, G_1]$ is a Green operator of order $0$. We let $r\ge 0$ be the class for the corresponding singular Green operator associated with the Green operator $[\mathbf{I}-2 \,\mbox{grad}({I}+G_2 K_{-1}\gamma_0)\mbox{div}\, G_1]$ (see the definition of singular Green operators in section 3). We first consider the case of $\mathbf{f}\in F\cap [H^s(\Omega)]^n$, where $s>r+\frac{1}{2}$.  Since $\Omega$ is a bounded domain, we know that the singular Green operator of order $0$ and class $r$ maps $H^s (\Omega)$ to $H^s(\Omega)$ if $s>r+\frac{1}{2}$ (see section 3) and that pseudodifferential operator of order $0$ maps $H^s(\Omega)$ to $H^s (\Omega)$ for all $s\in \mathbb{R}$ (see section 2),  which implies $[\mathbf{I}-2 \,\mbox{grad}({I}+G_2 K_{-1}\gamma_0)\mbox{div}\, G_1] \mathbf{f}\in [H^s(\Omega)]^n$ for $s>r+\frac{1}{2}$.
  It follows from the definition of $G_1$ that $\mathbf{u}= \mathbf{A}\mathbf{f}=G_1[\mathbf{I}-2 \,\mbox{grad}({I}+G_2 K_{-1}\gamma_0)\mbox{div}\, G_1] \mathbf{f}\in [H^{s+2}(\Omega)]^n$ and $\mathbf{u}\big|_{\partial \Omega}=0$. On the other hand, by the assumption of the lemma we have $\mathbf{u}=\mathbf{A} \mathbf{f} \in J$, which implies that $\mathbf{u}\in \mathcal{D}(S)$. Hence we have $S\mathbf{u}=\mathbf{f}\in J$. Since $\mathbf{f}$ also belongs to $F$, we immediately obtain $\mathbf{f}=0$. In other words, $\mathbf{u}$ is a solution of the following Stokes equations
   \begin{eqnarray} \label{2020.12.1-1}  \left\{ \begin{array}{ll}   -\mu \Delta {\mathbf{u}}+ \nabla p  ={\mathbf{0}}  \;\; &\mbox{in}\;\;  \Omega,\\
    \mbox{div}\; {\mathbf{u}}=0 \;\; & \mbox{in}\;\;  \Omega, \\
    {\mathbf{u}}=0 \;\;& \mbox{on}\;\;  \partial \Omega. \\
    \end{array}  \right.\end{eqnarray}
   By uniqueness of solutions of the Stokes equations (see \cite{Lad} or \cite{Te}), we get that $\mathbf{u}=0$, i.e., $\mathbf{A}\mathbf{f}=0$ if $\mathbf{f}\in F\cap [H^s(\Omega)]^n$, where $s>r+\frac{1}{2}$.
   Now, from (\ref{3-1}) and the compactness of $G_1$, we see that the operator $\mathbf{A}$ is compact on  $[L^2(\Omega)]^n$, and hence $\mathbf{A}_{JF}$ is a compact operator from $F$ to $J$. Since $F\cap [H^s(\Omega)]^n$ is dense in $F$ with respect to $L^2$-norm, we immediately get that $\mathbf{A}_{JF}\mathbf{f}=0$ for all $\mathbf{f}\in F$. That is, $\mathbf{A}_{JF}=0$.
 Thus the matrix (\ref{3-2}) representing the operator $\mathbf{A}$ has the form (\ref{3-3-00}).
       $\;\; \square$

 \vskip 0.33 true cm

 \noindent{\bf Remark 4.2.} \  {\it In \cite{Ko}, Kozhevnikov
 proved  \begin{gather*} \label{mm3-3-00}
   \mathbf{A} =\begin{pmatrix} S^{-1}  & {\mathbf{A}}_{JF} & 0\\
    0 &  {\mathbf{A}}_{{}_{FF}} & 0\\
     0 & {\mathbf{A}}_{EF} & 0 \end{pmatrix}. \end{gather*}
     Our Lemma 4.1 is an finer result for the presentation of the operator ${\mathbf{A}}$.}

  \vskip 0.24 true cm

   Since $\mathbf{A}$ is a compact (see proof of Lemma 4.1), self-adjoint operator with respect to the  $[L^2(\Omega)]^n$ inner product, we get that $\mathbf{A}$ has an orthonormal basis of eigenvectors $\{{\mathbf{w}}_k\}_{k=1}^\infty$ corresponding to eigenvalues $\{\tau_k\}_{k=1}^\infty$ which satisfy
 $\tau_k\to 0 \;\; \mbox{as}\;\; k\to +\infty$;  in addition, the number $0$ belongs to the spectrum of $\mathbf{A}$.
From (\ref{3-3-00}) of Lemma 4.1, we get that if $0\ne \tau$ is an eigenvector of $\mathbf{A}$, then exactly one of the following holds:

(i) \  either $\tau$ is an eigenvalue of $S^{-1}$,

(ii) \ or $\tau$ is an eigenvalue of ${\mathbf{A}}_{FF}$.

In fact,  if $\mathbf{f}\ne 0$ is an eigenvector according to eigenvalue $\tau\ne 0$, then ${\mathbf{A}}{\mathbf{f}}=\tau {\mathbf{f}}$. Thus, from (\ref{3-3-00}) we have two possible cases: (i) ${\mathbf{f}}_F=0$ (then $\tau$ and ${\mathbf{f}}_J$ are an eigenvalue and  eigenvector of $S^{-1}$), or ii) ${\mathbf{f}}_{F}\ne 0$ (then $\tau$ and ${\mathbf{f}}_F$ are an eigenvalue and  eigenvector of ${\mathbf{A}}_{FF}$).

\vskip 0.16 true cm

Thus we have showed the following:
\vskip 0.16 true cm

   \noindent  {\bf Lemma 4.3.} \    {\it   Let $ \tau\ne 0$ be an eigenvalue of $\mathbf{A}$. Then
   the space of all eigenvectors of $\mathbf{A}$ corresponding to $\tau$ is just the orthogonal sum of the space of all eigenvectors of $S^{-1}$ and the space of all eigenvectors of $\mathbf{A}_{{}_{FF}}$ corresponding to the same $\tau$.
        }

 \vskip 0.30 true cm

The previous properties allow us to define powers of $\mathbf{A}$.  We define ${\mathbf{A}}^2$ by its action on $\mathbf{u}\in \mathcal{D}(\mathbf{A})$:
\begin{eqnarray*} {\mathbf{A}}^2 \mathbf{u}=\sum_{k=1}^\infty \tau_k^2 \langle\mathbf{u},{\mathbf{w}}_k\rangle {\mathbf{w}}_k,\end{eqnarray*}
where $\langle\cdot, \cdot\rangle$ is the $[L^2(\Omega)]^n$ inner product.
 From the proof of Lemma 4.1, we know that  $\mbox{grad}\,G_2K_{-1} \gamma_0\,\mbox{div}\,\Phi$ is a Green operator (on $\Omega$) of order $-2$.
 Similar to the proof of Lemma 3.5, we can also get that the principal symbol of $\mbox{grad}\,G_2K_{-1} \gamma_0\,\mbox{div}\,\Phi$ is a sum of finitely many terms, each of which contains a factor of the form $\frac{\partial g^{jk}(x)}{\partial x_l}$ (as before, $(g^{jk}(x))$ is the inverse of the Riemannian matric $(g_{jk}(x))$ on $\Omega$).
  It follows from this, (\ref{3-1}), Lemma 3.4 and Lemma 3.5 that for all ${\mathbf{f}}\in [L^2(\Omega)]^n$,
  \begin{eqnarray}\label{4-1} \quad \;\quad \, {\mathbf{A}}^2 {\mathbf{f}}\! &=&\! \mu^{-2} \big\{G_1 [\mathbf{I} -2 \,\mbox{grad}({I} +G_2 K_{-1} \gamma_0 ) \mbox{div}\, G_1]\big\}^2 \mathbf{f}\\
     \!&=&\! \mu^{-2} \big\{ G_1^2 -2 G_1 \,\mbox{grad}(I+ G_2K_{-1}\gamma_0)\mbox{div}\, G_1^2 -2G_1^2 \, \mbox{grad} (I+G_2 K_{-1} \gamma_0)\mbox{div}\, G_1\nonumber\\
    && +4G_1 \,\mbox{grad}(I+G_2K_{-1}\gamma_0)\mbox{div}\, G_1^2\, \mbox{grad}\,(I+G_2K_{-1}\gamma_0)\mbox{div}\, G_1\big\}\mathbf{f}\nonumber \\
    \!&=&\! \mu^{-2} G_1  (\mathbf{I} +{\mathbf{P}}_{-1}+{\mathbf{M}}_{-1}) G_1 \mathbf{f},\nonumber\end{eqnarray}
    where ${\mathbf{P}}_{-1}$ and ${\mathbf{M}}_{-1}$ respectively are pseudodifferential operator of order $-1$ and singular Green operator of order $-1$ (both map $[H^2(\Omega)]^n$ into $[H^3(\Omega)]^n$).
Here we have used the fact that the operators $G_1$ and  $\mbox{grad}\, \mbox{div}$ can commute up to operators of lower order (also $\mbox{div}\, G_1\, \mbox{grad} = I+\Theta_{-1} +\Upsilon_{-1}$, where $\Theta_{-1}\in OPS^{-1}_{1,0}$ and $\Upsilon_{-1}$ is a singular Green operator of order $-1$, see (ii) of Lemma 3.5).
 Of course, the principal symbols of $\mathbf{P}_{-1}$ and $\mathbf{M}_{-1}$ are both a sum of finitely many terms, each of which contains a factor of the form $\frac{\partial g^{jk}(x)}{\partial x_l}$.

 Next, by (\ref{3-3-00}), we have $\mathbf{A}_{{}_{FF}} {\mathbf{f}}_{F}
= P_{{}_F} {\mathbf{A}}\mathbf{f}_{F}= P_{{}_F} \mathbf{A} \, \mbox{grad}\,p$, where $p\in H_0^1(\Omega)$ and $P_{{}_F}=\mbox{grad}\, G_1  \mbox{div}$ is given in (\ref{2-8}).
  It follows from (\ref{3-1}) that
  \begin{eqnarray*} {\mathbf{A}}_{{}_{FF}} \, \mbox{grad}\, p =-\mu^{-1} \mbox{grad}\, G_1 \, \mbox{div}\, G_1\big[ \mathbf{I}- 2\,\mbox{grad}({I} +G_2 K_{-1}\gamma_0)\mbox{div}\, G_1\big] \mbox{grad}\, p, \quad p\in H_0^1 (\Omega). \end{eqnarray*}
  Since $\mbox{div}\, G_1\, \mbox{grad}={I}+{\tilde{L'}}_{-1}+{\tilde{K'}}_{-1}$ (Lemma 3.5 (i)), where ${\tilde L'}_{-1}$ and ${\tilde K'}_{-1}$ respectively are pseudodifferential operator of order $-1$ and singular Green operator of order $-1$,  we have
  \begin{eqnarray} \label{m3-10} &&{\mathbf{A}}_{{}_{FF}}{\mathbf{f}}_{F}= {\mathbf{A}}_{{}_{FF}}\, \mbox{grad}\, p= -\mu^{-1} \mbox{grad}\, G_1[{I}-2I-{\tilde{P}}_{-1}-{\tilde{M}}_{-1}]p\\
  &&\quad \quad  \; \quad =\mu^{-1} \mbox{grad}\, G_1({I}+{\tilde{P}}_{-1}+{\tilde{M}}_{-1})p,\quad p\in H_0^1(\Omega),\nonumber\end{eqnarray}
  where ${\tilde{P}}_{-1}$ and ${\tilde{M}}_{-1}$ respectively are  pseudodifferential operator of order $-1$ and singular Green operators of order $-1$ which map $H_0^1(\Omega)$ to $H^2(\Omega)$.  In addition, the principal symbols of ${\tilde{P}}_{-1}$ and ${\tilde{M}}_{-1}$ are both a sum of finitely many terms, each of which contains a factor of the form $\frac{\partial g^{jk}(x)}{\partial x_l}$.
As pointed out before, the operators $\mathbf{A}$ and $
 \begin{pmatrix} S^{-1}  & 0 \\
    0 &  {\mathbf{A}}_{{}_{FF}} \\
    \end{pmatrix} $
have the same non-zero eigenvalues.
So we may assume ${\mathbf{f}}_{E}\equiv 0$ in the above coordinates representation $({\mathbf{f}}_J, {\mathbf{f}}_F, {\mathbf{f}}_E)^T$.
From now on, we restrict the operators ${\mathbf{A}}$ (or ${\mathbf{A}}^2$) and ${\mathbf{A}}_{{}_{FF}}$ on space $J\oplus F$ and still denote them by ${\mathbf{A}}$ (or ${\mathbf{A}}^2$) and $\mathbf{A}_{{}_{FF}}$.

\vskip 0.1 true cm

  Furthermore, let $ \mathbf{A}_{{}_{FF}} \mathbf{f}=\lambda\,\mathbf{f}$, where $\mathbf{f}\in F$ and $\lambda\ne 0$. Since ${\mathbf{f}}=\mbox{grad}\, p$, $\,p\in H_0^1(\Omega)$, we find by (\ref{m3-10}) that $ \mu^{-1} \mbox{grad}\, G_1({I}+{\tilde{P}}_{-1}+ {\tilde{M}}_{-1})p= \lambda\,\mbox{grad}\, p$. Applying the operator div and then $G_1$ to this equality, we get
  \begin{eqnarray} \label{m3-1}  \mu^{-1} G_1 ({I}+{\tilde{P}}_{-1}+ {\tilde{M}}_{-1})p=\lambda\, p.\end{eqnarray} Conversely, applying the operator grad to (\ref{m3-1}), we find by (\ref{m3-10}) that $ \mathbf{A}_{{}_{FF}} {\mathbf{f}} =\lambda\,{\mathbf{f}}$ with ${\mathbf{f}}=\mbox{grad}\, p$. Therefore, the number $\lambda\ne 0$ and the vector ${\mathbf{f}}=\mbox{grad}\, p \in F$ are an eigenvalue and corresponding eigenvector of the operator ${\mathbf{A}}_{{}_{FF}}$ if and only if the pair $(\lambda, p)$ constitutes an eigenvalue and eigenvector of the operator ${\tilde{A}}_{FF}$, where
  \begin{eqnarray} \label{4-3}
  {\tilde A}_{{}_{FF}}p := \mu^{-1}  G_1 (I+{\tilde{P}}_{-1}+ {\tilde{M}}_{-1})p,\quad \quad p\in H_0^1(\Omega),\end{eqnarray}
  and ${\tilde{P}}_{-1}$ and ${\tilde{M}}_{-1}$ are the same operators as in (\ref{m3-10}). Clearly, ${\tilde{A}}_{{}_{FF}} : H_0^1(\Omega)\to H^3(\Omega)$.

\vskip 0.26 true cm

  \noindent  {\bf Lemma 4.4.} \    {\it  The kernel spaces  of the operators ${\mathbf{A}}^2$ and ${\tilde{A}}_{{}_{FF}}$ are finite dimensional, and   \begin{eqnarray}\label{aa-3-1}\mbox{dim}\,(\mbox{ker}\, {\mathbf{A}}^2)= \mbox{dim}\, (\mbox{ker}\, \mathbf{A}_{FF})=\mbox{dim}\,(\mbox{ker}\, {\tilde{A}}_{{}_{FF}}).\end{eqnarray}}

  \vskip 0.10 true cm

  \noindent  {\it Proof.} \  Since the compact operators $G_1^{2}$ and $G_1$ are invertible, it follows from (\ref{4-1}) and (\ref{4-3}) that the kernel spaces of ${\mathbf{A}}^2$ and ${\tilde A}_{{}_{FF}}$ are finite-dimensional.
 If ${\mathbf{A}}_{{}_{FF}} {\mathbf{f}}_{{}_F}=0$, where $0\ne {\mathbf{f}}_{{}_F}= \mbox{grad}\, p\in F$, $\, p\in H_0^1(\Omega)$, then, by (\ref{m3-10})  we have \begin{eqnarray*}  \mbox{grad}\, G_1 ({I}+{\tilde{P}}_{-1}+{\tilde{M}}_{-1}) p =0 \quad \, \mbox{in}\,\, \Omega,\end{eqnarray*}
so that \begin{eqnarray}\label {b-4-1} G_1 (I+{\tilde{P}}_{-1}+{\tilde{M}}_{-1})p \equiv \mbox{const} \quad \,\, \mbox{in}\;\; \Omega.\end{eqnarray}
According to the definition of $G_1$ (see (\ref{2-6})), we get
  $(G_1({I}+{\tilde{P}}_{-1}+{\tilde{M}}_{-1})p)\big|_{\partial \Omega}=0$. It follows from (\ref{b-4-1}) that $G_1({I}+{\tilde{P}}_{-1}+{\tilde{M}}_{-1})p\equiv 0$ in $\Omega$,
i.e., ${\tilde{A}}_{{}_{FF}} p=0$ in $\Omega$. Conversely, if $0\ne p\in H_0^1(\Omega)$ and ${\tilde A}_{{}_{FF}}
p =0$, i.e., $\mu^{-1} G_1 (I+{\tilde{P}}_{-1}+{\tilde{M}}_{-1})p=0$ in $\Omega$, then \begin{eqnarray*} {\mathbf{A}}_{{}_{FF}} {\mathbf{f}}_{{}_F}=\mu^{-1} \mbox{grad}\, G_1(I+{\tilde{P}}_{-1}+{\tilde{M}}_{-1})p \equiv 0\quad \, \mbox{in}\,\, \Omega \end{eqnarray*}
with $\mathbf{f}_F=\mbox{grad}\, p$ (see (\ref{m3-10}).
 Clearly, $\mbox{grad}\, p\ne 0$. Thus \begin{eqnarray}  \label{c-3-10} \mbox{dim} (\mbox{ker}\, {\tilde{A}}_{{}_{FF}})=
 \mbox{dim}\, (\mbox{ker}\, {\mathbf{A}}_{{}_{FF}}).\end{eqnarray}
 We denote by $m_0$ the dimension of the above kernel space.

  Obviously, $\mbox{ker}\, {\mathbf{A}}_{{}_{FF}} \subset \mbox{ker}\, {\mathbf{A}}_{{}_{FF}}^2$. Now, let $\mathbf{f}\in \mbox{ker}\, {\mathbf{A}}_{{}_{FF}}^2$, i.e., ${\mathbf{A}}_{{}_{FF}}^2 \mathbf{f}=0$. We claim that ${\mathbf{A}}_{{}_{FF}}\mathbf{f}=0$ for such an $\mathbf{f}\in \mbox{ker}\, {\mathbf{A}}_{{}_{FF}}^2$. Suppose by contradiction that ${\mathbf{A}}_{{}_{FF}} \mathbf{f}\ne 0$. Let $\{\mbox{grad}\, p_j\}_{j=1}^\infty$  ($p_j\in H_0^1(\Omega)$) are orthonormal eigenvectors of ${\mathbf{A}}_{{}_{FF}}$  corresponding to all non-zero eigenvalues $\{\alpha_j\}_{j=1}^\infty$, and let $\{{\mathbf{r}}_j\}_{j=1}^{m_0}$ is an orthonormal basis of $\mbox{ker}\, {\mathbf{A}}_{{}_{FF}}$.  Let $\mathbf{f}=\sum_{j=1}^\infty \beta_j (\mbox{grad}\; p_j)+ \sum_{j=1}^{m_0}  {\tilde{\beta}}_j {\mathbf{r}}_j$.  Then \begin{eqnarray*} {\mathbf{A}}_{{}_{FF}}\mathbf{f}=\sum_{j=1}^\infty \alpha_j\beta_j (\mbox{grad}\; p_j)\ne 0,\end{eqnarray*} so that $(\alpha_1\beta_1 , \cdots, \alpha_j\beta_j, \cdots)\ne (0, \cdots, 0,\cdots)$. Furthermore, \begin{eqnarray*} {\mathbf{A}}^2_{{}_{FF}}\mathbf{f}=\sum_{j=1}^\infty \alpha_j^2\beta_j (\mbox{grad}\; p_j)\ne 0\end{eqnarray*} since $\{\mbox{grad}\, p_j\}$ is an orthonormal system. This contradicts the assumption ${\mathbf{A}}_{{}_{FF}}^2 \mathbf{f}=0$. Hence  the assertion $\mbox{ker}\, {\mathbf{A}}_{{}_{FF}}^2 \subset \mbox{ker}\,{\mathbf{A}}_{{}_{FF}}$ holds, so we have
   \begin{eqnarray} \label{c-3-12} \mbox{ker}\, {\mathbf{A}}_{{}_{FF}}^2 =\mbox{ker}\, {\mathbf{A}}_{{}_{FF}}.\end{eqnarray}
    Finally, since the operator $S^{-1}$ is invertible in $J$, we have that $\mbox{dim}(\mbox{ker}\, {\mathbf{A}}^2)= \mbox{dim}(\mbox{ker}\, {\mathbf{A}}_{{}_{FF}}^2)$, so by
     (\ref{c-3-10}) and (\ref{c-3-12}), we obtain $\mbox{dim}(\mbox{ker}\, {\mathbf{A}}^2)= \mbox{dim}(\mbox{ker}\, \mathbf{A}_{FF})=\mbox{dim}(\mbox{ker}\, {\tilde{A}}_{{}_{FF}})$. $\;\;\square$

\vskip 1.29 true cm

\section{Asymptotic expansion of trace of strongly continuous semigroup }

\vskip 0.45 true cm

 \noindent  {\it Proof of Theorem 1.1.} \  Since $\mbox{dim}(\mbox{ker}\, {\mathbf{A}}^2)=\mbox{dim}(\mbox{ker}\, {\tilde{A}}_{{}_{FF}})=m_0$ which is a finite number,  we may adjust the Jordan matrices  ${\mathbf{A}}^2$ and ${\tilde A}_{{}_{FF}}$ on the kernel spaces of these operators by replacing the zero eigenvalue by a common constant, saying $\varrho>0$. In other words, we add to ${\mathbf{A}}^2$ and ${\tilde A}_{{}_{FF}}$ some finite-dimensional operators of order $-\infty$ since their kernels consist of infinitely smooth functions.  Let us denote the finite-dimensional operators  added by $\mathbf{R}$  and ${\tilde R}_F$.

By inverting the operators ${\mathbf{A}}^2+\mathbf{R}$ and ${\tilde A}_{{}_{FF}}+{\tilde{R}_F}$
 and by comparing their symbols from (\ref{4-1}) and (\ref{4-3}) (or  by composition formula (\ref{17-9-12.1})  for pseudodifferential operators and the Wiener-Hopf algebra (\ref{2020.11.12-4})---(\ref{2020.11.12-6}) for singular Green symbols), we find that  \begin{eqnarray}\label{4-4}  &&({\mathbf{A}}^2+ \mathbf{R})^{-1} = \mu^2 G_1^{-2} \,\mathbf{I}+\mathbf{B}'_1+\mathbf{C}'_1,\\
\label{4-5}&& ({\tilde {A}}_{{}_{FF}} +{\tilde{R}_F})^{-1} =\mu G_1^{-1} +B'_2+C'_2,\end{eqnarray}
 where $\mathbf{B}'_1$ and $B'_2$ respectively are pseudodifferential operators of order $3$ and $1$, while ${\mathbf{C}'}_1$ and $C'_2$ respectively are singular Green operators of order $3$ and $1$. Clearly, the principal symbols of $\mathbf{B}'_{1}+\mathbf{C}'_{1}$ and $B'_2+C'_2$ are both the sums of finitely many terms, each of which contains a factor of the form $\frac{\partial g^{jk}(x)}{\partial x_l}$ (again $(g^{jk}(x))$ denotes the inverse of Riemannian matric $(g_{jk}(x))$ on $\Omega$).
Recall that $G_2 \gamma_0$ is a Green operator of order $-1$ and $G_1= \Phi -G_2 \gamma_0 \Phi = (I- G_2 \gamma_0) \Phi $. Then (\ref{4-4}) and (\ref{4-5}) can be rewritten as
\begin{eqnarray}\label{2020.9.16-4}  && ({\mathbf{A}}^2+ \mathbf{R})^{-1} = \mu^2\Phi^{-2} \mathbf{I}  -\mathbf{B}''_1 -\mathbf{C}''_1= \mu^2 \Delta^2\,\mathbf{I} -\mathbf{B}_1 -\mathbf{C}_1,\\
&& \label{2020.10.3-1}({\tilde {A}}_{{}_{FF}} +{\tilde{R}_F})^{-1} =\mu \Phi^{-1} +B''_2+C''_2= \mu \Delta +B_2+C_2,\end{eqnarray}
 where $\mathbf{B}_1$ and $B_2$ respectively are pseudodifferential operators of order $3$ and $1$, while $\mathbf{C}_1$ and $C_2$ respectively are singular Green operators of order $3$ and $1$. In addition, the principal symbols of $\mathbf{B}_{1}+\mathbf{C}_{1}$ and $B_2+C_2$ are both the sums of finitely many terms, each of which contains a factor of the form $\frac{\partial g^{jk}(x)}{\partial x_l}$.

 \vskip 0.1 true cm

The proof of (\ref{1-7}) is broken up into a number of steps.

\vskip 0.1 true cm

Step 1. \  We first calculate the asymptotic expansion for the integral of the trace of the strongly continuous semigroup $e^{-t(-\tilde{A}_{FF} -\tilde{R}_F)}$ as $t\to 0^+$.

Clearly, \begin{eqnarray} \label{2020.10.15-5} (-\mu \Delta - B_2 - C_2-\lambda)^{-1}= (-\mu \Delta -\lambda)^{-1} \Big[ I -(B_2 +C_2) (- \mu \Delta -\lambda)^{-1} \Big]^{-1}. \end{eqnarray}
 Putting $\Big[ I -(B_2 +C_2) ( -\mu \Delta -\lambda)^{-1} \Big]^{-1} = I+ L + \tilde{L}$, where $L$ (respectively, $\tilde{L}$) is a pseudodifferential (respectively, singular Green) operator of order $-1$, we have
 \begin{eqnarray*} \big[ I -(B_2 +C_2) (-\mu \Delta-\lambda)^{-1}\big] \big[I+L+\tilde{L}\big] =I, \end{eqnarray*}
 which implies $L+\tilde{L}= \big[ I- (B_2+C_2) (-\mu\Delta-\lambda)^{-1} \big]^{-1} (B_2+C_2) (-\mu\Delta -\lambda)^{-1}$. Thus
 \begin{eqnarray} \label{2020.10.15-1}  \;\;\;\quad\;\;\;\;\;  \big[ I \!-\!(B_2 \!+\!C_2) (\!-\mu \Delta-\lambda)^{-1}\big]^{-1} \! = \!I\!+\! \big[ I\!-\! (B_2\!+\!C_2) (\!-\mu\Delta-\lambda)^{-1} \big]^{-1} (B_2\!+\!C_2) (\!-\mu \Delta-\lambda)^{-1}.\end{eqnarray}
 By (\ref{2020.10.15-5}) and iterative formula (\ref{2020.10.15-1}), we get
 \begin{eqnarray} \label{2020.10.15-2}   (-\mu \Delta - B_2 - C_2-\lambda)^{-1}= (-\mu\Delta -\lambda)^{-1} \sum_{k=0}^\infty
 \big[ (B_2 +C_2) (-\mu\Delta -\lambda)^{-1}\big]^{k}.\end{eqnarray}
 This leads to
\begin{eqnarray*}  e^{-t(-\tilde{A}_{FF} -\tilde{R}_F)}f(x) \!\!\!\!\!\!\!&\!&\!\!= e^{-t (-\mu \Delta - B_2 -C_2)} f(x) \\
\!\!\!\!\!\!\!&\!&\!\!=\frac{i}{2\pi} \int_{\mathcal{C}} e^{-t\lambda} (-\mu \Delta - B_2 - C_2-\lambda)^{-1} f(x) \,d\lambda \\
\!\!\!\!\!\!\!&\!&\!\!= \frac{i}{2\pi} \int_{\mathcal{C}} e^{-t\lambda} (-\mu\Delta -\lambda)^{-1} \sum_{k=0}^\infty  \big[ (B_2+C_2) (-\mu \Delta -\lambda)^{-1} \big]^{k} f(x)\,d\lambda.\end{eqnarray*}

Let $x=(x'; x_n)$ be local coordinates for $\Omega$ near $\partial \Omega$.
 If $\mathfrak{E}$ is a local frame on $\partial \Omega$; extend $\mathfrak{E}$ to an $n$-dimensional local frame in a neighborhood of $\partial \Omega$ by parallel transport along the geodesic normal rays (see, section 3 or p.$\,$1101 of \cite{LU}).
  Let $\mathcal{M}=\Omega \cup (\partial \Omega)\cup \Omega^*$ be the (closed) double of $\Omega$, and $Q$ the double to $\mathcal{M}$ of
 the pseudodifferential operator $-\mu \Delta  -B_2-C_2$. The symbol (or the corresponding coefficients) of $-\mu \Delta -B_2-C_2$  on $\mathcal{M}$ jumps as $x$ crosses $\partial \Omega$, but $\frac{\partial u}{\partial t}=Qu$ still has a nice fundamental solution $K_2(t, x,y)$ of class $C^\infty [(0,\infty)\times (\mathcal{M} \setminus\partial \Omega)\times (\mathcal{M} \setminus\partial \Omega) ]
 \cap C^1 ((0,\infty)\times \mathcal{M}\times \mathcal{M})$, approximable even on $\partial \Omega$. Define $Q^{-}$ to be $Q\big|_{C^\infty(\bar \Omega)}$ subject to $u=0$ on $\partial \Omega$. Then the fundamental solution $K_2^{-}(t,x,y)$
 of $\frac{\partial u}{\partial t}=Q^{-}u$ can be expresses on $(0,\infty)\times \Omega\times \Omega$ as
 \begin{eqnarray} \label{c4-23} K_2^{-} (t,x,y) =K_2(t,x,y)- K_2(t,x,\overset{\ast} {y}),\end{eqnarray}
 $\overset{*} {y}$ being the double of $y\in \Omega$ (see, p.$\,$53 of \cite{MS}).
  Since the strongly continuous semigroup $(e^{-tQ})_{t\ge 0}$ can also be represented as  \begin{eqnarray*} e^{-t{Q}}f(x) =\frac{i}{2\pi } \int_{\mathcal{C}} e^{-t\lambda} ({Q}-\lambda I)^{-1} f(x)\, d\tau,\end{eqnarray*}
where $\mathcal{C}$ is a suitable curve in the complex plane in the positive direction around the spectrum of ${Q}$ (i.e., a contour around the positive real axis). Let $(q(x, \xi, \lambda); \psi(x', \xi',\xi_n,\eta_n,\lambda))$ is the full symbol of the Green operator $({Q}-\lambda I)^{-1}$.
Taking into account the fact that the Fourier transform of the Dirac delta function $\delta$ is $1$, we find that,
for any $x, y\in \mathcal{M}$,
\begin{align*} &K_2(t,x,y)= e^{-tQ}\delta(x-y)=\frac{i}{2\pi} \!\int_{\mathcal{C}}\! e^{-t\lambda}  ({Q}-\lambda I)^{-1} \delta(x- {y})d\lambda\\
& \qquad \quad \quad\;\, = \frac{i}{2\pi} \!\int_{\mathcal{C}}\! e^{-t\lambda}  \bigg\{ \frac{1}{(2\pi)^n} \int_{\mathbb{R}^n} e^{i(x-y)\cdot \xi} q (x, \xi, \lambda) (\hat{\delta}(\xi)) \, d\xi \\
& \qquad \qquad \quad \;\, +\! \frac{ 1}{(2\pi)^{n+1}}\! \int_{\mathbb{R}^{n\!-\!1}} \!\!d\xi' \int^{+}\! \!e^{i(x-y)\cdot \xi}d\xi_n\!
\! \int^+\!\!\psi(x', \xi',\xi_n,\eta_n,\lambda)(\hat{\delta} (\xi', \eta_n))\, d\eta_n\bigg\}
d\lambda\\
&\qquad \quad \;\,\quad =  \frac{i}{2\pi} \!\int_{\mathcal{C}}\! e^{-t\lambda}  \bigg\{ \frac{1}{(2\pi)^n} \int_{\mathbb{R}^n} e^{i(x-y)\cdot \xi} q (x, \xi, \lambda) \, d\xi \\
& \qquad \qquad \quad \;\, +\! \frac{ 1}{(2\pi)^{n+1}}\! \int_{\mathbb{R}^{n\!-\!1}} \!\!d\xi' \int^{+}\! \!e^{i(x-y)\cdot \xi}d\xi_n\!
\! \int^+\!\!\psi(x', \xi',\xi_n,\eta_n,\lambda)\, d\eta_n\bigg\}
d\lambda.\end{align*}
As pointed out before, the right-hand side of above equality is just the representation of the Schwartz kernel of the (Green) operator semigroup $(e^{-tQ})_{t\ge 0}$ (see section section 3, or (18.1.7) on p.$\,$69 in \cite{Ho3} or \cite{Gr}).
 Thus, we have the following Levi's sum for the elementary solution (cf. \S3 of \cite{MS}) \begin{eqnarray} \label{2020.10.16-6}  {{K}_2} (t,x,y)\!\!\!&&\!\!\!\!\!\!= e^{-tQ}\delta(x-y) = e^{-t(-\mu \Delta -B_2-C_2)}\delta(x-y)\\
\!\!\!\!&&\!\!\!\!\!\!= \frac{i}{2\pi} \int_{\mathcal{C}} e^{-t\lambda} (-\mu\Delta -\lambda)^{-1} \sum_{k=0}^\infty  \big[ (B_2+C_2) (-\mu \Delta -\lambda)^{-1} \big]^{k} \delta(x-y)\,d\lambda, \;\;\, \forall x, y\in \mathcal{M}.\nonumber\end{eqnarray}
In particular, for any $x\in \Omega$,
 \begin{eqnarray} \label{2020.10.17-1} \begin{aligned} {{K}_2} (t,x,x)&= \frac{i}{2\pi} \int_{\mathcal{C}} e^{-t\lambda} (-\mu\Delta -\lambda)^{-1} \sum_{k=0}^\infty  \big[ (B_2+C_2) (-\mu \Delta -\lambda)^{-1} \big]^{k} \delta(x-x)\,d\lambda, \\
{{K}_2} (t,x,\overset{\ast} {x})&= \frac{i}{2\pi} \int_{\mathcal{C}} e^{-t\lambda} (-\mu\Delta -\lambda)^{-1} \sum_{k=0}^\infty  \big[ (B_2+C_2) (-\mu \Delta -\lambda)^{-1} \big]^{k} \delta(x-\overset{\ast} {x})\,d\lambda.\end{aligned} \end{eqnarray}
For given (small enough) $\epsilon>0$, denote by $U_\epsilon(\partial \Omega)=\{z\in {\mathcal{M}}\big| \mbox{dist}\, (z, \partial \Omega)<\epsilon\}$ the $\epsilon$-neighborhood of $\partial \Omega$ in $\mathcal{M}$.
          It suffices for the proof of Theorem 1.1 (i.e., (\ref{1-7})) to check that
            \begin{align} \label{2020.10.17-3}&  \int_{\Omega'} \Big(\frac{i}{2\pi} \int_{\mathcal{C}} e^{-t\lambda} (-\mu\Delta -\lambda)^{-1}  \delta(x-x)\,d\lambda \Big)dx= \frac{|\Omega'|}{(4\pi \mu t)^{n/2}} +O(t^{1-\frac{n}{2}})\,\,  \mbox{as}\,\, t\!\to \!0^+, \,\, \forall \,\,\Omega'\subset \Omega,\end{align}

 \begin{align} \label{2020.10.17-7} &\quad  \int_{\Omega''}\! \Big(\frac{i}{2\pi} \!\int_{\mathcal{C}}\!\! e^{-t\lambda} (-\mu\Delta\! -\!\lambda)^{-1}  \delta(x\!-\!\overset{\ast} {x})\,d\lambda \Big)dx\!=\!O(t^{1\!-\!\frac{n}{2}}\!)\, \, \mbox{as}\,\,t\!\to \!0^+\!,\,\,
  \forall \, \Omega''\subset \!\Omega \setminus U_{\epsilon}(\partial \Omega),\qquad \qquad\qquad
  \end{align}

 \begin{align} \label{2020.10.25-1} &\; \;\; \;\int_{W}\!\!\Big(\!\frac{i}{2\pi} \!\!\int_{\mathcal{C}}\!\! e^{\!-\!t\lambda} (\!-\mu\Delta\! \!-\!\!\lambda\!)^{\!-\!1} \! \delta(x\!\!-\!\overset{\ast} {x})d\lambda \!\Big)\!dx\!=\! \frac{|W \cap \partial \Omega|}{4(\!4\pi \mu t)^{\!\frac{n\!-\!1}{2}}} \!+\!O(\!t^{\!1\!-\!\frac{n}{2}}\!)\, \, \mbox{as}\,\,t\!\to \!0^+\!,\,\,
  \forall \, W\!\!\subset \!\!U_{\!\epsilon}(\!\partial \Omega\!)\,\mbox{and}\, W\!\cap\! \partial \Omega \!\ne \!\emptyset,
  \end{align}

           \begin{align}   \label{2020.10.17-4}& \;\;\;\; \int_{\Omega'}\! \!\Big(\frac{i}{2\pi}\!\! \int_{\mathcal{C}} \!e^{-t\lambda}
            (\!-\mu\Delta\! -\!\lambda)^{\!-\!1} \!\sum_{k=1}^\infty  \!\big[ (B_2\!+\!C_2) (\!-\mu \Delta\! -\!\lambda)^{\!-\!1}\! \big]^k \delta(x\!-\!x)d\lambda\!\Big)dx\!=\!
           O(t^{1\!-\!\frac{n}{2}}) \,\,   \mbox{as}\,\, t\!\to\! 0^+\!, \,\, \forall \,\Omega'\!\subset\! \Omega,\end{align}

                \begin{align}  \label{2020.10.17-40} \;\;\;\; \int_{\Omega'}\!\! \Big(\frac{i}{2\pi} \!\int_{\mathcal{C}}\! e^{\!-\!t\lambda}
           (-\mu\Delta \!-\!\lambda)^{\!-\!1} \!\sum_{k=1}^\infty \! \big[ (\!B_2\!+\!C_2) (-\mu \Delta\! -\!\lambda)^{\!-\!1} \!\big]^k \! \delta(\!x\!-\!\overset{\ast} {x})d\lambda\Big)dx\!=\!
           O(t^{1\!-\!\frac{n}{2}}\!) \;   \mbox{as}\; t\!\to\! 0^+\!, \;\forall \Omega'\!\subset\! \Omega.\end{align}

\vskip 0.2 true cm

(a) \ \    We first consider   $\int_{\Omega'} \Big(\frac{i}{2\pi} \int_{\mathcal{C}} e^{-t\lambda} (-\mu\Delta -\lambda)^{-1}  \delta(x-x)\,d\lambda \Big)dx$
for any $\Omega'\subset \Omega$. Let $p(x, \xi, \lambda)\!\!:=\sum_{l\ge 0} p_{-2-l} (x,\xi,\lambda)\in S^{-2}_{1,0}$  be the full symbols of the pseudodifferential operators $(-\mu\Delta- \lambda)^{-1}$ in local coordinates.
   Then (see section 2)
 \begin{eqnarray*} p_{-2} =(\mu \sum_{j=1}^n \xi_j^2 -\lambda)^{-1},  \quad p_{-2-1} = \vartheta_{1,1} (x, \xi) p_{-2}^2,\quad p_{-2 -l} = \sum\limits_{k=1}^{2l} \vartheta_{l,k} (x, \xi)  p_{-2}^{k+1}, \cdots, \quad l> 1,\end{eqnarray*}
where $p_{-2-l}(x, \xi,\lambda)$ are homogeneous in $\xi$ of degree
$-2-l$ for $|\xi|>1$, and $\vartheta_{l,k}$ independent of $\lambda$ and homogeneous of degree $2k-l$ in $\xi$ for $|\xi|\ge 1$. In particular, each term of $\vartheta_{1,1} (x,\xi)$  has a factor of the form $\frac{\partial g^{jk}(x)}{\partial x_l}$, where $(g_{jk})$ is the Riemannian metric on $\Omega$.
               It follows that \begin{eqnarray*}\frac{i}{2\pi} \!\int_{\mathcal{C}} e^{-t\lambda} (-\mu\Delta -\lambda)^{-1} \delta(x-y)\,d\lambda=\frac{1}{(2\pi)^n} \int_{\mathbb{R}^n} e^{i(x-y)\cdot\xi} \bigg(\!\frac{i}{2\pi } \! \int_{\mathcal{C}} e^{-t\lambda} \sum_{l\ge 0} p_{-2-l}(x, \xi, \lambda) d\lambda\!\bigg) d\xi, \;\; \forall x,y\in \mathcal{M},\end{eqnarray*}
so \begin{eqnarray*}&&\frac{i}{2\pi} \int_{\mathcal{C}} e^{-t\lambda} (-\mu\Delta -\lambda)^{-1} \delta(x-x)\,d\lambda=\frac{1}{(2\pi)^n} \int_{\mathbb{R}^n} \! \bigg(\!\frac{i}{2\pi } \! \int_{\mathcal{C}} e^{-t\lambda} \sum_{l\ge 0} p_{-2-l}(x, \xi, \lambda) d\lambda\!\bigg) d\xi, \;\; \forall x\in \Omega,\\
   && \frac{i}{2\pi} \int_{\mathcal{C}} e^{-t\lambda} (-\mu\Delta -\lambda)^{-1} \delta(x-\overset{\ast} {x})\,d\lambda=\frac{1}{(2\pi)^n} \int_{\mathbb{R}^n} e^{i(x-\overset{\ast} {x})\cdot\xi} \bigg(\!\frac{i}{2\pi } \! \int_{\mathcal{C}} e^{-t\lambda} \sum_{l\ge 0} p_{-2-l}(x, \xi, \lambda) d\lambda\!\bigg) d\xi,  \;\; \forall x\in \Omega.
  \end{eqnarray*}
  For each fixed point $x$ in Riemannian manifold $(\bar \Omega, g)$, we may also take a geodesic normal coordinate system centered at this $x$
       such that (see p.$\,$555 of \cite{Ta2}) \begin{eqnarray} \label{18/7/14/1} g_{jk}(x)= \delta_{jk}, \; \; \frac{\partial g_{jk}}{\partial x_l}(x)
 =0  \;\;  \mbox{for all} \;\; 1\le j,k,l \le n. \end{eqnarray}
  Hence
\begin{eqnarray*}   v_{-2}(t,x, \xi) =e^{-t\mu \sum_{j=1}^n \xi_j^2}, \;\; v_{-3} (t, x, \xi)= t\, \vartheta_{1,1} (x, \xi)e^{-t\mu \sum_{j=1}^n \xi_j^2}, \;
           \end{eqnarray*}  \begin{eqnarray*}v_{-2-l} (t, x, \xi)= \sum_{k=1}^{2l}\frac{t^k}{k!} \vartheta_{l,k} (x,\xi) e^{-t\mu \sum_{j=1}^n \xi_j^2}, \quad l\ge 2,\end{eqnarray*}
   where $ v_{-2-l} (t, x, \xi)=\frac{i}{2\pi} \int_{\mathcal{C}}e^{-t \lambda}q_{-2-l} (x, \xi, \lambda)  d\lambda$, and $\sum_{l\ge 0} v_{-2-l} (t, x, \xi)$ is the full symbol of the operator $\frac{i}{2\pi} \int_{\mathcal{C}}e^{-t \lambda} \left(-\mu \Delta -\lambda\right)^{-1} d\lambda$.
 In particular,  $\vartheta_{1,1} (x, \xi) =0$ since $\frac{\partial g^{jk}(x)}{\partial x_l} =0$ by (\ref{18/7/14/1}).
This leads to
   \begin{eqnarray} \label{17-11-21}\frac{i}{2\pi} \int_{\mathcal{C}} e^{-t \lambda } (-\mu\Delta- \lambda)^{-1} \delta(x-x)
\!\!\!\!\!\!\!\! \!&& \! =\frac{1}{(2\pi)^n}  \int_{\mathbb{R}^n} e^{i(x-x)\cdot \xi}\bigg[e^{-t\mu \sum_{j=1}^n \xi_j^2} +t\vartheta_{1,1}(x,\xi)e^{-t\mu \sum_{j=1}^n \xi_j^2}\nonumber\\
 && \;\,\;+
  \sum_{l=2}^\infty \bigg(\sum_{k=1}^{2l}
    \frac{t^k}{k!}\big(\vartheta_{l,k}(x, \xi)\big) e^{-t\mu\sum_{j=1}^n \xi_j^2}\bigg)\bigg] d\xi\nonumber\\
\!\!\!\!\!\!\!\! && = \frac{1}{(2\pi)^n}  \int_{\mathbb{R}^n} \bigg[e^{-t\mu \sum_{j=1}^n \xi_j^2} +  \sum_{l=2}^\infty \bigg(\sum_{k=1}^{2l}
    \frac{t^k}{k!}\big(\vartheta_{l,k}(x, \xi)\big) e^{-t\mu\sum_{j=1}^n \xi_j^2}\bigg)\bigg] d\xi\nonumber\end{eqnarray}
    and   \begin{eqnarray} \label{2020.10.25-02}&&\frac{i}{2\pi} \int_{\mathcal{C}} e^{-t \lambda } (-\mu\Delta- \lambda)^{-1} \delta(x-\overset{\ast} {x})
 = \frac{1}{(2\pi)^n}  \int_{\mathbb{R}^n}e^{i(x-\overset{\ast} {x})\cdot \xi} \bigg[e^{-t\mu \sum_{j=1}^n \xi_j^2} +
 \nonumber\\
&& \quad\quad +  \sum_{l=2}^\infty \bigg(\sum_{k=1}^{2l}
    \frac{t^k}{k!}\big(\vartheta_{l,k}(x, \xi)\big) e^{-t\mu\sum_{j=1}^n \xi_j^2}\bigg)\bigg] d\xi.\nonumber\end{eqnarray}
            For any $\Omega'\subset \Omega$ and $\Omega''\subset \Omega\setminus U_\epsilon (\partial \Omega)$, we have that
  \begin{eqnarray*}  \!\!\!\! &&\!\!\!\! \int_{\Omega'}\! \left\{  \!\frac{1}{(2\pi)^n}\!\int_{\mathbb{R}^n}\!e^{-t\mu\sum_{j=1}^n \xi_j^2} d\xi\! \right\} dx\!=\!\frac{1}{(4\pi \mu t)^{n/2}}\!\int_{\Omega'}\!dx \!= \!\frac{|\Omega'|}{(4\pi \mu t)^{n/2}},\\
  \!\!\!\!&& \!\!\!\! \int_{\Omega''} \left\{  \frac{1}{(2\pi)^n}\int_{\mathbb{R}^n} e^{i(x-\overset{\ast}{x})\cdot \xi} e^{-t\mu\sum_{j=1}^n \xi_j^2} d\xi \right\} dx=\frac{1}{(4\pi \mu t)^{n/2}}  \int_{\Omega''} e^{-|x-\overset{\ast} {x}|^2/4\mu t}dx=O(t^{1-\frac{n}{2}})\;\;\mbox{as}\;\, t\to 0\end{eqnarray*}
    because of $|x-\overset{\ast}{x}|\ge 2\epsilon$.  Note that for a  symbol $\varphi (x,\xi)$ with homogeneous of degree $m$ in $\xi$, by variable substitution $\tilde{\xi}=\sqrt{\mu t} \, \xi$ one has
     \begin{eqnarray} \label{2020.11.15-7} && \bigg|\int_{{\mathbb{R}}^n} e^{i(x-y)\cdot \xi}\big(\varphi (x, \xi)\big) e^{-t\mu\sum_{j=1}^n \xi_j^2}d\xi\bigg|
   =  \bigg|(\mu t)^{-\frac{n+m}{2}}\int_{\mathbb{R}^n} e^{\frac{i(x-y)}{\sqrt{\mu t}}\cdot \tilde \xi} \phi(x, \tilde \xi)e^{-\tilde \xi^2} d\tilde \xi\bigg| \\
 && \qquad \qquad \quad  \le  (\mu t)^{-\frac{n+m}{2}} \int_{\mathbb{R}^n} \big| \phi(x, \tilde \xi)\big|e^{-\tilde \xi^2} d\tilde \xi\le c\,(\mu t)^{-\frac{n+m}{2}} \;\;\mbox{for any}\,\; x\in \Omega'\subset \Omega, \;y\in \mathcal{M}.\nonumber
   \end{eqnarray}
     From (\ref{2020.11.15-7}) we get that
         \begin{eqnarray*} \label{15-1}
   &&\int_{\Omega'} \bigg\{ \frac{1}{(2\pi)^n}\int_{\mathbb{R}^n} \bigg[   \sum_{l=2}^\infty \bigg(\sum_{k=1}^{2l}
    \frac{t^k}{k!}\big(\vartheta_{l,k}(x, \xi)\big) e^{-t\mu\sum_{j=1}^n \xi_j^2}\bigg)\bigg] d\xi\bigg\}dx =O(t^{1-n/2}) \;\;   \mbox{as}\,\, t\to 0^+,\nonumber\end{eqnarray*}
    and
    \begin{eqnarray}  \label{15-001}
   &&\int_{\Omega'} \bigg\{ \frac{1}{(2\pi)^n}\int_{\mathbb{R}^n} e^{i(x-\overset{\ast}{x})\cdot \xi}\bigg[  \sum_{l=2}^\infty \bigg(\sum_{k=1}^{2l}
    \frac{t^k}{k!}\big(\vartheta_{l,k}(x, \xi)\big) e^{-t\mu\sum_{j=1}^n \xi_j^2}\bigg)\bigg] d\xi\bigg\}dx =O(t^{1-n/2}) \;\;   \mbox{as}\,\, t\to 0^+.\nonumber\end{eqnarray}
Hence
   \begin{eqnarray*} \label{a-3-1-1} &&  \int_{\Omega'} \bigg[\frac{i}{2\pi} \int_{\mathcal{C}} e^{-t\lambda} (-\mu \Delta  -\lambda)^{-1} \delta (x-x) \bigg]dx = \frac{|\Omega'|}{(4\pi \mu t)^{n/2}} +O(t^{1-n/2}) \quad  \mbox{as}\;\, t\to 0^+\;\;\mbox{for any}\,\, \Omega'\subset \Omega,\nonumber\end{eqnarray*}
       \begin{eqnarray*} \label{a-3-1-002} &&  \int_{\Omega'} \bigg[\frac{i}{2\pi} \int_{\mathcal{C}} e^{-t\lambda} (-\mu \Delta  -\lambda)^{-1} \delta (x-\overset{\ast}{x}) \bigg]dx =O(t^{1-n/2}) \quad  \mbox{as}\;\, t\to 0^+\;\; \mbox{for any}\;\, \Omega'' \subset \Omega\setminus U_\epsilon (\partial \Omega),\nonumber\end{eqnarray*}
   and (\ref{2020.10.17-3}) and (\ref{2020.10.17-7}) are proved.

    (b)  \ \  To prove (\ref{2020.10.25-1}), we pick a self-double patch $W$ of $\mathcal{M}$ covering a patch $W\cap \partial \Omega$ endowed (see the diagram below, or p.$\,$53 of \cite{MS}) with local coordinates $x$ such that $\epsilon>x_n>0$ in $W\cap \Omega$, where
  $\epsilon>0$ is some fixed real number; $\,x_n=0$ on $W\cap \partial \Omega$;
  $\; x_n (\overset{*}{x})=-x_n(x)$; and the positive $x_n$-direction is perpendicular to $\partial \Omega$.
    \begin{figure}[h]
\centering
\begin{tikzpicture}[scale=1,line width=0.8]
\clip (-1,1.3) rectangle (6.5,7.4);

\path(70:7) coordinate(A);
\path(30:7) coordinate(B);
\path(70:4.5) coordinate(C);
\path(30:4.5) coordinate(D);

\draw (A) arc (70:30:7);
\draw (C) arc (70:30:4.5);
\draw (A)--(C);
\draw (B)--(D);

\draw[fill=black] (3.177,3.6) circle (1pt);
\draw[fill=black] (4.3,5) circle (1pt);
\draw[fill=black] (3.67,4.32) circle (1pt);

\draw[->,>=stealth] (4.3,5) .. controls (3.5,4.3) and (2.8,3) .. (2.6,2.5);

\node at (68:7.5) {$\Omega^{*}$};
\node at (2.2,2.65) {$\Omega$};
\node at (2.6,2) {$x_n>0$};
\node at (3.35,3.35) {$x$};
\node at (4.5,4.7) {$x^*$};
\node at (-0.5,6.7) {$\partial \Omega$};
\node at (-2,-2){};

\draw (-0.524,6.3709) arc (80.0001:25:8.5);
\end{tikzpicture}
\end{figure}
This products the following effect that
 \begin{eqnarray*} g_{jk} (\overset{*}{x})\!\!&=\!\!&- g_{jk} (x) \quad \, \mbox{for}\;\;
  j<k=n \;\;\mbox{or}\;\; k<j=n,\\  \!\!&=\!\!& g_{jk} (x) \;\;\mbox{for}\;\; j,k<n \;\;\mbox{or}\;\; j=k=n,\\
    g_{jk}(x)\!\!&=\!\!& 0 \;\; \mbox{for}\;\; j<k=n \;\;\mbox{or}\;\; k<j=n \;\;\mbox{on}\;\; \partial \Omega.\end{eqnarray*}
     For any small $n$-dimensional normal coordinate patch $W\subset U_\epsilon (\partial \Omega)$ covering a patch of $W\cap \partial \Omega$,
     without loss of generalization, we let $W\cap \Omega =(W\cap \partial \Omega) \times [0,\epsilon]$.
     Noting that $|x-\overset{\ast} {x}|=x_n-(-x_n)=2x_n$ we find by the method of pseudodifferential operator that
      \begin{eqnarray} \label{aa-3-1-2} &&  \int_{W\cap \Omega} \Big(\frac{i}{2\pi}\!\int_{\mathcal{C}}\! e^{-t\lambda} (-\mu\Delta -\lambda)^{-1} \! \delta(x\!-\!\overset{\ast} {x})d\lambda \Big)dx\\
      &&\quad \;\;=\! \int_0^\epsilon \!dx_n \!\int_{W\cap \partial \Omega}
 \frac{dx'}{(2\pi)^n}\!\int_{{\Bbb R}^{n}} \!e^{i(x-\overset{*}{x})\cdot \xi} \bigg[e^{-t\mu\sum_{j=1}^n \xi_j^2}
 + t\vartheta_{1,1} (x, \xi) e^{-t\mu \sum_{j=1}^n \xi_j^2}\nonumber
 \\
&&\quad \;\; \;\;\;   +\sum_{l=2}^\infty \bigg(\sum_{k=1}^{2l}
    \frac{t^k}{k!} \big( \vartheta_{l,k}(x, \xi)\big) e^{-t\mu\sum_{j=1}^n \xi_j^2}\bigg)\bigg] d\xi\nonumber\\
  &&\quad \;\;  =
 \int_0^\epsilon dx_n \int_{W\cap \partial \Omega}  \frac{dx'}{(2\pi)^n}\!\int_{-\infty}^\infty e^{2ix_n \xi_n}\bigg\{ \int_{{\Bbb R}^{n-1}} \! e^{i 0\cdot  \xi'}  \bigg[e^{-t\mu\sum_{j=1}^n \xi_j^2} \nonumber\\
&&\quad  \;\; \;\; \; +\sum_{l=2}^\infty \bigg(\sum_{k=1}^{2l}
    \frac{t^k}{k!} \big( \vartheta_{l,k}(x, \xi)\big) e^{-t\mu\sum_{j=1}^n \xi_j^2}\bigg)\bigg] d\xi' \bigg\}d\xi_n,\nonumber\end{eqnarray}
  where $\xi=(\xi', \xi_n)\in {\Bbb R}^n$. Here we have used the fact that $\vartheta_{1,1} (x, \xi)=0$
  since we can again take geodesic normal coordinate system at each fixed point $x$.
    It is easy to verify by a straightforward calculation that for any  fixed $\epsilon>0$,
  $$ \int_\epsilon^\infty \frac{1}{(4\pi \mu t)^{\frac{n}{2}}} e^{-\frac{x_n^2}{\mu t}}dx_n = O(t^{1-n/2})\quad \; \, \mbox{as} \;\, t\to 0^+$$
 (In fact, we can prove that for any positive integer $m\ge 1$, $ \;\int_\epsilon^\infty\!\! \frac{1}{(4\pi \mu t)^{\frac{n}{2}}} e^{-\frac{x_n^2}{\mu t}}dx_n \!=\! O(t^{m-n/2}) \; \, \mbox{as} \;\, t\to 0^+$).
  Therefore   \begin{eqnarray} \label{aaa-4-3}&&
   \int_0^\epsilon dx_n \int_{W\cap \partial \Omega}  \frac{dx'}{(2\pi)^n}\int_{-\infty}^\infty e^{2ix_n \xi_n}\bigg[ \int_{{\Bbb R}^{n-1}}  e^{-t\mu\sum_{j=1}^n \xi_j^2} d\xi' \bigg]d\xi_n\\
  && \qquad \quad  =  \int_0^\epsilon dx_n \int_{W\cap \partial \Omega}\frac{1}{(4\pi \mu t)^{\frac{n-1}{2}}} dx'\cdot \frac{1}{(2\pi)}\int_{-\infty}^\infty e^{2ix_n \xi_n}e^{-t\mu \xi_n^2}\,d\xi_n  \nonumber\\
  && \qquad \quad = \int_0^\epsilon \frac{1}{(4\pi \mu t)^{\frac{n}{2}}} e^{-\frac{x_n^2}{\mu t}} dx_n \int_{W\cap \partial \Omega} dx' \nonumber \\
   && \quad \qquad = \int_0^\infty  \frac{1}{(4\pi \mu t)^{\frac{n}{2}}} e^{-\frac{x_n^2}{\mu t}} dx_n \int_{W\cap \partial \Omega} dx'
   -\int_\epsilon^\infty  \frac{1}{(4\pi \mu t)^{\frac{n}{2}}} e^{-\frac{x_n^2}{\mu t}} dx_n \int_{W\cap \partial \Omega} dx'
      \nonumber \\
       && \quad \qquad =\frac{1}{4}\cdot
  \frac{|W\cap \partial \Omega|}{(4\pi \mu t)^{\frac{n-1}{2}}} +O(t^{1-n/2})\, \quad \; \;\mbox{as}\;\; t\to 0^+,\nonumber\end{eqnarray}
It follows from (\ref{2020.11.15-7})  that
    \begin{eqnarray}\label{aaa-4-2,} &&
   \int_0^\epsilon dx_n \int_{W\cap \partial \Omega} \frac{dx'}{(2\pi)^n} \int_{-\infty}^\infty e^{2ix_n \xi_n}\bigg\{ \int_{{\Bbb R}^{n-1}} \bigg[ \sum_{l=2}^\infty \bigg(\sum_{k=1}^{2l}
    \frac{t^k}{k!} \big( \vartheta_{l,k}(x, \xi)\big) e^{-t\mu\sum_{j=1}^n \xi_j^2}\bigg)\bigg] d\xi' \bigg\}d\xi_n\\
    && \qquad \qquad \qquad =O(t^{1-n/2})\quad \; \mbox{as}\;\; t\to 0^+.\nonumber\end{eqnarray}
Thus, by (\ref{aa-3-1-2})--(\ref{aaa-4-2,}) we get
 \begin{eqnarray}\label{aaa-4-4}  && \int_{(W\cap \Omega)\times [0,\epsilon]} \Big(\frac{i}{2\pi}\!\int_{\mathcal{C}}\! e^{-t\lambda} (-\mu\Delta -\lambda)^{\!-\!1} \! \delta(x\!-\!\overset{\ast} {x})d\lambda \Big)dx \nonumber\\
 &&\quad  \quad\, = \! \int_0^\epsilon\! dx_n \!\int_{W\cap \partial \Omega}  \frac{dx'}{(2\pi)^n}\!\int_{-\infty}^\infty \!e^{2ix_n \xi_n}\!\bigg\{ \!\int_{{\Bbb R}^{n-1}} \! \bigg[e^{-t\mu\sum_{j=1}^n \xi_j^2}
    \nonumber\\
 && \quad \quad \quad  +\sum_{l=2}^\infty \bigg(\sum_{k=1}^{2l}
    \frac{t^k}{k!} \big( \vartheta_{l,k}(x, \xi)\big) e^{-t\mu\sum_{j=1}^n \xi_j^2}\bigg)\bigg] d\xi' \bigg\}d\xi_n\nonumber\\
 && \qquad    =\frac{1}{4}\cdot
  \frac{|W \cap \partial \Omega|}{(4\pi \mu t)^{\frac{n-1}{2}}} +O(t^{1-\frac{n}{2}})\, \quad \; \mbox{as}\;\; t\to 0^+,\nonumber \end{eqnarray}
and (\ref{2020.10.25-1}) is proved.

(c) \ \ We now consider $\frac{i}{2\pi}\int_{\mathcal{C}} e^{-t\lambda} (-\mu\Delta -\lambda)^{-1} \sum_{k\ge 1}  \big[ (B_2+C_2) (-\mu \Delta -\lambda)^{-1} \big]^{k} \delta(x-y)\,d\lambda$.
For $k=1$, the Green operator $(-\mu\Delta -\lambda)^{-1}  (B_2+C_2) (-\mu \Delta -\lambda)^{-1}$ is the sum of the pseudodifferential operator $(-\mu\Delta -\lambda)^{-1}  B_2 (-\mu \Delta -\lambda)^{-1}$ and the singular Green operator $(-\mu\Delta -\lambda)^{-1}  C_2 (-\mu \Delta -\lambda)^{-1}$.
According to  symbol formula (\ref{17-9-12.1}), we know (see section 3) that the symbol of  $(-\mu\Delta -\lambda)^{-1}  B_2 (-\mu \Delta -\lambda)^{-1}$ is
\begin{eqnarray}\label{2020.11.17-10}   \sum_{\alpha\ge 0} \frac{(-i)^{|\alpha|}}{\alpha!} \partial^\alpha_{\xi} p(x, \xi,\lambda)\cdot \partial^\alpha_x \Big(\sum_{\beta\ge 0} \frac{(-i)^{|\beta|}}{\beta !}\partial^\beta_\xi b(x,\xi) \cdot \partial^\beta_x p(x,\xi,\lambda)\Big),\end{eqnarray}
where $p(x, \xi,\lambda)=\sum_{l\ge 0} p_{-2-l}  (x, \xi, \lambda)$ (which is given in (a)) and $b(x, \xi)=\sum_{l\le 1} b_l (x, \xi)$  are the symbols of pseudodifferential operators  $(-\mu\Delta -\lambda)^{-1}$ and $B_2$, respectively. From (\ref{2020.11.17-4}), we know that the symbol of the singular Green operator $ (\mu\Delta -\lambda)^{-1} C_2(\mu\Delta -\lambda)^{-1}$  is
 \begin{eqnarray} \label {2020.11.13-1} \quad\quad \;\; - \frac{1}{(2\pi)^2  } \!\int^+\! d\tau_1 \!\int^\oplus \frac{ p(\tau_1) g(\tau_1, \tau_2) p(\tau_2)}{ (\tau_1- \xi_n) (\tau_2- \eta_n) }d\tau_2,\end{eqnarray}
where $g(\xi_n,\eta_n)$ stands for the symbol $g(x', \xi', \xi_n,\eta_n)$ of the singular Green operator $C_2$ which has homogeneous of degree $1$.
For any $x\in \Omega$ and  $y\in \mathcal{M}$, we will consider  \begin{align*} & \frac{1}{(2\pi)^n} \!\int_{\mathbb{R}^n} \!  e^{i(x-y)\cdot \xi}\bigg\{ \! \frac{i}{2\pi}\! \int_{\mathcal{C}}e^{-t\lambda}\bigg[\! \sum_{\alpha\ge 0}\! \frac{\!(\!-i)^{|\alpha|}}{\alpha!} \partial^\alpha_{\xi} p(x, \xi,\lambda)\cdot \partial^\alpha_x \Big(\sum_{\beta\ge 0} \!\frac{\!(\!-i)^{|\beta|}}{\beta !}\partial^\beta_\xi b(x,\xi)\! \cdot \!\partial^\beta_x p(x,\xi,\lambda)\Big)\!\bigg]d\lambda\bigg\}\! d\xi\\
   & + \frac{1}{(2\pi)^{n\!+\!1}}\! \int_{\mathbb{R}^{n-1}}\! d\xi' \!\int^+\! e^{i(x-y)\cdot \xi} d\xi_n \!\int^+\!\bigg\{\!\frac{i}{2\pi} \! \int_{\mathcal{C}}\!e^{-t\lambda} \bigg[\!-\frac{1}{(2\pi)^2  } \!\int^+\! d\tau_1 \!\int^\oplus \frac{ p(\tau_1) g(\tau_1, \tau_2) p(\tau_2)}{ (\tau_1- \xi_n) (\tau_2- \eta_n) }d\tau_2\bigg] \!d\lambda\bigg\}  d\eta_n. \end{align*}
By further calculating (\ref{2020.11.17-10}), we see that  the symbol of homogeneous of degree $j\,$ ($j\le -3$) for the pseudodifferential operator  $(-\mu\Delta -\lambda)^{-1}  B_2 (-\mu \Delta -\lambda)^{-1}$ is a sum of finitely many terms of the form $(p_2(x$, $\xi$,$\lambda))^{m+1} r_{j+2(m+1)}(x, \xi)$, where  $r_{j+2(m+1)}(x$,$ \xi)\in S^{j+2(m+1)}_{1,0}$ and $m\ge 1$.
We also discuss such kind of $g(x', \xi'$, $\xi_n$, $\eta_n)$ such that
 $\frac{ p(\tau_1) g(\tau_1, \tau_2) p(\tau_2)}{(\tau_1- \xi_n) (\tau_2- \eta_n)}$ has the homogeneous order $j$. For each fixed $x\in \Omega$, we use a geodesic normal coordinate frame with center at this $x$. Without loss of generalization, we also assume that $p(x, \xi, \lambda)=p_{-2} (x,\xi, \lambda)$ in a geodesic normal coordinate neighborhood with the center at $x$. Applying residue theorem, we get  \begin{align} \label{2020.11.17-17} & \frac{i}{2\pi}\! \!\int_{\mathcal{C}} \!e^{-t\lambda}\bigg[  (p_{-2}(x, \xi,\lambda))^{m+1} \;r_{j+2(m+1)}(x, \xi)\bigg] d\lambda  = \frac{t^m}{m!} e^{-t\mu\sum_{j=1}^n \xi_j^2}\;r_{j+2(m+1)}(x, \xi)
  \end{align}
and
\begin{align}\label{2020.11.17-18,}& \frac{i}{2\pi}\! \!\int_{\mathcal{C}} \!e^{-t\lambda}\bigg[\!-\frac{1}{(2\pi)^2  } \!\!\int^+\!\! \!d\tau_1\! \!\int^\oplus\! \frac{ p_{-2}(\tau_1) g(\tau_1, \tau_2) p_{-2}(\tau_2)}{ (\tau_1- \xi_n) (\tau_2- \eta_n) }d\tau_2 \bigg] d\lambda\\
 &\qquad \quad =-\frac{1}{(2\pi)^2  }  \!\!\int^+\!\! \!d\tau_1\! \!\int^\oplus\!\bigg[\frac{i}{2\pi}\!  \int_{\mathcal{C}} e^{-t\lambda}\;\frac{ p_{-2}(\tau_1) g(\tau_1, \tau_2) p_{-2}(\tau_2)}{ (\tau_1- \xi_n) (\tau_2- \eta_n) }  d\lambda  \bigg]d\tau_2\nonumber
\\ & \qquad \quad  = -\frac{1}{(2\pi)^2  }  \!\!\int^+\!\! \!d\tau_1\! \!\int^\oplus\!\!\frac{ g(\tau_1, \tau_2) }{ (\tau_1\!-\! \xi_n) (\tau_2\!-\! \eta_n)(\tau_1^2\!-\!\tau_2^2) }\bigg[  e^{-t \mu (\sum_{j=1}^{n-1} \xi_j^2+\tau_2^2)}-  e^{-t \mu (\sum_{j=1}^{n-1}  \xi_j^2+\tau_1^2)} \bigg]d\tau_2,\nonumber
\end{align}
From (\ref{2020.11.15-7}), (\ref{2020.11.17-17}) and variable substitution $\tilde \xi=\sqrt{\mu t}\, \xi$, we have  \begin{align} \label{2020.11.18-2} &
 \bigg| \frac{1}{(2\pi )^n} \int_{\mathbb{R}^n}\! e^{i(x-y)\cdot \xi} \! \bigg[\frac{t^m}{m!} e^{-t\mu\sum_{j=1}^n \xi_j^2}\;r_{j+2(m+1)}(x, \xi)\bigg]d\xi\bigg| \\
&\qquad\quad =\bigg| \frac{t^m (\mu t)^{-\frac{n+j}{2}-m-1} }{(2\pi )^nm!} \,\int_{\mathbb{R}^n} \!e^{i\frac{(x-y)}{\sqrt{t\mu}}\cdot \tilde{\xi}}\bigg[ e^{-\sum_{j=1}^n \tilde \xi_j^2}\;r_{j+2(m+1)}(x, \tilde\xi)\bigg]d\tilde \xi\bigg|\nonumber \\
& \qquad \quad \le \frac{t^m (\mu t)^{-\frac{n+j}{2}-m-1} }{(2\pi )^nm!} \,\int_{\mathbb{R}^n} \!\bigg| e^{-\sum_{j=1}^n \tilde \xi_j^2}\;r_{j+2(m+1)}(x, \tilde\xi)\bigg|d\tilde \xi
=O(t^{-\frac{n+j}{2} -1})\;\;\mbox{as}\,\, t\to 0^+.\nonumber
  \end{align}
  Also, applying variable substitution $\tilde \tau_l=\sqrt{t\mu} \,\tau_l$, $\tilde \xi_n=\sqrt{t\mu}\, \xi_n$ and $\tilde \eta_n=\sqrt{t\mu}\, \eta_n$ we find from   (\ref{2020.11.15-7}),  (\ref{2020.11.17-18,}) and residue theorem that \begin{align}\label{2020.11.18-3} &\frac{1}{(2\pi)^{n+1}}\! \int_{\mathbb{R}^{n-1}}\! \!d\xi' \!\int^+\! \!e^{i(x-y)\cdot \xi} d\xi_n \!\int^+ \!\bigg\{ \!- \frac{1}{(2\pi)^2  }  \!\!\int^+\!\! \!d\tau_1\! \!\int^\oplus\!\frac{ g(\tau_1, \tau_2) }{ (\tau_1- \xi_n) (\tau_2- \eta_n)(\tau_1^2-\tau_2^2) }\\
& \qquad \qquad \qquad \times \bigg[  e^{-t \mu \big(\sum_{j=1}^{n-1} \xi_j^2+\tau_2^2\big)}-  e^{-t \mu \big(\sum_{j=1}^{n-1}  \xi_j^2+\tau_1^2\big)} \bigg]d\tau_2\bigg\} d\eta_n \nonumber\\
&\qquad \;\quad \;\quad=\frac{t^{-\frac{n+j}{2}-1} }{(2\pi)^{n+1}}\! \int_{\mathbb{R}^{n-1}}\! d\tilde\xi' \!\int^+\! e^{i\frac{(x-y)}{\sqrt{\mu t}}\cdot \tilde\xi} d\tilde\xi_n \!\int^+ \!\bigg\{ \! -\! \frac{1}{(2\pi)^2  }  \!\!\int^+\!\! \!d\tilde\tau_1\! \!\int^\oplus\!\frac{ g(\tilde\tau_1, \tilde\tau_2) }{ (\tilde\tau_1- \tilde\xi_n) (\tilde\tau_2- \tilde\eta_n)(\tilde\tau_1^2-\tilde\tau_2^2) }\nonumber\\
& \qquad \qquad \qquad \times \bigg[  e^{- (\sum_{j=1}^{n-1} \tilde\xi_j^2+\tilde\tau_2^2)}-  e^{- (\sum_{j=1}^{n-1} \tilde \xi_j^2+\tilde\tau_1^2)} \bigg]d\tilde\tau_2\bigg\} d\tilde\eta_n=O(t^{-\frac{n+j}{2}-1} )\;\;\mbox{as}\;\, t\to 0^+. \nonumber
\end{align}
Obviously, if $j<-3$, then the asymptotic remainder terms on the right-hand side in (\ref{2020.11.18-2}) and (\ref{2020.11.18-3}) are  $O(t^{1-\frac{n}{2}})$ as $t\to 0^+$. For $j=-3$ (i.e., $m=1$), the principal symbols $r_1(x, \xi)$ and $g_1(x', \xi', \xi_n,\eta_n)$ of  $B_2$ and $C_2$ are just the case as discussed before, each term of which contains a factor of the form $\frac{\partial g_{jk}}{\partial x_l}(x)$. Clearly, $r_1(x, \xi)$ and $g_1(x', \xi', \xi_n, \eta_n)$ vanish because we have taken a geodesic normal coordinate system center at each such $x$.
This implies that the right hand side in (\ref{2020.11.18-2}) and (\ref{2020.11.18-3}) are also  $O(t^{1-\frac{n}{2}})$.
Therefore we obtain that for all $x\in \Omega, y\in \mathcal{M}$,
\begin{align} & \frac{1}{(2\pi )^n} \int_{\mathbb{R}^n} e^{i(x-y)\cdot \xi}\bigg\{ \frac{i}{2\pi}\! \!\int_{\mathcal{C}} \!e^{-t\lambda}\bigg[  (p_{-2}(x, \xi,\lambda))^{m+1} \;r_{j+2(m+1)}(x, \xi)\bigg] d\lambda\bigg\} d\xi\nonumber \\
& \quad \quad \quad+ \frac{1}{(2\pi)^{n+1}}\! \int_{\mathbb{R}^{n-1}}\! d\xi' \!\int^+\! e^{i(x-y)\cdot \xi} d\xi_n \!\int^+\!\bigg\{\!\frac{i}{2\pi} \! \int_{\mathcal{C}}\!e^{-t\lambda} \bigg[\!-\frac{1}{(2\pi)^2  } \!\int^+\! d\tau_1 \!\int^\oplus \frac{ p(\tau_1) g(\tau_1, \tau_2) p(\tau_2)}{ (\tau_1- \xi_n) (\tau_2- \eta_n) }d\tau_2\bigg] \!d\lambda\bigg\}  d\eta_n\nonumber\\
&\qquad \quad=O(t^{1-\frac{n}{2}}) \quad \mbox{as}\;\, t\to 0^+,  \nonumber\end{align}
i.e.,  \begin{eqnarray} \label{2020.11.17-20} &\frac{i}{2\pi}\int_{\mathcal{C}} e^{-t\lambda} (-\mu\Delta -\lambda)^{-1}  \big[ (B_2+C_2) (-\mu \Delta \!-\!\lambda)^{-1} \big] \delta(x-y)\,d\lambda \! =\!O(t^{1-\frac{n}{2}}) \;\; \mbox{as}\;\, t\to 0^+.\nonumber\end{eqnarray}
 In particular, \begin{eqnarray} \label{2020.11.18-2.0}\quad\; \;\;\;\int_{\Omega'}\!\bigg\{\!\frac{i}{2\pi}\!\int_{\mathcal{C}} e^{-t\lambda} (-\mu\Delta -\lambda)^{-1}  \big[ (B_2+C_2) (-\mu \Delta -\lambda)^{-1} \big] \delta(x-x)\,d\lambda\!\bigg\}dx \! =\! O(t^{1-\frac{n}{2}}) \;\; \mbox{as}\;\, t\to 0^+\end{eqnarray}
and \begin{eqnarray}\label{2020.11.18-2.1}\;\;\;\;\quad\int_{\Omega'}\bigg\{\frac{i}{2\pi}\int_{\mathcal{C}} e^{-t\lambda} (-\mu\Delta -\!\lambda)^{-1}  \big[ (B_2+C_2) (-\mu \Delta -\!\lambda)^{-1} \big] \delta(x-\!\overset{\ast}{x})\,d\lambda\bigg\}dx  \!=\!O(t^{1-\!\frac{n}{2}}) \;\; \mbox{as}\;\, t\to 0^+
.\end{eqnarray}

(d) \  For each integer $k\ge 2$, let the full symbol of the Green operator $(-\mu\Delta- \lambda)^{-1} \big[ (B_2+C_2)(-\mu\Delta-\lambda)^{-1}\big]^k$ be $\big(\sum_{l\ge 0} \beta_{-2-k-l} (x, \xi, \lambda);  c_{-2-k-l} (x', \xi', \xi_n,\eta_n,\lambda)\big).$
Applying the symbol composition formula (\ref{17-9-12.1}), we find that the symbol with homogeneous of degree $j\,$ ($j\le -k - 2$) for the pseudodifferential  operator $(-\mu\Delta -\lambda)^{-1} (B_2)^k (-\mu \Delta -\lambda)^{-1} $ is a sum of finitely many terms of the form $ \big((p_{-2} (x, \xi, \lambda))^{m+1} r_{j+2(m+1)}(x, \xi)$, where  $r_{j+2(m+1)}(x, \xi)\in S^{j+2(m+1)}_{1,0}$ and $m\ge 1$.
We also discuss such kind of $c_j(x', \xi',\xi_n,\eta_n, \lambda)$ which has the homogeneous order $j$ for the singular Green operator
$ (-\mu \Delta-\lambda)^{-1} \big[(B_2+C_2) (-\mu \Delta-\lambda)^{-1}\big]^k - (-\mu \Delta-\lambda)^{-1} [B_2
 (-\mu \Delta-\lambda)^{-1}]^k$.
 Note that the symbol $c_{j}(x', \xi', \xi_n, \eta_n,\lambda)$ of the  singular Green operator $ (-\mu \Delta-\lambda)^{-1} \big[(B_2+C_2) (-\mu \Delta-\lambda)^{-1}\big]^k - (-\mu \Delta-\lambda)^{-1} [B_2
 (-\mu \Delta-\lambda)^{-1}]^k$ is given by (\ref{2020.12.4-4,}).
Completely similar to the argument of (c) (i.e., by applying residue theorem and variable substitution) we get
 that for all $x\in \Omega, y\in \mathcal{M}$,
\begin{align} \label{2020.11.24-30} & \frac{1}{(2\pi )^n} \int_{\mathbb{R}^n} e^{i(x-y)\cdot \xi}\bigg\{ \frac{i}{2\pi}\! \!\int_{\mathcal{C}} \!e^{-t\lambda}\bigg[  (p_{-2}(x, \xi,\lambda))^{m+1} \;r_{j+2(m+1)}(x, \xi)\bigg] d\lambda\bigg\} d\xi \\
& \qquad \quad \quad+ \frac{1}{(2\pi)^{n+1}}\! \int_{\mathbb{R}^{n-1}}\! d\xi' \!\int^+\! e^{i(x-y)\cdot \xi} d\xi_n \!\int^+\!\bigg\{ \frac{i}{2\pi} \int_{\mathcal{C}} e^{-t\lambda} c_j(x',\xi',\xi_n,\eta_n, \lambda) d\lambda \bigg\} d\eta_n\nonumber\\
&\qquad \quad=O(t^{1-\frac{n}{2}}) \quad \mbox{as}\;\, t\to 0^+.  \nonumber\end{align}
Consequently, for any $k\ge 1$, by (\ref{2020.11.18-2.0})---(\ref{2020.11.18-2.1}) and (\ref{2020.11.24-30}) we obtain  \begin{eqnarray} \label{2020.11.17-20} &\frac{i}{2\pi}\int_{\mathcal{C}} e^{-t\lambda} (-\mu\Delta -\lambda)^{-1}  \big[ (B_2+C_2) (-\mu \Delta -\lambda)^{-1} \big]^{k} \delta(x-y)\,d\lambda  =O(t^{1-\frac{n}{2}}) \quad \mbox{as}\;\, t\to 0^+,\nonumber\end{eqnarray}
which implies (\ref{2020.10.17-4}) and (\ref{2020.10.17-40}).

\vskip 0.35 true cm

 Finally, combining (\ref{c4-23}), (\ref{2020.10.17-1}) and (\ref{2020.10.17-3})---(\ref{2020.10.17-40}), we get
 \begin{eqnarray} \label{a-4-1-3} \int_{\Omega}e^{-t (-\tilde{A}_{FF} - \tilde{R}_F ) }\delta (x-x)\, dx \!\!\! &=\!\!\!&
 \int_{\Omega} K_2^{-} (t, x, x) \,dx =
 \int_{\Omega} \big[K_2 (t, x, x)- K_2 (t, x, \overset{\ast} {x})\big]dx \\
\!\!\! &=\!\!\!&(4\pi\mu  t)^{-n/2} \bigg[|\Omega|-\frac{1}{4} \sqrt{4\pi \mu t}  \,|\partial \Omega|+O(t)\bigg] \quad \, \mbox{as}\;\; t\to 0^+.\nonumber\end{eqnarray}

\vskip 0.16 true cm

Step 2. \ We will calculate the asymptotic expansion for the integral of the trace of the strongly continuous semigroup $e^{-t ({\mathbf{A}}^2+{\mathbf{R}})^{-1/2}}$ as $t\to 0^{+}$. Let us first consider the semigroup $ e^{-t ({\mathbf{A}}^2+{\mathbf{R}})^{-1}}$. From (\ref{2020.9.16-4}) we have
\begin{eqnarray*}
 e^{-t({\mathbf{A}}^2+\mathbf{R})^{-1}}\mathbf{f}(x)\!\!\!\!\!\!\!\!&& = e^{-t(\mu^2 \Delta^2 \, \mathbf{I} -\mathbf{B}_1-\mathbf{C}_1)}\mathbf{f}(x)\\
 \!\!\! \!\!\!\!\!\!&\!&\!\!  = \frac{i}{2\pi} \int_{\mathcal{C}} e^{-t\lambda} (\mu^2 \Delta^2\,\mathbf{I} -\mathbf{B}_1-\mathbf{C}_1-\lambda \mathbf{I})^{-1}\mathbf{f}(x) d\lambda\\
  \!\!\!\!\!\!\!\!&& = \frac{i}{2\pi} \int_{\mathcal{C}} e^{-t\lambda} (\mu^2\Delta^2 \,\mathbf{I}-\lambda\, \mathbf{I})^{-1} \sum_{k=0}^\infty
  \big[(\mathbf{B}_1 +\mathbf{C}_1)(\mu^2\Delta^2 \,\mathbf{I}-\lambda\, \mathbf{I})^{-1} \big]^k  \mathbf{f}(x)\, d\lambda,\nonumber\end{eqnarray*}
where $\mathbf{B}_1$ (respectively, $\mathbf{C}_1$) is a pseudodifferential (respectively, singular Green) matrix operator of order $3$.
Denote by $\mathbf{p} (x, \xi, \lambda) = \sum_{l\ge 0} \mathbf{p}_{-4-l} (x, \xi, \lambda)\in S^{-4}_{1,0}$ and $\big(\sum_{l\ge 0} {\boldsymbol{\beta}}_{-4-k-l} (x, \xi, \lambda); {\boldsymbol{\omega}}_{-4-k} (x', \xi', \xi_n,\eta_n,\lambda)\big)$ the full symbols of the pseudodifferential operators $(\mu^2 \Delta^2 \, \mathbf{I} -\lambda)^{-1}$ and Green operator $(\mu^2\Delta^2 \,\mathbf{I}-\lambda\, \mathbf{I})^{-1}  \big[(\mathbf{B}_1 +\mathbf{C}_1)(\mu^2\Delta^2 \,\mathbf{I}-\lambda\, \mathbf{I})^{-1} \big]^k $, respectively.
 In fact, \begin{eqnarray*} &&\mathbf{p}_{-4} =\Big(\mu^2\big( \sum_{j=1}^n \xi_j^2\big)^2 -\lambda\Big)^{-1}\mathbf{I},  \quad\;  \mathbf{p}_{-4-1} = \mathbf{h}_{1,1} (x, \xi) \mathbf{p}_{-4}^2,\\
 &&    \mathbf{p}_{-4 -l} = \sum\limits_{k=1}^{2l} \mathbf{h}_{l,k} (x, \xi)  \mathbf{r}_{-4}^{k+1},\, \cdots, \;\, l\ge 2,\end{eqnarray*}
where $\mathbf{p}_{-4-l}(x, \xi)$ are homogeneous in $\xi$ of degree
$-4-l$ for $|\xi|>1$, and $\mathbf{h}_{l,k}$ independent of $\lambda$ and homogeneous of degree $4k-l$ in $\xi$ for $|\xi|\ge 1$.
Clearly, for any $x, y\in \mathcal{M}$,
\begin{eqnarray}\label{2020.9.19-9}\boldsymbol{\Xi}(t,x,y)\!\!\!\!\!\!\!\!\! &\!&= e^{-t(\mu\Delta^2 \, \mathbf{I}-\mathbf{B}''_2 -\mathbf{C}''_2)}{\boldsymbol{\delta}}(x-y) \\
 \!\!\!\!\!\!\!\!&& = \frac{i}{2\pi} \int_{\mathcal{C}} e^{-t\lambda} (\mu^2\Delta^2 \,\mathbf{I}-\lambda\, \mathbf{I})^{-1} \sum_{k=0}^\infty
  \big[(\mathbf{B}_1 +\mathbf{C}_1)(\mu^2\Delta^2 \,\mathbf{I}-\lambda\, \mathbf{I})^{-1} \big]^k  \boldsymbol{\delta}(x-y)\, d\lambda.\nonumber\end{eqnarray}
  For any $k\ge 1$, according to the argument in section 3, the symbol of the Green operator  $(\mu^2\Delta^2 \,\mathbf{I}-\lambda\, \mathbf{I})^{-1}  \big[(\mathbf{B}_1 +\mathbf{C}_1)(\mu^2\Delta^2 \,\mathbf{I}-\lambda\, \mathbf{I})^{-1} \big]^k$ with homogeneous of degree $j\;$  ($j\le -4-k$) is
  $((\mathbf{p}_{-4} (x, \xi,\lambda))^{m+1}  {\boldsymbol{\iota}}_r (x, \xi)$;
${\boldsymbol{\omega}}_{j} (x', \xi',\xi_n,\eta_n$, $\lambda))$, where $\boldsymbol{\iota}_r(x,\xi)$
  is pseudodifferential operator of order $r\;$ ($r=4(m+1)+j$), and ${\boldsymbol{\omega}}_j(x', \xi',\xi_n, \eta_n,\lambda)$ is the singular Green symbol of order $j$.
   Completely similar to the discussion as in step 1, by taking geodesic normal coordinate with center at each pint $x\in \Omega$ we have that for any $\Omega'\subset \Omega$, $\,\Omega''\subset \Omega\setminus U_{\epsilon} (\partial \Omega)$ and $W= (W\cap \partial \Omega) \times [0, \epsilon]\subset U_{\epsilon} (\partial \Omega)$,
  \begin{eqnarray}
&& \int_{\Omega'}\!\bigg( \!\frac{i}{2\pi}\! \int_{\mathcal{C}}\! e^{-t\lambda} (\mu^2\Delta^2 \,\mathbf{I}-\lambda\, \mathbf{I})^{-1} \boldsymbol{\delta}(x\!-\!x)\, d\lambda\bigg) dx \!= \!\int_{\Omega'}\!\bigg[\!\frac{i}{2\pi}\! \int_{\mathcal{C}} \!e^{-t\lambda}  \bigg(\! \frac{1}{(2\pi)^n} \!\int_{\mathbb{R}^n} \!\mathbf{p}_{-4} (x, \xi, \lambda)d\xi \bigg)d\lambda\bigg]dx\\
&& \quad \;\quad  =n\int_{\Omega'}\bigg[ e^{-t\mu^2(\sum_{j=1}^n \xi_j^2)^2 } d\xi\bigg]dx=\frac{n}{(4\pi \mu t)^{n/4}}\; \frac{\Gamma (\frac{n}{4})\,|\Omega|}{2(4\pi \mu)^{\frac{n}{4}}\Gamma(\frac{n}{2}) } +O(t^{\frac{1}{2}-\frac{n}{4}}) \quad \;\, \mbox{as}\,\; t\to 0^+, \nonumber\end{eqnarray}
 \begin{eqnarray}
&&\int_{\Omega''}\bigg( \frac{i}{2\pi} \int_{\mathcal{C}} e^{-t\lambda} (\mu^2\Delta^2 \,\mathbf{I}-\lambda\, \mathbf{I})^{-1} \boldsymbol{\delta}(x-\overset{\ast}{x})\, d\lambda\bigg)dx=O(t^{\frac{1}{2}- \frac{n}{4}}) \quad \, \mbox{as}\;\, t \to 0^+, \end{eqnarray}
  \begin{eqnarray}
&&\int_{(W\cap \partial \Omega)\times [0,\epsilon]}\bigg( \frac{i}{2\pi} \int_{\mathcal{C}} e^{-t\lambda} (\mu^2\Delta^2 \,\mathbf{I}-\lambda\, \mathbf{I})^{-1} \boldsymbol{\delta}(x-\overset{\ast}{x})\, d\lambda\bigg)dx \\
&& \;\;=\int_{(W\cap \partial \Omega)\times [0,\epsilon]}\!\bigg[\!\frac{i}{2\pi} \int_{\mathcal{C}} e^{-t\lambda}  e^{i (x-\overset{\ast}{x})\cdot \xi} \bigg(\frac{1}{(2\pi)^n} \int_{\mathcal{C}} \mathbf{p}_{-4} (x, \xi, \lambda) d\xi \bigg) d\lambda\!\bigg)dx \nonumber\\
&& \;\;=
n \int_{W}\!\bigg[ e^{i (x-\overset{\ast}{x})\cdot \xi}  e^{-t\mu^2(\sum_{j=1}^n \xi_j^2)^2 } d\xi\bigg]dx\nonumber \\
&&aa\;\;=
-\frac{n}{(4\pi \mu t)^{n/4}}\bigg(\!  \frac{ t^{\frac{1}{4}}|(\partial \Omega)\cap W|}{4^{\frac{n}{2}+ \frac{1}{4}} \pi^{\frac{n}{4} -1} \mu^{\frac{n}{4}-\frac{1}{2}}\Gamma(\frac{n+1}{4})}
 + O(t^{\frac{1}{2}})\!\bigg) \; \, \mbox{as}\,\, t\to 0^+,\nonumber\end{eqnarray}
\begin{eqnarray}
&&\int_{\Omega'} \bigg(\frac{i}{2\pi} \int_{\mathcal{C}} e^{-t\lambda} (\mu^2\Delta^2 \,\mathbf{I}-\lambda\, \mathbf{I})^{-1} \sum_{k=1}^\infty
  \big[(\mathbf{B}_1 +\mathbf{C}_1)(\mu^2\Delta^2 \,\mathbf{I}-\lambda\, \mathbf{I})^{-1} \big]^k  \boldsymbol{\delta}(x-x)\, d\lambda\bigg)dx\\
 && \qquad \qquad  \qquad = O(t^{\frac{2-n}{4}}) \quad \mbox{as}\;\, t \to 0^+,\nonumber
\end{eqnarray}
\begin{eqnarray}
&& \int_{\Omega'} \bigg(\frac{i}{2\pi} \int_{\mathcal{C}} e^{-t\lambda} (\mu^2\Delta^2 \,\mathbf{I}-\lambda\, \mathbf{I})^{-1} \sum_{k=1}^\infty
  \big[(\mathbf{B}_1 +\mathbf{C}_1)(\mu^2\Delta^2 \,\mathbf{I}-\lambda\, \mathbf{I})^{-1} \big]^k  \boldsymbol{\delta}(x-\overset{\ast}{x})\, d\lambda\bigg)dx\\
  && \qquad \qquad = O(t^{\frac{2-n}{4}}) \quad \mbox{as}\;\, t \to 0^+.\nonumber
\end{eqnarray}
 Here we have used the trace of the matrix, by which the coefficients of asymptotic expansion are multiplied by $n$, and the fact that for any ${\mathbf{h}}(x, \xi)\in S^{m}_{1,0}$, \begin{eqnarray*} \int_{{\Bbb R}^n}
   \big({\mathbf{h}}(x, \xi)\big)e^{-t\mu^2(\sum_{k=1}^n \xi_k^2)^2} d\xi=O(t^{- (m+n)/4}).\end{eqnarray*}
 Hence   \begin{eqnarray} &&\int_{\Omega} {\boldsymbol{\Xi}}^{-} (t, x,x) dx =\int_{\Omega}  \Big({\boldsymbol{\Xi}} (t, x,x) -{\boldsymbol{\Xi}}(t, x,\overset{\ast}{x}) \Big)dx = \int_{\Omega}  \Big[\mbox{Tr}\big(e^{-t(\mu^2 \Delta^2 \mathbf{I} -\mathbf{B}_1-\mathbf{C}_1)}\big)\Big] dx
   \\
 &&\quad \;\;=\frac{n}{(4\pi \mu t)^{n/4}}\bigg[ \frac{\Gamma (\frac{n}{4})\,|\Omega|}{2(4\pi \mu)^{\frac{n}{4}}\Gamma(\frac{n}{2}) } - \frac{ t^{\frac{1}{4}}|\partial \Omega|}{4^{\frac{n}{2}+ \frac{1}{4}} \pi^{\frac{n}{4} -1} \mu^{\frac{n}{4}-\frac{1}{2}}\Gamma(\frac{n+1}{4})}
 + O(t^{\frac{1}{2}})\bigg]\quad \; \, \mbox{as}\,\, t\to 0^+.\nonumber\end{eqnarray}

 It is well-known (see, for example, (5.22) of p$\,$247 of \cite{Ta1}, \cite{GG}) that for $\lambda\ge 0$,
\begin{eqnarray*} e^{-t\sqrt{\lambda}} =\int_0^\infty \frac{t}{\sqrt{4\pi s^3}} e^{-t^2/(4s)} e^{-s\lambda} ds,\end{eqnarray*}
i.e., the Laplace transform of $\frac{t}{\sqrt{4\pi s^3}} e^{-t^2/(4s)}$ is $e^{-t\sqrt{\lambda}}$.
By applying the spectral theorem, we get that for all $t>0$,
\begin{eqnarray} \label{b.12}  e^{-t(\mu^2 \Delta^2\,\mathbf{I} -\mathbf{B}_1-\mathbf{C}_1 )^{1/2}}   = \int_0^\infty  \frac{t}{\sqrt{4\pi s^3}} e^{-t^2/(4s)} e^{-s (\mu^2 \Delta^2\,\mathbf{I} -\mathbf{B}_1-\mathbf{C}_1)}ds. \end{eqnarray}
 Therefore, we have  \begin{eqnarray} \label{4-1-2}&& \quad\; \int_{\Omega} \bigg\{  \mbox{Tr}\big( e^{-t({\mathbf{A}}^2+ \mathbf{R})^{-1/2}}\big) \bigg\}dx  = \int_{\Omega} \bigg\{  \mbox{Tr} \bigg(e^{-t(\mu^2 \Delta^2 \mathbf{I} -\mathbf{B}_1-\mathbf{C}_1)^{1/2}}\bigg)\bigg\}dx \\
  &&  \qquad \qquad =   \int_0^\infty \frac{t}{\sqrt{4\pi s^3}} e^{-t^2/(4s)}\bigg\{\!\int_\Omega \bigg[
\mbox{Tr}\bigg( e^{-s (\mu^2 \Delta^2\,\mathbf{I} -\mathbf{B}_1-\mathbf{C}_1)}\bigg)\bigg] dx \bigg\}ds \nonumber\\
&& \qquad \qquad  =  \int_0^\infty \frac{t}{\sqrt{4\pi s^3}} e^{-t^2/(4s)}\frac{n}{(4\pi\mu  s)^{n/4}}\bigg[ \frac{\Gamma (\frac{n}{4})\,|\Omega|}{2(4\pi \mu)^{\frac{n}{4}}\Gamma(\frac{n}{2}) } - \frac{ s^{\frac{1}{4}}\,|\partial \Omega|}{4^{\frac{n}{2}+ \frac{1}{4}} \pi^{\frac{n}{4} -1} \mu^{\frac{n}{4}-\frac{1}{2}}\Gamma(\frac{n+1}{4})}
 + O(s^{\frac{1}{2}})\bigg]
 ds\nonumber\\
&& \qquad \qquad  = \frac{n}{(4\pi\mu  t)^{n/2}} \bigg[|\Omega|-\frac{1}{4} \sqrt{4\pi\mu t} \; |\partial \Omega|+O(t)\bigg] \quad \, \mbox{as}\;\; t\to 0^+.\nonumber\end{eqnarray}

\vskip 0.16 true cm

 Step 3.  \  Recall that the kernel spaces of ${\mathbf{A}}^2$ and ${\tilde{A}}_{FF}$ have the same dimensional number
  $m_0$.  As in proof of Lemma 4.4, let $\{\boldsymbol{\phi_j}\}_{j=1}^\infty$ (respectively, $\{{\tilde{p}}_j\}_{j=1}^\infty$) are the orthonormal eigenvectors corresponding to all non-zero eigenvalues $\{\alpha_j\}_{j=1}^\infty$ (respectively, $\{{\tilde{\beta}}_j\}_{j=1}^\infty$) of $\mathbf{A}$ (respectively,
  $-{\tilde{A}}_{{}_{FF}}$), and $\{{\boldsymbol{\psi}}_j\}_{j=1}^{m_0}$ (respectively, $\{{\tilde{q}}_j\}_{j=1}^{m_0}$) is an orthonormal basis of $\mbox{ker}\,\mathbf{A}$ (respectively, of $\mbox{ker}\,{\tilde{A}}_{{}_{FF}}$). Then the operators ${\mathbf{A}}^2+\mathbf{R}$ and $-{\tilde{A}}_{{}_{FF}}-{\tilde{R}_F}$ have the following representations:
 \begin{eqnarray*}   {\mathbf{A}}^2 +\mathbf{R}= \sum_{j=1}^\infty \alpha_j^2 {\boldsymbol{\phi}}_j +  \sum_{j=1}^{m_0} \varrho^2  {{\boldsymbol{\psi}}}_j, \\
  -{\tilde{A}}_{{}_{FF}}-{\tilde{R}_F}=  \sum_{j=1}^\infty {\tilde{\beta}}_j {\tilde{p}}_j +  \sum_{j=1}^{m_0} \varrho  {\tilde{q}}_j, \end{eqnarray*}
  where $\varrho>0$ is a constant as pointed out before.
   Furthermore, the semigroups can be represented as  \begin{eqnarray}  e^{-t({\mathbf{A}}^2 +\mathbf{R})^{-\frac{1}{2}}}= \sum_{j=1}^\infty e^{-t \frac{1}{|\alpha_j|}} {\boldsymbol{\phi}}_j +  \sum_{j=1}^{m_0} e^{-t \frac{1}{\varrho}} {\boldsymbol{\psi}}_j, \\
 e^{-t(- {\tilde{A}}_{{}_{FF}}-{\tilde{R}_F})^{-1}}=  \sum_{j=1}^\infty e^{-t \frac{1}{{\tilde{\beta}}_j}} {\tilde{p}}_j +  \sum_{j=1}^{m_0} e^{-t\frac{1}{\varrho}}  {\tilde{q}}_j; \end{eqnarray}
  in addition, the  fundamental solutions of $e^{-t({\mathbf{A}}^2 +\mathbf{R})^{-\frac{1}{2}}}$  and $e^{-t( -{\tilde{A}}_{{}_{FF}}-{\tilde{R}_F})^{-1}}$, respectively, are
 \begin{eqnarray*}  {\mathbf{K}}_1(t,x,y) =  \sum_{j=1}^\infty e^{-t \frac{1}{|\alpha_j|}} {\boldsymbol{\phi}}_j(x) \otimes  {\boldsymbol{\phi}}_j(y) +  \sum_{j=1}^{m_0}   e^{-t \frac{1}{\varrho}}{\boldsymbol{\psi}}_j(x)\otimes {\boldsymbol{\psi}}_j(y)\;\qquad \quad \qquad  \\
 \qquad \qquad\mbox{and} \quad  K_2(t,x,y)=  \sum_{j=1}^\infty e^{-t \frac{1}{{\tilde\beta}_j}} {\tilde{p}}_j (x)\cdot {\tilde{p}}_j(y)+  \sum_{j=1}^{m_0} e^{-t\frac{1}{\varrho}} {\tilde{q}}_j(x)\cdot {\tilde{q}}_j (y)\qquad \qquad\qquad  \end{eqnarray*}
 with uniform convergence on compact figures of $(0,\infty) \times \Omega \times \Omega$, and the  integrals of the traces of the semigroups  $e^{-t({\mathbf{A}}^2+\mathbf{R})^{-\frac{1}{2}}}$ and
  $e^{-t(-{\tilde{A}}_{{}_{FF}}-\tilde{R}_F)^{-1}}$ are easily evaluated as (see, for example, \cite{Min}, \cite{CLN})
\begin{eqnarray*} \int_\Omega \Big[\mbox{Tr}\big(e^{-t ({\mathbf{A}}^2+\mathbf{R})^{-\frac{1}{2}}}\big)\Big] dx \!\!\!&\!=\!\!\!&\!\!\sum_{j=1}^\infty e^{-t \frac{1}{|\alpha_j|}} \int_\Omega |{\boldsymbol{\phi}}_j(x)|^2 dx +  \sum_{j=1}^{m_0} e^{-t \frac{1}{\varrho}} \int_\Omega |{{\boldsymbol{\psi}}}_j(x)|^2 dx \qquad \\
 \! \!\!&\! =\!\!\!\! &\!\!\sum_{j=1}^\infty e^{-t \frac{1}{|\alpha_j|}} +  \sum_{j=1}^{m_0} e^{-t \frac{1}{\varrho}} ,\\
  \int_\Omega \Big[ \mbox{Tr}\big(e^{-t (-{\tilde{A}}_{{}_{FF}} -{\tilde{R}_F})^{-1}}\big)\Big]dx  &=& \sum_{j=1}^\infty e^{-t \frac{1}{{\tilde\beta}_j}} \int_\Omega |{\tilde{p}}_j(x)|^2 dx +  \sum_{l=1}^{m_0} e^{-t \frac{1}{\varrho}} \int_\Omega |{\tilde{q}}_j(x)|^2 dx \quad \\
\!\!\!&\! =\!\!\!&\!\! \sum_{j=1}^\infty e^{-t \frac{1}{{\tilde\beta}_j}} +  \sum_{j=1}^{m_0} e^{-t \frac{1}{\varrho}}.\end{eqnarray*}
  Recall that  \begin{gather*}   \begin{pmatrix} S^{-1}  & 0 \\
    0 &  \mathbf{A}_{{}_{FF}}  \end{pmatrix} \end{gather*}
  is defined on $J\oplus  F$, and its all non-zero eigenvalues are the same as that of $\mathbf{A}$. From Lemma 4.3 we know that the space of all eigenvectors of the operator $\begin{pmatrix} S^{-1}  & 0 \\
    0 &  \mathbf{A}_{{}_{FF}}  \end{pmatrix}$ corresponding to non-zero eigenvalue $\tau$ is just the orthogonal sum of the space of all eigenvectors of $S^{-1}$ and the space of all eigenvectors of $\mathbf{A}_{{}_{FF}}$ corresponding to the same eigenvalue $\tau$.
       Also, $\mathbf{A}_{{}_{FF}}$ (defined on $F$) can equivalently reduces to $\tilde{A}_{{}_{FF}}$ (defined on $H^1_0(\Omega)$). In view of \begin{eqnarray*} {\mathbf{K}}(t, x,y)=
  \sum_{j=1}^\infty e^{-t\lambda_j}{\mathbf{u}}_j(x)\otimes {\mathbf{u}}_j(y),\end{eqnarray*}
 we have \begin{eqnarray*}\int_\Omega \Big[\mbox{Tr}\big(e^{-tS}\big)\Big]dx = \int_{\Omega}\Big[ \mbox{Tr}\big( {\mathbf{K}}(t, x,x)\big)\Big]dx= \sum_{j=1}^\infty e^{-t\lambda_j} \int_\Omega |{\mathbf{u}}_j(x)|^2 dx = \sum_{j=1}^\infty e^{-t\lambda_j},\end{eqnarray*}
 where $\{{\mathbf{u}}_j\}_{j=1}^\infty$ are the orthonormal eigenvectors corresponding to the Stokes eigenvalues $\{\lambda_j\}_{j=1}^\infty$.
    This implies
   \begin{eqnarray} \label{4-1-3} \int_\Omega \Big[  \mbox{Tr}\big( e^{-tS}\big)
\Big] dx&=& \sum_{j=1}^\infty e^{-t \lambda_j}\\
&=& \sum_{j=1}^\infty e^{-t \frac{1}{|\alpha_j|}}- \sum_{j=1}^\infty e^{-t \frac{1}{\tilde\beta_j}}\nonumber\\
&=&\sum_{j=1}^\infty e^{-t \frac{1}{|\alpha_j|}} +  \sum_{j=1}^{m_0} e^{-t \frac{1}{\varrho}}- \sum_{j=1}^\infty e^{-t \frac{1}{{\tilde\beta}_j}}
- \sum_{j=1}^{m_0} e^{-t \frac{1}{\varrho}}\nonumber\\
&=& \int_\Omega\Big[\mbox{Tr}\big( e^{-t({\mathbf{A}}^2+\mathbf{R})^{-1/2}}\big)\Big]dx  - \int_\Omega \Big[\mbox{Tr}\big(
e^{-t(-{\tilde{A}}_{{}_{FF}}-{\tilde{R}})^{-1}}\big)\Big] dx.\nonumber\end{eqnarray}
By (\ref{4-1-3}), (\ref{2020.9.16-4}), (\ref{2020.10.3-1}), (\ref{a-4-1-3}) and (\ref{4-1-2}) we finally obtain  \begin{eqnarray} \label{4-1-4} \sum_{j=1}^\infty e^{-t \lambda_j} \!\!&=\!\!& \int_{\Omega} \bigg\{ \mbox{Tr}\big( e^{-t({\mathbf{A}}^2+\mathbf{R})^{-1/2}}\big)\bigg\}dx - \int_\Omega \big(
e^{-t(-{\tilde{A}}_{{}_{FF}}-{\tilde{R}})^{-1}}\big)dx \\
\!\!&=\!\!& \int_\Omega \Big[ \mbox{Tr}\big( e^{-t(\mu^2 \Delta^2\,\mathbf{I}-\mathbf{B}_1-\mathbf{C}_1)^{1/2}}\big)\Big]dx  - \int_{\Omega} \big(
e^{-t(-\mu \Delta -B_2-C_2)}\big)dx \nonumber\\
\!\!&=\!\!& \frac{n}{(4\pi \mu t)^{n/2}} \bigg[ |\Omega|- \frac{1}{4} \sqrt{4\pi\mu  t} \,|\partial \Omega| +O(t)\bigg]\nonumber\\
\!&& \!- \frac{1}{(4\pi \mu t)^{n/2}} \bigg[ |\Omega|- \frac{1}{4} \sqrt{4\pi\mu  t} \,|\partial \Omega| +O(t)\bigg] \nonumber\\
\!\!&=\!\!& \frac{n-1}{(4\pi \mu t)^{n/2}} \bigg[ |\Omega|- \frac{1}{4} \sqrt{4\pi\mu  t} \,|\partial \Omega| +O(t)\bigg] \quad \;\, \mbox{as}\;\; t\to 0^+. \nonumber \end{eqnarray}\qed

\vskip 0.32 true cm

\noindent{\bf Remark 5.1.} \   {\it The result in step 2 of above proof can also be obtained by an alternative method: since $({\mathbf{A}}^2+\mathbf{R})^{-1/2}=(\mu^2 \Delta^2\,\mathbf{I}-\mathbf{B}_1-\mathbf{C}_1)^{1/2}$, from the symbol composition formulas (\ref{17-9-12.1}) (for the pseudodifferential operators) and (\ref{2020.11.12-4})---(\ref{2020.11.12-6}) (for the Green operators) we get the principal symbol of $(\mu^2 \Delta^2 \,\mathbf{I}-\mathbf{B}_1-\mathbf{C}_1)^{1/2}$ being $\mu \sum_{j=1}^n \xi_j^2 \, \mathbf{I}+{\boldsymbol{\Phi}}(x, \xi) +{\boldsymbol{\Psi}}(x, \xi)$, where ${\boldsymbol{\Phi}}(x,\xi)\in S^{1}_{1,0}$ and ${\boldsymbol{\Psi}}(x,\xi)$ is the symbol of a singular Green operator of degree $1$ and class $r$.
Similar to the discussion in step 1 we can also get (\ref{4-1-2}).}

\vskip 0.35 true cm

 Now, we  use the Stokes spectral invariants which have been obtained from Theorem 1.1 to finish the proof of Corollary  1.2.

\vskip 0.35 true cm

 \noindent  {\it Proof of Corollary 1.2.} \  Since the Stokes spectrum for the domain $\Omega$ coincides with that for the ball $B_r$, by Theorem 1.1 we know that the first two coefficients  $\frac{(n-1)|\Omega|}{(4\pi \mu t)^{n/2}}$ and $\frac{-(n-1)|\partial \Omega|}{4(4\pi \mu t)^{(n-1)/2}}$
 of the  asymptotic expansion in (\ref{1-7}) are
  the Stokes spectral invariants, i.e., $|\Omega|= |B_r|$ and $|\partial \Omega|=|\partial B_r|$. Thus $\frac{|\partial \Omega|}{|\Omega|^{(n-1)/n}} =
 \frac{|\partial B_r|}{|B_r|^{(n-1)/n}}$. Note that for any $r>0$, $\frac{|\partial B_r|}{|B_r|^{(n-1)/n}} =
 \frac{|\partial B_1|}{|B_1|^{(n-1)/n}}$.  According to the classical isoperimetric inequality (which states that for any bounded domain $\Omega\subset {\Bbb R}^n$ with smooth boundary, the following inequality holds:
 \begin{eqnarray*}\frac{|\partial \Omega|}{|\Omega|^{(n-1)/n}} \ge
 \frac{|\partial B_1|}{|B_1|^{(n-1)/n}}.\end{eqnarray*}
 Moreover, equality obtains if and only if $\Omega$ is a ball), we immediately get that $\Omega$ is a ball and $\Omega=B_r$.  $\;\; \square$

\vskip 0.28 true cm

\noindent{\bf Remark 5.2.} \   {\it  By applying the Tauberian theorem (see, for example, Theorem 15.3 of p.$\,$30 of \cite{Kor}) for the first term on the right side of (\ref{1-7}) (i.e., $\sum_{k=1}^\infty e^{-t\lambda_k}=\int_0^\infty e^{-t\tau} dN(\tau)= \frac{(n-1)|\Omega|}{(4\pi\mu t)^{n/2}}+o(t^{-n/2})$ as $t\to 0^+$), we can get the Weyl-type law   $N_S(\tau)  =\frac{(n-1) \omega_n\,|\Omega|}{(2\pi)^{n}\mu^{n/2}} \tau^{n/2} +o(\tau^{n/2})\;\; \mbox{as}\;\; \tau\to +\infty$  for the Stokes eigenvalues, i.e.,
\begin{eqnarray} \label {m4-01} \lambda_k\sim \mu \bigg(\frac{(2\pi)^n k}{(n-1)\omega_n |\Omega|}\bigg)^{2/n}\quad \, \mbox{as}\;\; k\to \infty,\end{eqnarray}
which is just the result proved by Metivier in 1978 (see \cite{Me}).}

\vskip 0.22 true cm

\noindent{\bf Remark 5.3.} \  {\it In 1986, Girault and Raviart (see \cite{GiR} or \cite{CL}) proved that all the Stokes eigenvalues coincide with all the buckling eigenvalues in the two-dimensional case (The definition of the $k$-th buckling eigenvalue $\Lambda_k$ as well as the corresponding $k$-th buckling eigenfuction $\psi_k$ for a clamped plate $\Omega\subset \mathbb{R}^n$, see (\ref{2020.12.8-1})). In \cite{Liu}, the author of this paper proved  the following asymptotic formula for $n$-dimensional buckling eigenvalues:
\begin{eqnarray}\label{17-11-27-6} \Lambda_k \sim \mu  \bigg(\frac{(2\pi)^n k}{\omega_n |\Omega|}\bigg)^{2/n} \quad \; \,\mbox{as}\;\; k\to  +\infty.\end{eqnarray}
Combining this asymptotic formula and (\ref{m4-01}), we immediately see that if $n>2$, then the Stokes spectrum can not coincide with the buckling spectrum for any $n$-dimensional smooth domain $\Omega$. In other words, Girault-Raviart's equivalent spectral conclusion holds only when $n=2$ for the Stokes eigenvalue problem and the buckling eigenvalue problem.}

 \vskip 0.58 true cm

\vskip 0.98 true cm

\centerline {\bf  Acknowledgments}

\vskip 0.39 true cm

 The author would like to sincerely thank the anonymous referees and the Editor for many important
comments and valuable suggestions on the original and revised versions of this manuscript, which motivated the
author to discuss and calculate the symbols of the needed Green operators by a new way. This research was supported by NNSF of China (11171023/A010801) and NNSF of China (11671033/A010802).

  \vskip 1.68 true cm

\vskip 0.32 true cm

\end{document}